\pdfoutput=1
\documentclass[9pt,twoside,epsf]{article}
\usepackage{stmaryrd}
\usepackage{amsmath,latexsym,amssymb,amsfonts,amsbsy}
\usepackage{bm}
\usepackage{graphicx,epstopdf}
\usepackage[caption=false]{subfig}
\usepackage{cases}
\usepackage{multirow}
\usepackage{booktabs}
\usepackage{verbatim}
\usepackage{epstopdf}

\newfont{\bb}{msbm10}
\def\Bbb#1{\mbox{\bb #1}}

\def\T{\top}

\def\diag{{\rm diag}}

\def\rank{{\rm rank}}

\usepackage{colortbl}
\usepackage{ragged2e}
\usepackage[margin=2em,labelsep=space,skip=0.5em,font=normalfont]{caption}
\usepackage{algorithm}  
\usepackage{algpseudocode} 
\floatname{algorithm}{\color{black} Algorithm}
\algrenewcommand{\algorithmiccomment}[1]{\quad{\color{red}\%\ #1}}
\numberwithin{algorithm}{section}
\makeatletter
\newenvironment{breakalgo}[2]{%
  \captionsetup{margin=0pt,justification=RaggedRight,singlelinecheck=false}%
  \par\noindent%
  \medskip%
  \rule{\linewidth}{0.8pt}%
  \vspace{-0.5\baselineskip}%
  \noindent\captionof{algorithm}{#1}\label{#2}%
  \vspace{-0.7\baselineskip}%
  \noindent\rule{\linewidth}{.4pt}%
  \vspace{-0.3\baselineskip}%
}{%
  \vspace{-.75\baselineskip}%
  \rule{\linewidth}{.4pt}%
  \medskip%
}

\makeatother

\newtheorem{method}{Method}[section]

\newtheorem{remark}{Remark}[section]
\newtheorem{theorem}{Theorem}[section]

\newtheorem{lemma}{Lemma}[section]

\newcommand{\reals}{\makebox{{\Bbb R}}}
\newcommand{\ceals}{\makebox{{\Bbb C}}}

\newcommand{\Z}{\makebox{{\Bbb Z}}}

\newcommand{\bbb}[1]{\text{\bf #1}}

\makeatletter 

\@addtoreset{equation}{section}
\makeatother 

\begin{document}
\cleardoublepage
\pagestyle{myheadings}

\bibliographystyle{plain}

\title{Diagonal and normal with Toeplitz-block splitting iteration method for space fractional coupled nonlinear Schrödinger equations with repulsive nonlinearities\thanks{Corresponding author: Xi Yang (yangxi@nuaa.edu.cn); Other authors: Fei-Yan Zhang (zhangfy@nuaa.edu.cn), Chao Chen (chenchao3@nuaa.edu.cn).}}

\author{
Fei-Yan Zhang\thanks{School of Mathematics, Nanjing University of Aeronautics and Astronautics, Nanjing 211106, China.}\ \thanks{Wuxi Taihu University, Wuxi 214064, China.},
Xi Yang\footnotemark[2]\ \thanks{Key Laboratory of Mathematical Modelling and High Performance Computing of Air Vehicles (NUAA), MIIT, Nanjing 211106, China.},
Chao Chen\footnotemark[2]
}

\maketitle

\begin{abstract}
By applying the linearly implicit conservative difference scheme proposed in [D.-L. Wang, A.-G. Xiao, W. Yang. J. Comput. Phys. 2014;272:670-681], the system of repulsive space fractional coupled nonlinear Schrödinger equations leads to a sequence of linear systems with complex symmetric and Toeplitz-plus-diagonal structure. In this paper, we propose the diagonal and normal with Toeplitz-block splitting iteration method to solve the above linear systems. The new iteration method is proved to converge unconditionally, and the optimal iteration parameter is deducted. Naturally, this new iteration method leads to a diagonal and normal with circulant-block preconditioner which can be executed efficiently by fast algorithms. In theory, we provide sharp bounds for the eigenvalues of the discrete fractional Laplacian and its circulant approximation, and further analysis indicates that the spectral distribution of the preconditioned system matrix is tight. Numerical experiments show that the new preconditioner can significantly improve the computational efficiency of the Krylov subspace iteration methods. Moreover, the behavior of the corresponding preconditioned GMRES method exhibits a linear dependence on the space mesh size, which weakens as the fractional order parameter decreases.
	
\bigskip

{\bf Key words.} circulant matrix, coupled nonlinear Schrödinger equations, fractional derivative, preconditioning, repulsive nonlinearity, Toeplitz matrix

\bigskip

{\bf MSC codes.} 65F08, 65F10, 65M06, 65M22

\end{abstract}

\section{Introduction}\label{sec-Introduction}
The system of Schrödinger equations is a basic physical model describing nonrelativistic quantum mechanical behaviors. It is well known that Feynman and Hibbs derived the standard Schrödinger equation for the evolution of microscopic particles based on the Brownian path integrals. By replacing Brownian motion to the L\'{e}vy-like process, the standard Schrödinger equation is generalized to the fractional Schrödinger equation \cite{Laskin01,Laskin02}. Since the system of fractional Schrödinger equations plays important roles in physical applications \cite{GuoXu06,MS95,ZVB95}, the related theories of fractional Schrödinger equations have been studied extensively, including the existence and uniqueness of the solution \cite{BGuo2008}, the well-posedness of one dimensional (1D) and higher dimensional problems \cite{GuoHou2010,GuoHou2013}, the construction of the ground state solution \cite{Secchi2013,LiZW2019}, etc. Due to the nonlocal nature of the fractional derivative operator, it is difficult to obtain the closed-form exact solution of a system of fractional Schrödinger equations. In order to overcome this conundrum, a number of numerical methods have been established and studied, e.g., finite element methods \cite{LM2018JCP,LM2017NUMA}, spectral methods \cite{DSW2016CMA,WY2019ANM,FF2014SIAM}, collocation methods \cite{Amore2010JMP,Bhrawy2017ANM}, and finite difference methods \cite{WDL2013JCP,WDL2014JCP,WPD2015JCP,ZRP2019SCM,ZX2014SISC}, etc.

In this paper, the following system is considered, i.e., the space fractional coupled nonlinear Schrödinger (CNLS) equations
\begin{align} \label{equ1}
	\begin{cases}
		{\imath}u_{t}-\gamma(-\Delta)^{\frac{\alpha}{2}}u+\rho(|u|^2+\beta|v|^2)u =0,\\
		{\imath}v_{t}-\gamma(-\Delta)^{\frac{\alpha}{2}}v+\rho(|v|^2+\beta|u|^2)v =0,\
	\end{cases}
	x\in\reals,\ 0 < t \le \mbox{T}
\end{align}
with the initial conditions $u(x,0)=u_{0}(x)$ and $v(x,0)=v_{0}(x)$. Here, $\imath=\sqrt{-1}$, $(-\Delta)^{\frac{\alpha}{2}}$ is the 1D fractional Laplacian with $1<\alpha \le 2$, and the parameters $\gamma>0, \beta \ge 0, \rho$ are real constants. When $\alpha = 2$, the system (\ref{equ1}) reduces to the classical CNLS equations \cite{nonlinearoptic2009}. When $\beta=0$, the system (\ref{equ1}) reduces to the decoupled nonlinear Schrödinger (DNLS) equations \cite{BGuo2008}. When $\rho=0$, the system (\ref{equ1}), without the nonlinear term, describes the free particles \cite{Laskin2002,Luchko2013}. We only consider the case of $\rho<0$ \cite{BWZ2012arXiv,Carr2000PRA,JS1999CPAM}, i.e., the system (\ref{equ1}) with repulsive nonlinearities, and leave the attractive nonlinearity case ($\rho>0$ \cite{BWZ2012arXiv,BWZ2003SINA,Saito2001PRL}) as a future work.

For solving the system (\ref{equ1}) numerically, the infinite interval $\reals$ is truncated into a bounded computational interval $\mbox{a} \le x \le \mbox{b}$, and the Dirichlet boundary conditions are adopted, i.e.,
\begin{gather*}
		u(\mbox{a},t)=u(\mbox{b},t)=0,\qquad    v(\mbox{a},t)=v(\mbox{b},t)=0, \qquad    0 \le t \le \mbox{T}.	
\end{gather*}
Then, the linearly implicit conservative difference (LICD) scheme proposed in \cite{WDL2014JCP} based on the fractional centered difference formula \cite{riesz21} can be applied to the truncated system (\ref{equ1}), which is stable and convergent, and leads to a sequence of complex linear systems. The coefficient matrices are dense, non-Hermitian, and of the form $\imath I+D-T$, where $D$ is diagonal and nonpositive, and $T$, the discrete fractional Laplacian, is Toeplitz and symmetric positive definite. Hence, these matrices can be treated as either Toeplitz-plus-diagonal matrices or complex symmetric matrices.

When we treat $\imath I+D-T$ as a Toeplitz-plus-diagonal matrix, since there is no fast direct solver available for dense Toeplitz-plus-diagonal linear systems, the preconditioned Krylov subspace iteration methods \cite{GolubBook,sad22} may be a class of reasonable options. Chan and Ng \cite{TPB12} proposed a banded preconditioner for the Toeplitz-plus-band matrix (the Toeplitz-plus-diagonal matrix is just a special case). The main drawback of this preconditioner is that the generating function of the corresponding Toeplitz matrix should be known in advance, but unfortunately, it is unavailable in general. Ng and Pan \cite{Cpd14} constructed the approximate inverse circulant-plus-diagonal (AICD) preconditioner for the Toeplitz-plus-diagonal matrix. The AICD preconditioner works effectively and efficiently if the elements of the Toeplitz matrix decay away exponentially from the main diagonal, the nonzero elements of the diagonal matrix vary mildly and smoothly, and the interpolation points with respect to the diagonals is appropriately selected. Bai et al. established a diagonal and circulant splitting (DCS) preconditioner for the Toeplitz-plus-diagonal matrix, which significantly accelerates the convergence rate of the preconditioned Krylov subspace iteration methods; see \cite{BaiLuPan2017} for the 1D case, and \cite{BaiLu2020} for the higher dimensional cases. The above preconditioners work well under a basic assumption that the system matrix is Hermitian positive definite. Nevertheless, $\imath I+D-T$ is obviously non-Hermitian, not coinciding with the basic assumption.

When we treat $\imath I+D-T$ as a complex symmetric matrix, we can consider two kinds of efficient methods for the corresponding complex symmetric linear system. The first kind is the class of the alternating-type iteration methods, including the Hermitian and skew-Hermitian splitting (HSS \cite{HSS10}) iteration method and the preconditioned HSS (PHSS \cite{BaiGP2004}) iteration method for a general non-Hermitian positive definite linear system, and the modified HSS (MHSS \cite{MHSS18}) iteration method and the preconditioned MHSS (PMHSS \cite{BaiBC2011,BaiBCW2013}) iteration method for a complex symmetric linear system who's coefficient matrix has nonnegative definite real and imaginary parts, and at least one of them is positive definite. The second kind is the C-to-R iteration method \cite{Axelsson2000}, which reforms $\imath I+D-T$ to a real block two-by-two matrix, and performs a block Gaussian elimination, finally solves a Schur-complement linear subsystem. Unfortunately, these methods can not take full use of the Toeplitz structure of $\imath I+D-T$, and thus not leading to convincing numerical behavior.

To deal with the matrix $\imath I+D-T$, Ran et al. proposed the partially inexact HSS (PIHSS \cite{Pihss08}) iteration method, and the HSS-like iteration method \cite{HSS-like09}. Numerical results in \cite{Pihss08} and \cite{HSS-like09} show that the PIHSS iteration method outperforms the HSS iteration method in terms of computing time, and the HSS-like iteration method has better behavior in terms of both iteration counts and computing time compared with the PIHSS iteration method. In \cite{RanWW2017}, Ran et al. came to a conclusion that the HSS-like preconditioner is more efficient than the HSS preconditioner and the AICD preconditioner when they are in conjunction with the Krylov subspace iteration methods. By taking the Toeplitz-plus-diagonal structure into account, Wang et al. constructed an efficient variant of the PMHSS iteration method \cite{W-PMHSS11} which naturally leads to an efficient PMHSS preconditioner. Numerical results in \cite{W-PMHSS11} show that the PMHSS preconditioner outperforms the HSS-like preconditioner. In a word, the PMHSS preconditioner suggested by Wang et al. is the most efficient one available in the literature.

In order to effectively utilize the structure of $\imath I+D-T$, we consider a real block two-by-two equivalent form and its diagonal and normal with Toeplitz-block (DNTB) splitting. This new splitting leads to the DNTB iteration method, which is a parameterized iteration method and unconditionally converges for any positive iteration parameter, and the optimal iteration parameter is deducted. The DNTB iteration method naturally admits the DNTB preconditioner, and the implementation of this preconditioner requires to solve a diagonal linear subsystem and a block normal linear subsystem with Toeplitz-block. In practice, an even more effective and efficient preconditioner can be constructed by approximating the Toeplitz-block with a circulant matrix. Then, we obtain the diagonal and normal with circulant-block (DNCB) preconditioner, which can be implemented efficiently based on the fast Fourier transform (FFT) \cite{GolubBook}. Theoretically, we provide sharp bounds for the eigenvalues of the discrete fractional Laplacian $T$ and its circulant approximation $C$, and further analysis indicates that the eigenvalues of the DNCB preconditioned system matrix are clustered around 1. Numerical experiments show that the DNCB preconditioner can significantly improve the computational efficiency of the Krylov subspace iteration methods, and outperforms a circulant improved version of the PMHSS preconditioner. In addition, the corresponding convergence behavior of the DNCB preconditioner exhibits a linear dependence on the space mesh size, which weakens as the fractional order parameter decreases.

This paper is organized as follows. In Section 2, the complex linear system with the coefficient matrix $\imath I+D-T$ is derived by applying the LICD scheme to the system of repulsive space fractional CNLS equations. In Section 3, we construct the DNTB iteration method and study its asymptotic convergence theory. The DNCB preconditioner is presented and analyzed in Section 4, and the implementation details and the computational complexities of the preconditioners involved in the experiments are described in Section 5. Numerical results are reported in Section 6, and finally, the paper is ended with some concluding remarks in Section 7.
\label{sec-introduction}

\section{Discretization of the space fractional CNLS equations }
In this section, we apply the LICD scheme \cite{WDL2014JCP} to discrete the space fractional CNLS equations (\ref{equ1}), and then derive the complex symmetric linear system on each time level.

We denote by $M$ and $N$ the prescribed positive integers, and let $\tau=\mbox{T}/N$ and $h=(\mbox{b}-\mbox{a})/(M+1)$ be the temporal step size and the spatial step size respectively. The time levels are $t_{n}=n\tau$ for $n=0,1,...,N$, and the spatial discrete points are $x_{j}=\mbox{a}+jh$ for $j=0,1,...,M+1$. The numerical solutions of (\ref{equ1}) are denoted by $u^{n}_{j}\approx u(x_{j},t_{n})$ and $v^{n}_{j}\approx v(x_{j},t_{n})$.

The fractional Laplacian $(-\Delta)^\frac{\alpha}{2}$ is equivalent to the Riesz fractional derivative in 1D case \cite{fLapReviewJCP2020}, and the latter can be discretized by the fractional center difference \cite{riesz21,MDuman2012} on the uniform spatial grids in the bounded interval $[\mbox{a},\mbox{b}]$ as
\begin{gather*}
	(-\Delta)^\frac{\alpha}{2}u(x_{j},t)=\frac{1}{h^{\alpha}}\sum^M_{k=1}c_{j-k}u_{k}+\mathcal{O}(h^2),
\end{gather*}
where the coefficients read that
$c_{k}=(-1)^{k}\Gamma(\alpha+1)/[\Gamma(\alpha/2-k+1) \Gamma(\alpha/2+k+1)]$, $\forall\ k\in\Z$,
with $\Gamma(\cdot)$ the gamma function. In addition,
$c_{k}$ satisfies the properties \cite{MDuman2012} as follows
\begin{gather} \label{propertiesOfCk}
c_{0}\ge 0,\
c_{k}=c_{-k}\le0,\ \forall\ k\ge 1,\ \mbox{and}\
\sum^{+\infty}_{k=-\infty,k\ne0}|c_{k}|=c_{0}.
\end{gather}
By adopting the LICD scheme, the following discrete space fractional CNLS equations derived from the system (\ref{equ1}) can be obtained
\begin{equation} \label{discretizedCNLS}
\begin{cases}
\imath\frac{u^{n+1}_j-u^{n-1}_j}{2\tau}-\frac{\gamma}{h^\alpha} \sum^M_{k=1}c_{j-k}\hat{u}^n_k+\rho(|u^n_j|^2+\beta|v^n_j|^2) \hat{u}^n_j=0, \\
\imath\frac{v^{n+1}_j-v^{n-1}_j}{2\tau}-\frac{\gamma}{h^\alpha} \sum^M_{k=1}c_{j-k}\hat{v}^n_k+\rho(|v^n_j|^2+\beta|u^n_j|^2) \hat{v}^n_j=0,
\end{cases}
\end{equation}
where
$\hat{u}^n_j=(u^{n+1}_j+u^{n-1}_j)/2$, $\hat{v}^n_j=(v^{n+1}_j+v^{n-1}_j)/2$, $j=1,2,\ldots,M,n=1,2,\ldots,N-1$, and
the initial and boundary conditions read that $u^0_j=u_0(x_j)$, $v^0_j=v_0(x_j)$, $u^{n}_0=u^{n}_{M+1}=0$, and $v^{n}_0=v^{n}_{M+1}=0$. In addition, the discrete space fractional CNLS equations (\ref{discretizedCNLS}) can also be written as the matrix-vector form as follows
\begin{align}\label{discretizedCNLSMaxForm}
    \begin{cases}
        A^{n+1}_u\bbb{u}^{n+1} =\bbb{b}^{n+1}_u, \\
        A^{n+1}_v\bbb{v}^{n+1} =\bbb{b}^{n+1}_v,
    \end{cases} &\forall \ n \ge 1,
\end{align}
where the coefficient matrices $A^{n+1}_u$, $A^{n+1}_v \in \ceals^{M\times M}$ read that
\begin{align*}
  \begin{cases}
    A^{n+1}_u=\imath I+ D^{n+1}_u-T, \\
    A^{n+1}_v=\imath I+ D^{n+1}_v-T.
  \end{cases}
\end{align*}
Here, $I \in \reals^{M \times M}$ represents the identity matrix, the diagonal matrices $D^{n+1}_u$, $D^{n+1}_v\in\reals^{M\times M}$ read that
\begin{align*}
  \begin{cases}
    D^{n+1}_u=\diag\{d^{n+1}_{u,1},d^{n+1}_{u,2},\ldots,d^{n+1}_{u,M}\}, \\
    D^{n+1}_v=\diag\{d^{n+1}_{v,1},d^{n+1}_{v,2},\ldots,d^{n+1}_{v,M}\},
  \end{cases}
\end{align*}
 where $d^{n+1}_{u,j}=\rho\tau(|u^n_j|^2+\beta |v^n_j|^2)$ and $d^{n+1}_{v,j}=\rho\tau(|v^n_j|^2+\beta |u^n_j|^2)$ for $j=1,2,\ldots,M$, and $T\in\reals^{M\times M}$ represents a Toeplitz matrix of the form
\begin{align}\label{equ5}
	T=\mu
	\begin{bmatrix}
		c_0 & c_{-1} & \cdots & c_{2-M} & c_{1-M} \\
		c_1 & c_0 & \cdots & c_{3-M} & c_{2-M} \\
		\vdots & \vdots & \ddots & \vdots & \vdots \\
		c_{M-2} & c_{M-3} & \cdots & c_0 & c_{-1} \\
		c_{M-1} & c_{M-2} & \cdots & c_1 & c_0 \\
	\end{bmatrix}
\end{align}
with $\mu=\frac{\gamma\tau}{h^\alpha}$. The solutions $\bbb{u}^{n+1}$, $\bbb{v}^{n+1}$ and the right-hand-sides $\bbb{b}^{n+1}_u$, $\bbb{b}^{n+1}_v$ of (\ref{discretizedCNLSMaxForm}) read that
\begin{align*}
\begin{cases}
  \bbb{u}^{n+1}={[u^{n+1}_1,...,u^{n+1}_M]}^{\T}, \\
  \bbb{v}^{n+1}={[v^{n+1}_1,...,v^{n+1}_M]}^{\T},
\end{cases}
	 \quad \mbox{and} \quad
\begin{cases}
  \bbb{b}^{n+1}_u={[b^{n+1}_{u,1},...,b^{n+1}_{u,M}]}^{\T}, \\
  \bbb{b}^{n+1}_v={[b^{n+1}_{v,1},...,b^{n+1}_{v,M}]}^{\T},
\end{cases}
\end{align*}
where
\begin{align*}
\begin{cases}
  b^{n+1}_{u,j} \ =\ \imath u^{n-1}_j+\mu \sum^M_{k=1}c_{j-k}u^{n-1}_k-d^{n+1}_{u,j}u^{n-1}_j, \\
  b^{n+1}_{v,j} \ =\ \imath v^{n-1}_j+\mu \sum^M_{k=1}c_{j-k}v^{n-1}_k-d^{n+1}_{v,j}v^{n-1}_j,
\end{cases}
	& \forall\ 1\le j\le M.
\end{align*}

Obviously, the coupled complex linear systems in (\ref{discretizedCNLSMaxForm}) share the following unified form
\begin{align} \label{equ3}
A\ \bbb{x} &= \bbb{b},
\end{align}
with a complex symmetric coefficient matrix $A=\imath I+D-T\in\ceals^{M\times M}$, since $D \in\reals^{M \times M}$ is diagonal, and $T \in \reals^{M \times M}$ is Toeplitz and symmetric positive definite. There may be two cases of $A\in\ceals^{M\times M}$ for $\rho\neq 0$. Since we already know $\gamma>0$, $\beta\ge 0$, it holds as follows.

\begin{itemize}
	\item If $\rho<0$, $D$ is negative semi-definite. Thus, together with the fact that $T$ is positive definite, it follows that $A$ is negative definite.
	\item If $\rho>0$, $D$ is positive semi-definite. Thus, together with the fact that $T$ is positive definite, it follows that $A$ is indefinite.
\end{itemize}

In this paper, the case of $\rho<0$ will be studied, and we leave the study of the case $\rho>0$ as a future work. In the latter case, a sequence of complex symmetric and non-Hermitian indefinite linear systems (\ref{discretizedCNLS}) need to be solved efficiently, which leads to a great challenge in algorithm design.

\section{The DNTB iteration method}
\quad In this section, we construct a structured iteration method for solving the complex linear system (\ref{equ3}) and provide the related convergence theory. Denote by $\bbb{x}=y+\imath z\in\ceals^{M}$ the solution, and $\bbb{b}=p+\imath q\in\ceals^{M}$ the right-hand-side, where $y$, $z$, $p$, $q\in\reals^{M}$ are real vectors. It can be easily verified that the complex linear system $(\ref{equ3})$ is equivalent to the following real non-symmetric positive definite block linear system

\begin{eqnarray}\label{positiveBlockForm}
	\mathcal{R}x\equiv
	\begin{bmatrix}
		T-D & -I\\
		I & T-D\\
	\end{bmatrix}
	\begin{bmatrix}
		z \\
		y\\
	\end{bmatrix}
	=
	\begin{bmatrix}
		-q \\
		-p \\
	\end{bmatrix}
	\equiv f,
\end{eqnarray}
where $\mathcal{R} \in \reals^{2M \times 2M}$, $x, f\in \reals^{2M}$. The system matrix $\mathcal{R}$ admits a diagonal and normal with Toeplitz-block (DNTB) splitting of the form

\begin{align}\label{rdt}
	\mathcal{R} = \mathcal{B} + \mathcal{H} 	
\end{align}
with
\begin{equation}\label{equ4}
	\mathcal{B}	=
	\begin{bmatrix}
		-D & 0 \\
		0 & -D \\	
	\end{bmatrix}	\qquad \text{and} \qquad
	\mathcal{H} =	
	\begin{bmatrix}
		T & -I \\
		I & T \\	
	\end{bmatrix}.	
\end{equation}

Motivated by the DNTB splitting and the alternating direction implicit (ADI) iteration \cite{ADI1955,ADI}, we can construct the following DNTB iteration method for solving the block linear system (\ref{positiveBlockForm}).

\begin{method}[The DNTB iteration method]
	\label{DNTBmethod}
	Let  $x^{(0)} \in\reals^{2M\times 2M}$ be an arbitrary initial guess. For $k=0,1,2,...$ until the sequence of iterates $\{x^{(k)}\}_{k\ge 0}$ converges, compute the next iterate $x^{(k+1)}$
	according to the following procedure:
	\begin{equation}\label{equ12}
		\begin{cases}
			(\omega I+\mathcal{B})x^{(k+\frac{1}{2})}=(\omega I-\mathcal{H})x^{(k)}+f, \\
			(\omega I+\mathcal{H})x^{(k+1)}=(\omega I-\mathcal{B})x^{(k+\frac{1}{2})}+f, \
		\end{cases}
	\end{equation}
	where $\omega$ is a given positive parameter.
\end{method}


The DNTB iteration can be reformulated into a one-step iteration scheme as follows
\begin{align*}
	x^{(k+1)}=\mathcal{L}_{\omega} x^{(k)}+\mathcal{F}^{-1}_{\omega} f,
\end{align*}
with
\begin{align}
	\mathcal{L}_{\omega}=(\omega I+\mathcal{H})^{-1}(\omega I-\mathcal{B})(\omega I+\mathcal{B})^{-1}(\omega I-\mathcal{H}).\label{lw}
\end{align}
This iteration scheme can also be derived from the following splitting
\begin{equation}
	\mathcal{R} \ =\ \mathcal{F}_{\omega}-\mathcal{G}_{\omega},
\end{equation}
with
\begin{align}\label{FG}
	\begin{cases}
		\mathcal{F}_{\omega} \ =\ \frac{1}{2\omega}
		(\omega I+\mathcal{B})(\omega I+\mathcal{H}),\\
		\mathcal{G}_{\omega} \ =\ \frac{1}{2\omega}
		(\omega I-\mathcal{B})(\omega I-\mathcal{H}).
	\end{cases}
\end{align}
Hence, the iteration matrix reads that $\mathcal{L}_{\omega}=\mathcal{F}_{\omega}^{-1}\mathcal{G}_{\omega}$.

The convergence property of the DNTB iteration method is summarized in the following theorem.
\begin{theorem}\label{rhoo}
	Let $\mathcal{R}\in\reals^{2M\times 2M}$ be a real non-symmetric positive definite block matrix as defined in (\ref{positiveBlockForm}). Let $\mathcal{B}, \mathcal{H}\in\reals^{2M\times 2M}$ be the matrices in the DNTB splitting (\ref{rdt}). Let $\omega$ be a positive parameter, and the iteration matrix $\mathcal{L}_{\omega}$ of the DNTB iteration scheme (\ref{equ12}) is given by (\ref{lw}).
	Then, it holds that the spectral radius of $\mathcal{L}_{\omega}$ is bounded as follows
	\begin{eqnarray}
		\rho {(\mathcal{L}_{\omega})} \le  \sigma(\omega) < 1,
	\end{eqnarray}
	with
	\begin{align}\label{sigma}
		\sigma(\omega)= \underset{ \lambda_i \in \lambda{(D)}}{\mathop{\max}}\left\lvert{\frac{\omega+\lambda_i}{\omega-\lambda_i}}\right \rvert  \underset{ \mu_i \in  \lambda{(T)}}{\mathop{\max}}\sqrt {\frac {(\omega-\mu_i)^2+1}{(\omega+\mu_i)^2+1}},
	\end{align}
	where $\lambda{(D)}$ and $\lambda{(T)}$ are the spectral sets of $D$ and $T$, respectively.
\end{theorem}
{\em Proof.}
Due to the facts that $\mathcal{B}$ is positive semi-definite, $\mathcal{H}$ is positive definite, and $\omega>0$, it follows that the matrices $\omega I+\mathcal{B}$ and $\omega I+\mathcal{H}$ are invertible.

Similarity transformation of $\mathcal{L}_{\omega}$ leads to a fact
\begin{align*}
	\widetilde{\mathcal{L}}_{\omega}=(\omega I+\mathcal{H})\mathcal{L}_{\omega}(\omega I+\mathcal{H})^{-1}
	=\mathcal{U}_{\omega}\mathcal{V}_{\omega}
\end{align*}
with
\begin{align*}
	\mathcal{U}_{\omega}=(\omega I-\mathcal{B})(\omega I+\mathcal{B})^{-1} \quad \text{and} \quad
	\mathcal{V}_{\omega}=(\omega I-\mathcal{H})(\omega I+\mathcal{H})^{-1}.
\end{align*}
Since $\mathcal{B} \ge 0$, the estimate of $\|\mathcal{U}_{\omega}\|_2$ reads that
\begin{align*}
	\left \|\mathcal{U}_{\omega} \right \|_2 = \left \|(\omega I-\mathcal{B})(\omega I+\mathcal{B})^{-1}\right \|_2
	=\underset{\lambda_i\in\lambda(D)}{\mathop{\max}}\left\lvert{\frac{\omega+\lambda_i}{\omega-\lambda_i}}\right \rvert \le 1.
\end{align*}
Due to the fact that $\mathcal{H}^\T \mathcal{H}=\mathcal{H}\mathcal{H}^\T$, it holds that
\begin{align*}
	\mathcal{V}_{\omega}^\T \mathcal{V}_{\omega} &=(\omega I+\mathcal{H}^\T)^{-1}(\omega I-\mathcal{H}^\T)(\omega I-\mathcal{H})(\omega I+\mathcal{H})^{-1}\\
	&=(\omega I-\mathcal{H}^\T)(\omega I-\mathcal{H})(\omega I+\mathcal{H}^\T)^{-1}(\omega I+\mathcal{H})^{-1}\\
	&=
	\begin{bmatrix}
		(\omega I-T)^2+I & 0 \\
		0 & (\omega I-T)^2+I \
	\end{bmatrix}
	\begin{bmatrix}
		(\omega I+T)^2+I & 0 \\
		0 & (\omega I+T)^2+I \
	\end{bmatrix}^{-1}.
\end{align*}
Since $T \in \reals^{M\times M}$ is symmetric positive definite, it reads that
\begin{align*}
	\left \|\mathcal{V}_{\omega} \right \|_2 = \underset{\mu_i \in \lambda(T)}{\mathop{\max}} \sqrt {\frac {(\omega-\mu_i)^2+1}{(\omega+\mu_i)^2+1}}  < 1.
\end{align*}
Therefore, the spectral radius of $\mathcal{L}_{\omega}$ is bounded as follows
\begin{gather*}
	\rho (\mathcal{L}_{\omega})=\rho (\widetilde{\mathcal{L}}_{\omega}) \le \left \|\mathcal{U}_{\omega} \right \|_2 \left \|\mathcal{V}_{\omega} \right \|_2 = \sigma(\omega)<1.
\end{gather*}
$\hfill\square$

\begin{remark}	
	Here are some remarks on the upper bound $\sigma(\omega)$ provided by Theorem \ref{rhoo} and the selection of the parameter $\omega$:
	\begin{itemize}
		\item [1.] The convergence rate of Method \ref{DNTBmethod} is bounded by $\sigma(\omega)$, depending on both the spectra of $D$ and $T$, but not depending on the spectra of $\mathcal{R}$ and the eigenvectors of any of the above matrices.
		\item [2.] The upper bound $\sigma(\omega)$ can be relaxed to a new bound $\widehat{\sigma}(\omega)$ as follows
		\begin{align}\label{DNTBupperBoundHat}
		\sigma(\omega) \le \underset{\lambda_{\min} \le \lambda \le \lambda_{\max}}{\mathop{\max}}\left\lvert{\frac{\omega+\lambda}{\omega-\lambda}}\right \rvert \underset{\mu_{\min} \le \mu \le \mu_{\max}}{\mathop{\max}} \sqrt {\frac {(\omega-\mu)^2+1}{(\omega+\mu)^2+1}} \equiv \hat{\sigma}(\omega),
		\end{align}
		where $\lambda_{\min}\le\lambda_{\max}\le 0$, and $0<\mu_{\min}\le\mu_{\max}$ are the lower/upper bounds of the eigenvalues of $D$ and $T$, respectively. Obviously, the bound $\widehat{\sigma}(\omega)$ has two divisors, i.e.,
		\begin{align*}
		\widehat{\sigma}(\omega)=\sigma_1 (\omega)\sigma_2 (\omega)
		\end{align*}
		with
		\begin{align*}
		\begin{cases}
		\sigma_1 (\omega)=\underset{\lambda_{\min} \le \lambda \le \lambda_{\max}}{\mathop{\max}}\left\lvert{\frac{\omega+\lambda}{\omega-\lambda}}\right \rvert,\\
		\sigma_2 (\omega)=\underset{\mu_{\min} \le \mu \le \mu_{\max}}{\mathop{\max}} \sqrt {\frac {(\omega-\mu)^2+1}{(\omega+\mu)^2+1}}.
		\end{cases}
		\end{align*}
		The divisor $\sigma_1 (\omega)$ is minimized at $\omega^\star_1$ (Corollary 2.3 in \cite{HSS10}), i.e.,
		\begin{align*}
		\omega^\star_1=\arg\underset{\omega>0}{\mathop{\min}}\left \{\underset{\lambda_{\min} \le \lambda \le \lambda_{\max}}{\mathop{\max}}\left\lvert\frac{\omega+\lambda}{\omega-\lambda}\right\rvert\right \}=\sqrt{\lambda_{\min}\lambda_{\max}}		
		\end{align*}
		and
		\begin{align*}
		\sigma_1 (\omega^\star_1)=\frac{\sqrt {-\lambda_{\min}}-\sqrt {-\lambda_{\max}}}{\sqrt {-\lambda_{\min}}+\sqrt {-\lambda_{\max}}}.
		\end{align*}
		The divisor $\sigma_2 (\omega)$ is minimized at $\omega^\star_2$ (Theorem 2.2 in \cite{On-suc17}), i.e.,
		\begin{align*}
		\omega^\star_2=\arg\underset{\omega>0}{\mathop{\min}} \left\{ \underset{\mu_{\min} \le \mu \le \mu_{\max}}{\mathop{\max}}\sqrt {\frac {(\omega-\mu)^2+1}{(\omega+\mu)^2+1}} \right\}=
		\begin{cases}
		\begin{aligned}
		\sqrt{\mu_{\min}\mu_{\max}-1}\mbox{, if $1 < \sqrt{\mu_{\min}\mu_{\max}}$},\\
		\sqrt{\mu_{\min}^2+1}\mbox{, if $1 \ge \sqrt{\mu_{\min}\mu_{\max}}$},
		\end{aligned} 	
		\end{cases}	
		\end{align*}
		and
		\begin{align*}
		\sigma_2(\omega^\star_2)=
		\begin{cases}
		\begin{aligned}
		\left(\frac{\mu_{\min}+\mu_{\max}-2\sqrt{\mu_{\min}\mu_{\max}-1}}{\mu_{\min}+\mu_{\max}+2\sqrt{\mu_{\min}\mu_{\max}-1}}\right)^{\frac{1}{2}} \mbox{, if $1 < \sqrt{\mu_{\min}\mu_{\max}}$},\\
		\left(\frac{\sqrt{\mu_{\min}^{2}+1}-\mu_{\min}}{\sqrt{\mu_{\min}^{2}+1}+\mu_{\min}}\right)^{\frac{1}{2}}\mbox{, if $1 \ge \sqrt{\mu_{\min}\mu_{\max}}$},	
		\end{aligned}
		\end{cases}
		\end{align*}	
		
		Obviously, if the values $\omega^\star_1$ and $\omega^\star_2$ are known in advance, it is reasonable to search the optimal value of $\omega$ minimizing $\widehat{\sigma}(\omega)$ between them.		
	\end{itemize}
	
\end{remark}

The optimal iteration parameter $\omega_{opt}$ which minimizes the upper bound $\widehat{\sigma}(\omega)$ is determined in the following theorem.

\begin{theorem} \label{DNTBoptimalParameter}
	Let $\omega>0$ be an arbitrary constant. Let $\widehat{\sigma}(\omega)$ be the upper bound of the convergence rate of the DNTB iteration method in (\ref{DNTBupperBoundHat}). Let $g(\lambda,\mu;\omega)$ be a function with respect to $\omega$ for given constants $\lambda \le 0$ and $\mu > 0$, i.e.,
\begin{align*}
  g(\lambda,\mu;\omega) &= \frac{\omega+\lambda}{\omega-\lambda} \sqrt{\frac{(\omega-\mu)^2+1}{(\omega+\mu)^2+1}}.
\end{align*}
After introducing the constants $\lambda_\star=\sqrt{\lambda_{\min}\lambda_{\max}}$, $\widetilde{\mu}_\star=\sqrt{\mu_{\min}\mu_{\max}-1}$, $\widehat{\mu}_\star=\sqrt{\mu_{\min}^2+1}$, and defining the functions with respect to $\omega$ as $g_1(\omega)=g(\lambda_{\max},\mu_{\max};\omega)$, $g_2(\omega)=g(\lambda_{\max},\mu_{\min};\omega)$, $g_3(\omega)=-g(\lambda_{\min},\mu_{\min};\omega)$, the optimal value of $\omega$ minimizing the upper bound $\widehat{\sigma}(\omega)$ can be determined as follows.
		\begin{itemize}
		\item When $1 < \mu_{\min} \mu_{\max}$:	
		\begin{itemize}	
			\item[\textit{(a)}] if $\lambda_\star<\widetilde{\mu}_\star$, $\widetilde{\mu}_\star<\widehat{\mu}_\star$, it holds
			$\omega_{\text{\rm opt}} = \arg \min_{\omega\in\{\omega_1,\omega_2\}} \widehat{\sigma}(\omega)$, where
			$\omega_1 =  \arg \min_{\omega\in[\lambda_\star,\widetilde{\mu}_\star]} g_1(\omega)$, and
			$\omega_2 =  \arg \min_{\omega\in[\widetilde{\mu}_\star,\widehat{\mu}_\star]} g_2(\omega)$.	
			\item[\textit{(b)}]	if $\lambda_\star<\widetilde{\mu}_\star$, $\widetilde{\mu}_\star \ge \widehat{\mu}_\star$, it holds
			$\omega_{\text{\rm opt}} = \arg \min_{\omega\in[\lambda_\star,\widetilde{\mu}_\star]} g_1(\omega)$.	
			\item[\textit{(c)}]	if $\lambda_\star \ge \widetilde{\mu}_\star$, $\widetilde{\mu}_\star \ge \widehat{\mu}_\star$, it holds
			$\omega_{\text{\rm opt}} = \arg \min_{\omega\in[\widetilde{\mu}_\star,\lambda_\star]} g_3(\omega)$.
			\item[\textit{(d)}]	if $\lambda_\star \ge \widetilde{\mu}_\star$, $\widetilde{\mu}_\star<\widehat{\mu}_\star$, there are two cases:
			\begin{itemize}
				\item[\textit{(d1)}] if $\lambda_\star < \widehat{\mu}_\star$, it holds
				$
				\omega_{\text{\rm opt}} = \arg \min_{\omega\in[\lambda_\star,\widehat{\mu}_\star]} g_2(\omega).
				$
				\item[\textit{(d2)}] if $\lambda_\star \ge \widehat{\mu}_\star$, it holds
				$
				\omega_{\text{\rm opt}} = \arg \min_{\omega\in[\widehat{\mu}_\star,\lambda_\star]} g_3(\omega).
				$	
			\end{itemize}
		\end{itemize}
		\item When $1 \ge \mu_{\min} \mu_{\max}$:
		\begin{itemize}
			\item[\textit{(e)}] if $\lambda_\star<\widehat{\mu}_\star$, it holds
			$
			\omega_{\text{\rm opt}} = \arg \min_{\omega\in[\lambda_\star,\widehat{\mu}_\star]} g_2(\omega).
			$
			\item[\textit{(f)}]	if $\lambda_\star \ge \widehat{\mu}_\star$, it holds
			$
			\omega_{\text{\rm opt}} = \arg \min_{\omega\in[\widehat{\mu}_\star,\lambda_\star]} g_3(\omega).
			$
		\end{itemize}
	\end{itemize}	
\end{theorem}
{\em Proof.}
For $\omega>0$, the monotonicity of $\left|\frac{\omega+\lambda}{\omega-\lambda}\right|$ with respect to $\lambda\le 0$ leads to
\begin{align*}
  \sigma_1(\omega) &= \max \left\{ \left| \frac{\omega+\lambda_{\min}}{\omega-\lambda_{\min}} \right|, \left| \frac{\omega+\lambda_{\max}}{\omega-\lambda_{\max}} \right| \right\},
\end{align*}
and the monotonicity of $\frac{(\omega-\mu)^2+1}{(\omega+\mu)^2+1}$ with respect to $\mu>0$ leads to
\begin{align*}
  \sigma_2(\omega) &=
  \begin{cases}
    \max \left\{ \sqrt{\frac {(\omega-\mu_{\min})^2+1}{(\omega+\mu_{\min})^2+1}}, \sqrt{\frac {(\omega-\mu_{\max})^2+1}{(\omega+\mu_{\max})^2+1}} \right\}, & \mbox{if } 1< \mu_{\min}\mu_{\max}, \\
    \sqrt{\frac {(\omega-\mu_{\min})^2+1}{(\omega+\mu_{\min})^2+1}}, & \mbox{if } 1\ge \mu_{\min}\mu_{\max}.
  \end{cases}
\end{align*}
Hence, the optimal value $\omega_{\text{\rm opt}}$ of the positive parameter $\omega$ minimizing the upper bound $\widehat{\sigma}(\omega)$ in (\ref{DNTBupperBoundHat}) can be obtained by considering the piecewise form and the monotonicity of $\widehat{\sigma}(\omega)$ and the different cases stated in Theorem \ref{DNTBoptimalParameter}.
$\hfill\square$

\begin{remark}\label{rmk:TBDNoptimalParameter}
  The exact formulae of $\omega_{\text{\rm opt}}$ can be determined by finding the zeroes of the following equation
  \begin{align*}
    \frac{\text{d}}{\text{d}\omega} g(\lambda,\mu;\omega) & =0 \quad\text{with $\omega\in[\omega_{\text{\rm L}}, \omega_{\text{\rm U}}]\subset (0,+\infty)$},
  \end{align*}
  where the interval $[\omega_{\text{\rm L}}, \omega_{\text{\rm U}}]$ represents any of the cases included in Theorem \ref{DNTBoptimalParameter}. The above equation can be equivalently reformulated into a quartic equation, i.e.,
  \begin{align}\label{quartic_equation}
    \Upsilon\ \omega^4 + \Theta\ \omega^2 + \Xi & = 0 \quad\text{with $\omega\in[\omega_{\text{\rm L}}, \omega_{\text{\rm U}}]\subset (0,+\infty)$},
  \end{align}
  where $\Upsilon = \mu - \lambda$, $\Theta = 2\lambda(\mu^2-1) - \mu(\mu^2+\lambda^2+1)$, and $\Xi = \lambda(\mu^2+1)(\lambda\mu-\mu^2-1)$.
  Denoting by $\Delta=\Theta^2-4\Upsilon\Xi$ the discriminant of (\ref{quartic_equation}), then we can determine the optimal value $\omega_{\text{\rm opt}}$ of $\omega$ minimizing $g_i(\omega)$ ($i=1$, $2$, $3$) as follows:
  \begin{itemize}
    \item When $\Delta<0$, the quartic equation (\ref{quartic_equation}) has no real positive zero. Therefore, it holds $\omega_{\text{\rm opt}} =\arg \min_{\omega\in\{\omega_{\text{\rm L}}, \omega_{\text{\rm U}}\}} g_i(\omega)$.
    \item When $\Delta\ge 0$, the following cases should be considered:
    \begin{itemize}
      \item if the quartic equation (\ref{quartic_equation}) has no positive zero in the interval $[\omega_{\text{\rm L}}, \omega_{\text{\rm U}}]$, it holds \[\omega_{\text{\rm opt}} =\arg \min_{\omega\in\{\omega_{\text{\rm L}}, \omega_{\text{\rm U}}\}} g_i(\omega);\]
      \item if the quartic equation (\ref{quartic_equation}) has one positive zero $\omega_0\in[\omega_{\text{\rm L}}, \omega_{\text{\rm U}}]$, it holds \[\omega_{\text{\rm opt}} =\arg \min_{\omega\in\{\omega_{\text{\rm L}}, \omega_0, \omega_{\text{\rm U}}\}} g_i(\omega);\]
      \item if the quartic equation (\ref{quartic_equation}) has two positive zeroes $\widehat{\omega}_0$, $\widetilde{\omega}_0\in[\omega_{\text{\rm L}}, \omega_{\text{\rm U}}]$, it holds \[\omega_{\text{\rm opt}} =\arg \min_{\omega\in\{\omega_{\text{\rm L}}, \widehat{\omega}_0, \widetilde{\omega}_0, \omega_{\text{\rm U}}\}} g_i(\omega).\]
    \end{itemize}
  \end{itemize}
  Obviously, the above idea for determining $\omega_{\text{\rm opt}}$ is based on the constants $\Upsilon$, $\Theta$, $\Xi$, $\omega_{\text{\rm L}}$, $\omega_{\text{\rm U}}$, which can be derived from the values of $\lambda_{\min}$, $\lambda_{\max}$, $\mu_{\min}$, $\mu_{\max}$. In addition, $\lambda_{\min}$ and $\lambda_{\max}$ can be easily obtained since they are eigenvalues of the diagonal matrix $D$, and $\mu_{\min}$ and $\mu_{\max}$ are the extreme eigenvalues of the discrete fractional Laplacian $T$, which can be estimated according to the eigenvalue bounds stated in Lemma \ref{lambda}.
\end{remark}

\section{Preconditioning}

The matrix $\mathcal{F}_{\omega}$ in (\ref{FG}) can serve as a preconditioner of Krylov subspace iteration methods for solving the linear system (\ref{positiveBlockForm}). We call $\mathcal{F}_{\omega}$ as the DNTB preconditioner. The main computational workload for implementing the preconditioner $\mathcal{F}_{\omega}$ is to solve the related generalized residual (GR) equation, which further requires to solve two linear subsystems, including $(\omega I+\mathcal{B})x=r$ and $(\omega I+\mathcal{H})x=r$. The former can be solved directly in $\mathcal{O}(M)$ flops since $\mathcal{B}$ is diagonal. The latter requires to solve a Schur complement linear subsystem with coefficient matrix of the form $(\omega I + T) + (\omega I + T)^{-1}\in \reals^{M\times M}$, which can not be handled by fast algorithms since $T$ is a Toeplitz matrix.

To remedy the imperfection of $\mathcal{F}_{\omega}$, we consider adopting the circulant preconditioning technique to reduce the computational costs, i.e., replacing the discrete fractional Laplacian $T$ by a circulant approximation. Then, we obtain a new preconditioner as follows
\begin{align*}\label{FA}
	\widetilde{\mathcal{F}}_{\omega}=\frac{1}{2\omega}
	(\omega I+\mathcal{B})(\omega I+\mathcal{C}),
\end{align*}
where $\mathcal{C}=\bigl[ \begin{smallmatrix}
                C & -I \\
                I & C
              \end{smallmatrix}\bigr]$,
and one can take $C\in \reals^{M\times M}$ as any circulant approximation of $T$. Obviously, $\widetilde{\mathcal{F}}_{\omega}$ is the product of a scalar $1/(2\omega)$, a diagonal matrix $\omega I+\mathcal{B}$, and a normal matrix $\omega I+\mathcal{C}$ with circulant-blocks. Thus, $\widetilde{\mathcal{F}}_{\omega}$ is called the diagonal and normal with circulant-block (DNCB) preconditioner. In the sequel, the theoretical analysis of $\widetilde{\mathcal{F}}_{\omega}$ is based on the Strang's circulant approximation \cite{Cir-pre20} of even order for the sake of simplicity. For odd order case, similar theories can also be developed. When $M$ is even, it reads that
\begin{equation}\label{equ6}
	C=\mu
	\begin{bmatrix}
		c_0 & c_{1} & \cdots & c_{M/2-1} & 0 & c_{M/2-1} & \cdots & c_{1} \\
		c_{1} & c_0 & \ddots & \ddots & c_{M/2-1} & 0 & \ddots & \vdots \\
		\vdots & \ddots & \ddots & \ddots & \ddots & \ddots & \ddots & c_{M/2-1} \\
		c_{M/2-1} & \ddots & \ddots & \ddots & \ddots & \ddots & \ddots & 0\\
		0 & \ddots & \ddots & \ddots & \ddots & \ddots & \ddots & c_{M/2-1}\\
		c_{M/2-1} & 0 & \ddots & \ddots & \ddots & \ddots  & \ddots & \vdots\\
		\vdots & \ddots & \ddots & \ddots & \ddots & \ddots & c_0 & c_{1} \\
		c_{1} & \cdots & c_{M/2-1} & 0 & c_{M/2-1} & \ddots & c_{1} & c_0 \\
	\end{bmatrix}.
\end{equation}

In order to obtain sharp eigenvalue bounds for the discrete fractional Laplacian $T$ and its Strang's circulant approximation $C$, we first prove the following lemma.
\begin{lemma}\label{lemma1}
	Let $c_j=(-1)^j\Gamma(\alpha+1)/[\Gamma(\alpha/2-j+1)\Gamma(\alpha/2+j+1)], k_0 \ge 3,$ and $1<\alpha < 2$. Then, it holds that
	\begin{gather*}
 \frac{(1-\frac{1+\alpha}{5+\alpha/2})^{5+\alpha/2} e^{\alpha+1} \Gamma(\alpha+1) \sin(\pi \alpha/2)}{\pi \alpha (k_0+1/2)^{\alpha}}<
	\sum_{j=k_0+1}^{\infty}|c_j| < \frac{\sqrt{2} e^{13/12} \Gamma(\alpha+1) \sin(\pi \alpha /2)}{\pi \alpha (k_0-1)^{\alpha}}.
	\end{gather*}
\end{lemma}
{\em Proof.} According to the property of Gamma function $\Gamma(s+1)=s\Gamma(s), \forall s \in \ceals \backslash \{0,-1,-2,...\}$, it reads that
\begin{align*}
\Gamma(\alpha / 2-j+1)=\frac{\Gamma(\alpha / 2)}{\displaystyle\prod_{l=2}^j(\alpha / 2-l+1)}
\end{align*}
and
\begin{align*}
\left \vert \displaystyle\prod_{l=2}^j(\alpha / 2-l+1) \right \vert=\displaystyle\prod_{l=2}^j(l-1-\alpha / 2)=\frac{\Gamma(j-\alpha / 2)}{\Gamma(1-\alpha / 2)},
\end{align*}
together with the fact $\Gamma(\alpha / 2)\Gamma(1-\alpha / 2)=\pi / \sin(\pi \alpha/2)$, it follows that
\begin{align*}
\left \vert c_j \right \vert = \frac{\Gamma(\alpha+1)\sin(\pi \alpha / 2)}{\pi} \frac{\Gamma(j-\alpha / 2)}{\Gamma(j+1+\alpha / 2)}.
\end{align*}
Based on the following fact of Gamma function
\begin{align*}
\Gamma(s)=\sqrt{2\pi}s^{s-1/2}e^{-s}(1+O(1 / s)), \quad for \quad s \rightarrow +\infty,
\end{align*}
where $\mathcal{O}(1 / s)$ is positive and monotonically decreasing for $s > 0$, and $1+\mathcal{O}(1 / s)<e^{13/12}$ for $s \ge 2$ \cite{robbin23}, it reads that
\begin{align*}
\frac{\Gamma(j-\alpha / 2)}{\Gamma(j+1+\alpha / 2)}=(1-\frac{1+\alpha}{j+1+\alpha / 2})^{j+1+\alpha / 2}(1-\frac{1+\alpha}{j+1+\alpha / 2})^{-1/2}(j-\alpha / 2)^{-(\alpha+1)}e^{\alpha+1}\frac{1+O(\frac{1}{j-\alpha / 2})}{1+O(\frac{1}{j+1+\alpha / 2})}.
\end{align*}
Due to the facts $(1-\frac{\alpha+1}{s})^s$ monotonically increasing and satisfying $(1-\frac{\alpha+1}{s})^s \le e^{-(\alpha+1)}$ for $s \ge \alpha+1$, $1 < (1-\frac{1+\alpha}{j+1+\alpha/2})^{-1 / 2} \le \sqrt{2}$ for $j \ge 4$, and $(j-1 / 2)^{-(\alpha+1)} < (j-\alpha / 2)^{-(\alpha+1)} \le (j-1)^{-(\alpha+1)}$ for $j > 1$, it holds that
\begin{align*}
(1-\frac{1+\alpha}{5+\alpha/2})^{5+\alpha / 2} e^{\alpha+1} (j-1 / 2)^{-(\alpha+1)} < \frac{\Gamma(j-\alpha / 2)}{\Gamma(j+1+\alpha / 2)} < \sqrt{2} e^{13/12}  (j-1)^{-(\alpha+1)}.
\end{align*}
Then, the estimate of $|c_j|$ for $j \ge 4$ reads that
\begin{align*}
\frac{(1-\frac{1+\alpha}{5+\alpha/2})^{5+\alpha / 2} e^{\alpha+1} \Gamma(\alpha+1)\sin(\pi \alpha / 2)}{\pi (j-1 / 2)^{\alpha+1}}\le |c_j| \le \frac{\sqrt{2} e^{13/12}\Gamma(\alpha+1)\sin(\pi \alpha / 2)}{\pi (j-1)^{\alpha+1}}.
\end{align*}
The above estimate together with the facts $k_0 \ge 3$, $\displaystyle \sum_{j=k_0+1}^{+\infty}(j-1 / 2)^{-(\alpha+1)}  > \frac{1}{\alpha(k_0+1 / 2)^\alpha}$, and $\displaystyle \sum_{j=k_0+1}^{+\infty}(j-1)^{-(\alpha+1)} <\frac{1}{\alpha(k_0-1)^\alpha}$ result in the estimate of $\sum_{j=k_0+1}^{\infty}|c_j|$.

$\hfill\square$

Based on Lemma \ref{lemma1}, we can derive the eigenvalue bounds of $T$ and $C$.

\begin{lemma}\label{lambda}
	Let $T$ be the Toeplitz matrix in (\ref{equ5}), and let $C$ be the Strang's circulant approximation of $T$ in (\ref{equ6}). Denote by	
	\begin{align*}
	\theta = \frac{(1-\frac{1+\alpha}{5+\alpha/2})^{5+\alpha/2} e^{1+\alpha} \Gamma(\alpha+1) \sin(\pi \alpha /2)}{ \pi \alpha},
	\end{align*}
	let $M$ be even, then it holds that
	\begin{align*}
	\frac{2\gamma \tau \theta }{(\text{\rm b}-\text{\rm a})^{\alpha}}   <  \lambda_T <  \frac{2 \gamma \tau}{h^{\alpha}} \left[\frac{\Gamma(\alpha+1)}{\Gamma(\alpha/2+1)^2}- \frac{\theta h^{\alpha}}{(\text{\rm b}-\text{\rm a})^{\alpha}}\right], \quad \text{for $M
		\ge 4$},
	\end{align*}
	and
	\begin{align*}
	\frac{2^{\alpha+1}\gamma \tau \theta }{(\text{\rm b}-\text{\rm a})^{\alpha} } < \lambda_C <  \frac{2 \gamma \tau}{h^{\alpha}} \left[\frac{\Gamma(\alpha+1)}{\Gamma(\alpha/2+1)^2}- \frac{2^{\alpha} \theta h^{\alpha}}{(\text{\rm b}-\text{\rm a})^{\alpha}}\right], \quad \text{for $M
		\ge 8$},
	\end{align*}
   where $\lambda_T$ and $\lambda_C$ are eigenvalues of $T$	and $C$, respectively.
\end{lemma}
{\em Proof.} According to the Gerschgorin disk theorem \cite{GolubBook}, the properties of $c_k$, and Lemma \ref{lemma1}, it holds that
\begin{align*}
\mu(c_0-2\sum_{k=1}^{M-1}|c_k|) \le &\lambda_T \le \mu(c_0+2\sum_{k=1}^{M-1}|c_k|)\\
2\mu\sum_{k=M}^{+\infty}|c_k|   \le &\lambda_T \le 2 \mu(c_0-\sum_{k=M}^{+\infty}|c_k|)\\
\frac{2\gamma \tau \theta }{(\mbox{b}-\mbox{a})^{\alpha}}   <  &\lambda_T < 2 \frac{\gamma \tau}{h^{\alpha}} \left[\frac{\Gamma(\alpha+1)}{\Gamma(\alpha/2+1)^2}- \frac{\theta h^{\alpha}}{(\mbox{b}-\mbox{a})^{\alpha}}\right],  \quad \text{for $M \ge 4$},
\end{align*}
and
\begin{align*}
\mu(c_0-2\sum_{k=1}^{M/2-1}|c_k|) \le &\lambda_C \le \mu(c_0+2\sum_{k=1}^{M/2-1}|c_k|)\\
2\mu\sum_{k=M/2}^{+\infty}|c_k|   \le &\lambda_C \le 2 \mu(c_0-\sum_{k=M/2}^{+\infty}|c_k|)\\
\frac{2^{\alpha+1}\gamma \tau \theta }{(\mbox{b}-\mbox{a})^{\alpha} } < &\lambda_C < 2 \frac{\gamma \tau}{h^{\alpha}} \left[\frac{\Gamma(\alpha+1)}{\Gamma(\alpha/2+1)^2}- \frac{2^{\alpha} \theta h^{\alpha}}{(\mbox{b}-\mbox{a})^{\alpha}}\right],  \quad \text{for $M \ge 8$}.
\end{align*}
$\hfill\square$

\begin{remark}
	According to Lemma \ref{lambda}, the upper bounds of the spectral condition numbers of $T$ and $C$ read that
	\begin{align*}
	 \kappa(T) \le \frac{(M+1)^{\alpha}\Gamma(\alpha+1)}{\theta\Gamma(\alpha/2+1)^2} - 1 \quad \text{and} \quad \kappa(C) \le \frac{(M+1)^{\alpha}\Gamma(\alpha+1)}{\theta 2^{\alpha} \Gamma(\alpha/2+1)^2} - 1.
    \end{align*}
   Obviously, the magnitudes of $\kappa(T)$ and $\kappa(C)$  are of the same order, say $\mathcal{O}((M+1)^\alpha)$.
\end{remark}

Since the DNCB preconditioned system matrix $\widetilde{\mathcal{F}}^{-1}_{\omega}\mathcal{R}$ can be factorized as follows
\begin{align}\label{two}
	\widetilde{\mathcal{F}}^{-1}_{\omega}\mathcal{R} &=\underbrace{\widetilde{\mathcal{F}}^{-1}_{\omega} \mathcal{F}_{\omega}}  \underbrace{\mathcal{F}_{\omega}^{-1} \mathcal{R}},
\end{align}
the property of $\widetilde{\mathcal{F}}^{-1}_{\omega}\mathcal{R}$ can be derived from the properties of both $\widetilde{\mathcal{F}}^{-1}_{\omega} \mathcal{F}_{\omega}$ and the DNTB preconditioned system matrix $\mathcal{F}_{\omega}^{-1} \mathcal{R}$. Firstly, the results of $\mathcal{F}^{-1}_{\omega}\mathcal{R}$ is stated in the following theorem.

\begin{theorem}\label{FR}
Let $\mathcal{R} \in \mathbb{R}^{2M \times 2M}$ be the real non-symmetric positive definite matrix in (\ref{positiveBlockForm}), $\omega > 0$ be a prescribed parameter, and $\sigma(\omega)$ be defined by (\ref{sigma}). Then, the eigenvalues of the DNTB preconditioned system matrix $\mathcal{F}_{\omega}^{-1}\mathcal{R}$ are located in a circle of radius $\sigma(\omega)$ centered at $1$.
\end{theorem}
{\em Proof.}
The result can be derived easily from the relation $\mathcal{L}_{\omega} =  I-\mathcal{F}_{\omega}^{-1}\mathcal{R}$ and Theorem \ref{rhoo}.
$\hfill\square$

Secondly, by taking advantages of Lemmas \ref{lemma1} and \ref{lambda}, we can summarize the property of $\widetilde{\mathcal{F}}^{-1}_{\omega}\mathcal{F}_{\omega}$ as follows.

\begin{theorem}\label{ffief}
	Let $1<\alpha<2$, and $M \ge 8$ be even. Defining the constants
	\begin{align*}
		\theta = \frac{(1-\frac{1+\alpha}{5+\alpha/2})^{5+\alpha/2} e^{1+\alpha} \Gamma(\alpha+1) \sin(\pi \alpha/2)}{\pi \alpha }
		\qquad and \qquad
		\theta_{0} = \frac{\sqrt{2} e^{13/12} \Gamma(1+\alpha) \sin(\pi \alpha/2)}{ \pi \alpha},
	\end{align*}
	let $\epsilon$ be a small positive constant satisfying $2^{\alpha} \mu \theta_0 / (M-2)^{\alpha}<\epsilon \le \mu \theta_{0}$, and $k_0 = \lceil (\mu \theta_{0} / \epsilon)^{1/\alpha} \rceil +1 $ where $\lceil \centerdot \rceil$ rounds a real number towards positive infinity. Then, there exist two matrices $\mathcal{E}_{\omega c} \in \mathbb{R}^{2M \times 2M}$ and $\mathcal{F}_{\omega c} \in \mathbb{R}^{2M \times 2M}$, satisfying $\rank(\mathcal{E}_{\omega c})=4k_0$ and
	\begin{gather*}
		\| \mathcal{E}_{\omega c} \|_2 \le \frac{M^{1/2} \mu}{\sqrt{1+\left[\omega+\frac{2^{\alpha+1} \gamma \tau \theta}{(\text{\rm b}-\text{\rm a})^{\alpha}}\right]^2}}\left[\frac{c_0}{2}-\frac{\theta}{(M-1/2)^{\alpha}}\right], \quad
		\| \mathcal{F}_{\omega c} \|_2 \le \frac{M^{1/2} \epsilon}{\sqrt{1+\left[\omega+\frac{2^{\alpha+1} \gamma \tau \theta}{(\text{\rm b}-\text{\rm a})^{\alpha}}\right]^2}}.
	\end{gather*}
	Such that
	\begin{gather*}
		\widetilde{\mathcal{F}}_{\omega}^{-1}\mathcal{F}_{\omega} = I+\mathcal{E}_{\omega c}+\mathcal{F}_{\omega c}.
	\end{gather*}
\end{theorem}
{\em Proof.}  According to the structure of $T$ and $C$, and the fact that $M$ is even, it reads that
\begin{align*}
T-C=\widehat{E}+\widehat{F}
\end{align*}
where
\begin{gather*}
\widehat{E}=\mu
\begin{bmatrix}
0 & 0 & \widehat{E}_{13} \\
0 & 0 & 0\\
\widehat{E}^{\T}_{13} & 0 & 0 \\
\end{bmatrix},\quad
\widehat{F}=\mu
\begin{bmatrix}
0 & \widehat{F}_{12} & 0 \\
\widehat{F}^{\T}_{12} & 0 & 0 \\
0 & 0 & 0\\
\end{bmatrix}.
\end{gather*}
and
\begin{gather*}
	\widehat{F}_{12}=
	\begin{bmatrix}
		c_{M/2} & c_{M/2+1}-c_{M/2-1} & \cdots &  c_{M-k_0-1}-c_{k_0+1}\\
		0 & \ddots & \ddots &  \vdots\\
		\vdots & \ddots & \ddots & c_{M/2+1}-c_{M/2-1} \\
		\vdots & \ddots & \ddots  & c_{M/2}\\
		0 & \cdots & \cdots & 0\\
		\vdots &  &  & \vdots\\
		0 & \cdots & \cdots  & 0\\
	\end{bmatrix}
	\in \mathbb{R}^{(M/2)\times\left(M/2-k_0\right)},
\end{gather*}
\begin{gather*}
	\widehat{E}_{13}=
	\begin{bmatrix}
		c_{M-k_0}-c_{k_0} & \cdots & \cdots  & c_{M-1}-c_1\\
		\vdots & \ddots & \ddots &  \vdots\\
		c_{M/2+1}-c_{M/2-1} & \ddots & \ddots & \vdots\\
		c_{M/2} & \ddots & \ddots  & \vdots\\
		0 & \ddots & \ddots &  c_{M-k_0}-c_{k_0}\\
		\vdots & \ddots & \ddots &  c_{M/2+1}-c_{M/2-1}\\
		0 & \cdots & 0 &  c_{M/2}\\
	\end{bmatrix}
	\in \mathbb{R}^{(M/2)\times k_0}.
\end{gather*}
with $2\le k_0 \le M/2$.

Obviously, the form of $\widehat{E}$ guarantees $\rank(\widehat{E})=2k_0$. In addition, the structures of $\widehat{E}$ and $\widehat{F}$ together with Lemma \ref{lemma1} lead to the following $\ell_{\infty}$-norm estimates
\begin{align*}
	\| \widehat{E} \|_\infty &=\mu \max\left\{{ \| \widehat{E}_{13}  \|_\infty, \| \widehat{E}^\T_{13}  \|_\infty}\right\}=\mu \| \widehat{E}_{13}  \|_\infty\\
	& \le \mu \sum^{M-1}_{j=1}|c_j|=\mu\left(\frac{c_0}{2}-\sum_{j=M}^{\infty}|c_{j}|\right)\\
	& < \mu\left[\frac{c_0}{2}-\frac{\theta}{({M-1/2)}^{\alpha}}\right],
\end{align*}
\begin{align*}
	\| \widehat{F} \|_\infty &=\mu \max\left\{{ \| \widehat{F}_{12}  \|_\infty, \| \widehat{F}^\T_{12}  \|_\infty}\right\}=\mu \| \widehat{F}_{12}  \|_\infty\\
	& \le \mu \sum^{M-k_0-1}_{j=k_0+1}|c_j|<\mu \sum^\infty_{j=k_0+1}|c_j|\\
	& < \frac{\mu \theta_{0}}{(k_0-1)^{\alpha}} < \epsilon.
\end{align*}
Straightforward computations lead to
\begin{align*}
	\widetilde{\mathcal{F}}^{-1}_{\omega} \mathcal{F}_{\omega}-I & =
	\begin{bmatrix}
		\omega I+C & -I\\
		I & \omega I+C
	\end{bmatrix}^{-1}
	\begin{bmatrix}
		T-C & 0\\
		0 & T-C
	\end{bmatrix}\\
	& = \mathcal{E}_{\omega c}+\mathcal{F}_{\omega c},
\end{align*}
where
\begin{gather*}
	\mathcal{E}_{\omega c}=
	\begin{bmatrix}
		\omega I+C & -I\\
		I & \omega I+C
	\end{bmatrix}^{-1}
	\begin{bmatrix}
		\widehat{E} & 0\\
		0 & \widehat{E}
	\end{bmatrix}
	\quad \text{and} \quad
	\mathcal{F}_{\omega c}=
	\begin{bmatrix}
		\omega I+C & -I\\
		I & \omega I+C
	\end{bmatrix}^{-1}
	\begin{bmatrix}
		\widehat{F} & 0\\
		0 & \widehat{F}
	\end{bmatrix}.
\end{gather*}
Due to the equivalence of vector norms and the $\ell_{\infty}$-norm estimates of $\widehat{E}$ and $\widehat{F}$, we can obtain the $\ell_2$-norm estimates of $\mathcal{E}_{\omega c}$ and $\mathcal{F}_{\omega c}$ as follows
\begin{align*}
	\left \|\mathcal{E}_{\omega c}\right \|_2 &\le
	\left \|
	\begin{bmatrix}
		\omega I+C & -I\\
		I & \omega I+C
	\end{bmatrix}^{-1}
	\right \|_2
	\left \|\widehat{E}\right \|_2 \\
	&\le
	\left \|
	\begin{bmatrix}
		\omega I+C & -I\\
		I & \omega I+C
	\end{bmatrix}^{-1}
	\right \|_2
	 M^{1/2} \| \widehat{E} \|_{\infty}\\
	&\le
	\frac{M^{1/2} \mu}{\sqrt{1+\left[\omega+\frac{2^{\alpha+1} \gamma \tau \theta }{(\text{\rm b}-\text{\rm a})^{\alpha}}\right]^2}}\left[\frac{c_0}{2}-\frac{\theta}{(M-1/2)^{\alpha}}\right],
\end{align*}

\begin{align*}
	\left \|\mathcal{F}_{\omega c}\right \|_2 &\le
	\left \|
	\begin{bmatrix}
		\omega I+C & -I\\
		I & \omega I+C
	\end{bmatrix}^{-1}
	\right \|_2
	 \left \|\widehat{F}\right \|_2 \\
	&\le
	\left \|
	\begin{bmatrix}
		\omega I+C & -I\\
		I & \omega I+C
	\end{bmatrix}^{-1}
	\right \|_2
	 M^{1/2} \| \widehat{F} \|_{\infty}\\
	&<
	\frac{M^{1/2} \epsilon}{\sqrt{1+\left[\omega+\frac{2^{\alpha+1} \gamma \tau \theta }{(\text{\rm b}-\text{\rm a})^{\alpha}}\right]^2}}.
\end{align*}
$\hfill\square$

\begin{remark}
According to Theorem \ref{ffief}, we know that $\widetilde{\mathcal{F}}^{-1}_{\omega} \mathcal{F}_{\omega}-I=\mathcal{E}_{\omega c}+\mathcal{F}_{\omega c}$, where $\mathcal{E}_{\omega c}$ is a low rank matrix, and $\mathcal{F}_{\omega c}$ is a small norm matrix. On one hand, the matrix $I + \mathcal{F}_{\omega c}$ is a small perturbation of the identity matrix $I$. Specifically, the Bauer-Fike theorem \cite{BauerFike1960} leads to a fact that the eigenvalues of $I + \mathcal{F}_{\omega c}$ are clustered in a small disk centered at 1. On the other hand, since the matrix $\mathcal{E}_{\omega c}$ has bounded $\ell_2$-norm and low rank, the matrix $\widetilde{\mathcal{F}}^{-1}_{\omega} \mathcal{F}_{\omega}$ can be considered as a low rank correction of $I + \mathcal{F}_{\omega c}$. Then, one can expect that only a small number of the eigenvalues of $\widetilde{\mathcal{F}}^{-1}_{\omega} \mathcal{F}_{\omega}$ drift relatively further away from those of $I + \mathcal{F}_{\omega c}$. Hence, most of the eigenvalues of $\widetilde{\mathcal{F}}^{-1}_{\omega} \mathcal{F}_{\omega}$ are clustered around $1$.
\end{remark}

Finally, together with the relation (\ref{two}), and Theorems \ref{FR} and \ref{ffief}, we turn to consider the eigenvalue clustering property of  $\widetilde{\mathcal{F}}^{-1}_{\omega}\mathcal{R}$.
\begin{theorem}\label{eigDistrDNCBPrecSysMatrix}
	Let $1<\alpha < 2$, and $M \ge 8$ be even, and	
	\begin{align*}
		\xi=\frac{2M^{1/2}\sqrt{1+\left\{\omega+2\left[\frac{\Gamma(\alpha+1)\mu}{\Gamma(\alpha/2+1)^2}-\frac{\theta \gamma \tau}{ (\text{\rm b}-\text{\rm a})^{\alpha}}\right]\right\}^2}}{\sqrt{1+\left[\omega+\frac{2^{\alpha+1} \gamma \tau \theta}{(\text{\rm b}-\text{\rm a})^{\alpha}}\right]^2}\sqrt{1+\left[\omega +\frac{2 \gamma \tau \theta }{(\text{\rm b}-\text{\rm a})^{\alpha}}\right]^2}}.
	\end{align*}
	Let $\epsilon$ be a small positive constant satisfying $2^{\alpha} \mu \theta_0 / (M-2)^{\alpha}\le \epsilon \le \mu \theta_{0}$, and $k_0 = \lceil ( \mu \theta_{0} / \epsilon)^{1 /\alpha} \rceil +1 $. Then, there exist two matrices $\mathcal{M}_{\omega} \in \mathbb{R}^{2M \times 2M}$ and $\mathcal{N}_{\omega} \in \mathbb{R}^{2M \times 2M}$ satisfying $\rank(\mathcal{M}_{\omega})=4k_0$, $\| \mathcal{M}_{\omega} \|_2 \le \mu \xi \bigl[\frac{c_0}{2}-\frac{\theta}{(M-1/2)^{\alpha}}\bigr]$, and $\| \mathcal{N}_{\omega} \|_2 \le \xi \epsilon$, such that
	\begin{gather*}
		\widetilde{\mathcal{F}}^{-1}_{\omega}\mathcal{R} = \mathcal{F}_{\omega}^{-1} \mathcal{R}+\mathcal{M}_{\omega}+\mathcal{N}_{\omega}.
	\end{gather*}		
\end{theorem}
{\em Proof.}
It can be easily verified that
\begin{align*}
	\mathcal{F}_{\omega}^{-1}\mathcal{R} &= I-\mathcal{F}_{\omega}^{-1}\mathcal{G}_{\omega}\\
	& = (\omega I+\mathcal{H})^{-1}(I-\widetilde{\mathcal{L}}_{\omega})(\omega I+\mathcal{H})
\end{align*}
with $\widetilde{\mathcal{L}}_{\omega}=\mathcal{U}_{\omega} \mathcal{V}_{\omega}$.

According to the proof of Theorem \ref{rhoo}, we know that $\| \widetilde{\mathcal{L}}_{\omega} \|_2 < 1 $, it then holds $\|I- \widetilde{\mathcal{L}}_{\omega}\|_2 \le \| I \|_2 + \|\widetilde{\mathcal{L}}_{\omega}\|_2 < 2$. The lower/upper bound estimates of the eigenvalues of $T$ provided in Lemma \ref{lambda} lead to the following fact
\begin{align*}
	\|(\omega I+\mathcal{H})^{-1}\|_2 \| \omega I+\mathcal{H} \|_2
	& \le \frac{\displaystyle \max_{\lambda_{i} \in \lambda(T)} \sqrt{1+(\omega+\lambda_i)^2}}{\displaystyle \min_{\lambda_{i} \in \lambda(T)} \sqrt{1+(\omega+\lambda_i)^2}}\\
	& \le \frac{\sqrt{1+\left\{\omega+2\left[\frac{\Gamma(\alpha+1)\mu}{\Gamma(\alpha / 2+1)^2}-\frac{\theta \gamma \tau }{(\text{\rm b}-\text{\rm a})^{\alpha}}\right]\right\}^2}}{\sqrt{{1+\left[\omega+\frac{2\gamma \tau \theta }{(\text{\rm b}-\text{\rm a})^{\alpha}} \right]^2}}}.
\end{align*}
Then, the $\ell_2$-norm of $\mathcal{F}^{-1}_{\omega} \mathcal{R}$ is bounded as follows
\begin{align*}
	\left \| \mathcal{F}^{-1}_{\omega} \mathcal{R} \right \|_2 & \le
	\left \|(\omega I+\mathcal{H})^{-1}\right \|_2 \left \| \omega I+\mathcal{H} \right \|_2 \left \| I-\widetilde{\mathcal{L}}_{\omega} \right \|_2\\
	& \le
	\frac{2 \sqrt{1+\left\{\omega+2\left[\frac{\Gamma(\alpha+1)\mu}{\Gamma(\alpha/2+1)^2}-\frac{\theta \gamma \tau }{(\text{\rm b}-\text{\rm a})^{\alpha}}\right]\right\}^2}}{\sqrt{{1+\left[\omega+\frac{2\gamma \tau \theta }{(\text{\rm b}-\text{\rm a})^{\alpha}} \right]^2}}}.
\end{align*}
According to the relation (\ref{two}) and Theorem \ref{ffief}, we have
\begin{align*}
	\widetilde{\mathcal{F}}_{\omega}^{-1}\mathcal{R} & =(I+\mathcal{E}_{\omega c}+\mathcal{F}_{\omega c})  \mathcal{F}_{\omega}^{-1} \mathcal{R}\\
	& = \mathcal{F}_{\omega}^{-1} \mathcal{R} + \mathcal{M}_{\omega}+ \mathcal{N}_{\omega}
\end{align*}
with $ \mathcal{M}_{\omega}=\mathcal{E}_{\omega c} \mathcal{F}_{\omega}^{-1} \mathcal{R}$ and
$\mathcal{N}_{\omega} = \mathcal{F}_{\omega c} \mathcal{F}_{\omega}^{-1} \mathcal{R}$. Moreover, it holds that $\rank(\mathcal{M}_{\omega})=4k_0$, and the matrices $\mathcal{M}_{\omega}$ and $\mathcal{N}_{\omega}$ satisfy
\begin{align*}
	\left\| \mathcal{M}_{\omega} \right\|_2 & \le \left\| \mathcal{E}_{\omega c} \right\|_2 \left\| \mathcal{F}_{\omega}^{-1} \mathcal{R} \right\|_2 \\
	& \le \mu \xi \left[\frac{c_0}{2}-\frac{\theta}{(M-1/2)^{\alpha}}\right]
\end{align*}
and
\begin{align*}
	\left\| \mathcal{N}_{\omega} \right\|_2 & \le \left\| \mathcal{F}_{\omega c} \right\|_2 \left\| \mathcal{F}_{\omega}^{-1} \mathcal{R} \right\|_2 \\
	& \le \xi \epsilon.
\end{align*}
$\hfill\square$

\begin{remark}
According Theorem \ref{FR}, we know that the eigenvalues of $\mathcal{F}_{\omega}^{-1}\mathcal{R}$ are located in a disk of radius $\sigma(\omega)$ centered at 1. Theorem \ref{eigDistrDNCBPrecSysMatrix} shows that, on one hand, $\mathcal{N}_{\omega}$ is a small norm matrix, thus, if $\mathcal{F}_{\omega}^{-1}\mathcal{R}$ is diagonalizable, the eigenvalues of $\mathcal{F}_{\omega}^{-1}\mathcal{R}+\mathcal{N}_{\omega}$ can be considered as small perturbations of the eigenvalues of $\mathcal{F}_{\omega}^{-1}\mathcal{R}$. On the other hand, $\mathcal{M}_{\omega}$ is a low rank matrix, we may conclude that most of the eigenvalues of $\widetilde{\mathcal{F}}_{\omega}^{-1}\mathcal{R}$ are distributed around the eigenvalues of $\mathcal{F}_{\omega}^{-1}\mathcal{R}+\mathcal{N}_{\omega}$. In summary, it may be true that most of the eigenvalues of $\widetilde{\mathcal{F}}_{\omega}^{-1}\mathcal{R}$ are clustered around those of $\mathcal{F}_{\omega}^{-1}\mathcal{R}$. Therefore, the eigenvalues of $\widetilde{\mathcal{F}}_{\omega}^{-1}\mathcal{R}$ are relatively clustered.
\end{remark}

\section{Implementation and complexity}
As stated in Section \ref{sec-Introduction}, it is reasonable to adopt the preconditioned Krylov subspace iteration methods for solving the discrete space fractional CNLS equations since they are efficient and competitive. Furthermore, a well-performed preconditioner is the efficient variant of the PMHSS preconditioner proposed in \cite{W-PMHSS11}, which is specifically designed for the block equivalent form (\ref{positiveBlockForm}) of the discrete space fractional CNLS equations. Therefore, the PMHSS-type and DNTB-type preconditioners are supposed to be tested numerically in Section \ref{sec-numerical-experiments}, and the implementation details and computational complexities of these preconditioners are discussed here.

The efficient variant of the PMHSS preconditioner provided by \cite{W-PMHSS11} reads that
\begin{align*} \label{preconditioner-PMHSS}
\mathcal{F}_{\mbox{\tiny PMHSS}} & =
\begin{bmatrix}
I & I \\
-I & I
\end{bmatrix}^{-1}
\begin{bmatrix}
\omega I + T & 0 \\
0 & \omega I + T
\end{bmatrix}
\begin{bmatrix}
\widehat{D} & 0 \\
0 & \widehat{D}
\end{bmatrix},
\end{align*}
where $\widehat{D}=(\omega I + D)^{-1}[(\omega+1)I + D]$ is diagonal, and the parameter $\omega$ satisfies $\omega >\|D\|_{\infty}$ such that $\widehat{D}$ is well defined. In actual implementation, the matrix $\widehat{D}$ can be calculated and stored in advance. The implementation of $\mathcal{F}_{\mbox{\tiny PMHSS}}$ requires to solve two linear subsystems with dense Toeplitz matrix $\omega I + T$. Obviously, direct methods will be very expensive in flops and storage, and iteration methods will lead to inner-outer schemes. Both the above ideas will affect the computational efficiency for solving the block linear system (\ref{positiveBlockForm}). Therefore, the implementation details of $\mathcal{F}_{\mbox{\tiny PMHSS}}$ are not included.

In order to improve the implementing cost of $\mathcal{F}_{\mbox{\tiny PMHSS}}$ in preconditioned Krylov subspace iteration methods, a circulant approximation of the discrete fractional Laplacian $T$ can be applied to obtain a more efficient preconditioner as follows, i.e., the circulant improved PMHSS (CPMHSS) preconditioner

\begin{equation}
\begin{aligned}
\mathcal{F}_{\mbox{\tiny CPMHSS}} &   =
\begin{bmatrix}
I & I \\
-I & I
\end{bmatrix}^{-1}
\begin{bmatrix}
\omega I + C & 0 \\
0 & \omega I + C
\end{bmatrix}
\begin{bmatrix}
\widehat{D} & 0 \\
0 & \widehat{D}
\end{bmatrix}  \\ \label{preconditioner-CPMHSS}
& =
\begin{bmatrix}
I & I \\
-I & I
\end{bmatrix}^{-1}
\begin{bmatrix}
F & 0 \\
0 & F
\end{bmatrix}^{-1}
\begin{bmatrix}
\widehat{\Lambda} & 0 \\
0 & \widehat{\Lambda}
\end{bmatrix}
\begin{bmatrix}
F & 0 \\
0 & F
\end{bmatrix}
\begin{bmatrix}
\widehat{D} & 0 \\
0 & \widehat{D}
\end{bmatrix},
\end{aligned}
\end{equation}
where $F \in \mathbb{C}^{M \times M}$ is the discrete Fourier transform (DFT), $\Lambda =\diag(Fc) \in \mathbb{C}^{M \times M} $ with $c \in \mathbb{C}^{M }$ being the first column of $C$, and $\widehat{\Lambda}=\omega I+\Lambda$ can be calculated in advance. By taking advantage of the algorithms of fast Fourier transform (FFT) and inverse FFT (IFFT), the action of $F$ or $F^{-1}$ on a vector can be accomplished in $\mathcal{O}(M\log M)$ flops. With notations $x=(x_1^{\T},x_2^{\T})^{\T}\in\reals^{2M}$ with $x_1$, $x_2\in\reals^M$, and $r=(r_1^{\T},r_2^{\T})^{\T}\in\reals^{2M}$ with $r_1$, $r_2\in\reals^M$, we can list the preconditioning procedure with respect to $\mathcal{F}_{\mbox{\tiny CPMHSS}}$ in Algorithm \ref{alg-CPMHSS}. Obviously, this algorithm requires six $M$-vector operations (including diagonal linear subsystem solves and vector additions, which can be accomplished in $\mathcal{O}(M)$ flops), and four $M$-vector FFT/IFFT operations ($\mathcal{O}(M\log M)$ flops in total).

\begin{breakalgo}{Solve the GR equation $\mathcal{F}_{\mbox{\tiny CPMHSS}}\,x=r$}{alg-CPMHSS}
	\begin{algorithmic}[1]
		\State $x_1 = r_1 + r_2$ and $x_2 = -r_1 + r_2$; \% \texttt{ Two $M$-vector operations}
		\State $x_j=\text{FFT}(x_j)$, $j=1$, $2$; \%  \texttt{ Two $M$-vector FFTs}
		\State Solve $\widehat{\Lambda}x_j = x_j$, $j=1$, $2$; \% \texttt{ Two $M$-vector operations}
		\State $x_j=\text{IFFT}(x_j)$, $j=1$, $2$; \% \texttt{ Two $M$-vector IFFTs}
		\State Solve $\widehat{D} x_j = x_j$, $j=1$, $2$. \% \texttt{ Two $M$-vector operations}
	\end{algorithmic}
\end{breakalgo}

For the DNTB preconditioner $\mathcal{F}_{\omega}$, since a scalar factor will not affect its property in nature, we simplify the form of $\mathcal{F}_{\omega}$ by ignoring the scalar factor ${1}/{(2\omega)}$ as follows
\begin{align*} 
\mathcal{F}_{\mbox{\tiny DNTB}} =
\begin{bmatrix}
\omega I - D & 0 \\
0 & \omega I - D
\end{bmatrix}
\begin{bmatrix}
\omega I + T & -I \\
I & \omega I + T
\end{bmatrix}.
\end{align*}
Obviously, the main computational cost for implementing $\mathcal{F}_{\mbox{\tiny DNTB}}$ is to solve a block linear subsystem with $\bigl[\begin{smallmatrix}
	\omega I+T & -I \\
	I & \omega I+T
	\end{smallmatrix}\bigr]$, both direct and iterative methods are not cheap. Therefore, the implementation details of $\mathcal{F}_{\mbox{\tiny DNTB}}$ is not discussed.

In a similar fashion, the DNCB preconditioner $\widetilde{\mathcal{F}}_{\omega}$ can be simplified by removing the scalar factor ${1}/{(2\omega)}$ and keeping its properties essentially the same. Thus, a practical form of $\widetilde{\mathcal{F}}_{\omega}$ reads that
\begin{equation}
\begin{aligned}
\mathcal{F}_{\mbox{\tiny DNCB}} &=
\begin{bmatrix}
\omega I - D & 0 \\
0 & \omega I - D
\end{bmatrix}
\begin{bmatrix}
\omega I+C & -I \\
I & \omega I+C
\end{bmatrix} \\  \label{preconditioner-CBDNS}
&=\begin{bmatrix}
\widetilde{D} & 0 \\
0 & \widetilde{D}
\end{bmatrix}
\begin{bmatrix}
F & 0 \\
0 & F
\end{bmatrix}^{-1}
\begin{bmatrix}
I & 0 \\
L_{21} & I
\end{bmatrix}
\begin{bmatrix}
U_{11} & -I \\
0 & U_{22}
\end{bmatrix}
\begin{bmatrix}
F & 0 \\
0 & F
\end{bmatrix},
\end{aligned}
\end{equation}
where $\widetilde{D}=\omega I - D$, $L_{21}=(\omega I+\Lambda)^{-1}$, $U_{11}=\omega I+\Lambda$, $U_{22}=\omega I+\Lambda+(\omega I+\Lambda)^{-1} \in \mathbb{R}^{M \times M}$ are diagonal, and they can be calculated and stored before the iteration begins. With the same notations in Algorithm \ref{alg-CPMHSS}, we can describe the preconditioning procedure with respect to $\mathcal{F}_{\mbox{\tiny DNCB}}$ in Algorithm \ref{alg-CBDNS}. Obviously, this algorithm requires seven $M$-vector operations ($\mathcal{O}(M)$ flops in total), and four $M$-vector FFT/IFFT operations ($\mathcal{O}(M\log M)$ flops in total).

\begin{breakalgo}{Solve the GR linear system $\mathcal{F}_{\mbox{\tiny DNCB}}\,x=r$}{alg-CBDNS}
	\begin{algorithmic}[1]
		\State Solve $\widetilde{D}x_j = r_j$, $j=1$, $2$; \% \texttt{ Two $M$-vector operations}
		\State $x_j=\text{FFT}(x_j)$, $j=1$, $2$; \%  \texttt{ Two $M$-vector FFTs}
		\State $x_2 = x_2 - L_{21}x_1$; \% \texttt{ Two $M$-vector operations}
		\State Solve $U_{22}x_2 = x_2$ and $U_{11}x_1 = x_1+x_2$. \% \texttt{ Three $M$-vector operations}
		\State $x_j=\text{IFFT}(x_j)$, $j=1$, $2$. \% \texttt{ Two $M$-vector IFFTs}
	\end{algorithmic}
\end{breakalgo}

According to the above discussions, the computational workloads for implementing $\mathcal{F}_{\mbox{\tiny CPMHSS}}$ and $\mathcal{F}_{\mbox{\tiny DNCB}}$ are both dominated by four $M$-vector FFT/IFFT operations when $M$ is large, i.e., $\mathcal{O}(M\log M)$ flops are required at each iteration of the related preconditioned Krylov subspace iteration methods. Therefore, the computational efficiency of a preconditioned Krylov subspace iteration method in conjunction with $\mathcal{F}_{\mbox{\tiny CPMHSS}}$ or $\mathcal{F}_{\mbox{\tiny DNCB}}$ mainly depends on the corresponding convergence rate.

\section{Numerical experiments}
\label{sec-numerical-experiments}

A large number of numerical experiments are presented to show the properties of the DNCB preconditioner, and the effectiveness and efficiency of the related preconditioned GMRES method for the discrete space fractional CNLS equations in the case of $\rho < 0$. Specifically, the LICD scheme applied to the space fractional CNLS equations requires the initial values at the initial and the first time levels to start-up, the former is given by the initial condition, and the latter can be obtained by a second order implicit conservative scheme, e.g., Crank-Nicolson difference scheme \cite{WDL2013JCP}. In all the experiments, the related block linear systems (\ref{positiveBlockForm}) at the 2nd time levels of the discrete fractional CNLS equations are tested in different settings, and the initial guess of the (preconditioned) GMRES method is zero vector.

The parameter $\alpha$ of the fractional Laplacian is selected to be $\alpha=1.1:0.1:2$, and the number $M$ of inner spatial discrete points ranges from $800$ to $102400$ with $M$ doubling consecutively. The temporal step size $\tau$ is simply fixed to be $\tau=0.01$. We denote by `IT' the number of iterations, and `CPU' the computing time in seconds. In all the experiments, the (preconditioned) GMRES method without restart is tested, and terminated either the spectral norm relative residual of the (preconditioned) system of linear equations reduced below $10^{-6}$ or the number of iterations exceeding $1000$.

Before showing the numerical experiments, we provide some discussions about the discrete mass and energy conservations of the LICD scheme. Theoretically, Wang et al. \cite{WDL2014JCP} concluded that the discrete mass and energy are conserved for the LICD scheme in accurate arithmetic. As we all know, floating point arithmetic in computers always has rounding errors. Nevertheless, Wang et al. \cite{WDL2014JCP} observed the conservation phenomenon of the LICD scheme in their experiments when the linear systems are solved by a direct method at each time level under finite precision. In other words, as long as we solve the linear systems by any convergent iteration method with sufficient precision, the discrete mass and energy conservation of the LICD scheme can still be maintained. Therefore, we will not conduct similar experiments as shown in \cite{WDL2014JCP} to verify the conservation laws. Instead, for the LICD scheme, we provide the plots of the numerical solution computed by the DNCB preconditioned GMRES method and its corresponding small error with the exact solution computed by the Gaussian elimination, see Figures \ref{fig:singal-1.1}-\ref{fig:singal-2} for the decoupled case, and Figures \ref{fig:couple-alp=1.1-beta=1}-\ref{fig:couple-alp=2-beta=1} for the coupled case. This error can be further reduced by adopting a smaller tolerance in the iteration method until the conservation of the discrete mass and energy can be observed under a give precision.

\subsection{The decoupled nonlinear case}

Let $\gamma=1$, $\rho=-2$, $1<\alpha\le 2$, we consider the following truncated space fractional DNLS equations
\begin{eqnarray}\label{equ:couple}
\imath u_t-\gamma(-\Delta)^{\frac{\alpha}{2}}u + \rho\vert u \vert^2u=0, \qquad  -20\le x \le 20,\quad 0 < t \le 2,
\end{eqnarray}
with the initial and boundary conditions
\begin{eqnarray}\label{equ:couple initial}
u(x,0)=\text{sech} (x)  e^{2\imath x}, u(-20,t)=u(20,t)=0.
\end{eqnarray}
The LICD scheme applied to (\ref{equ:couple}) and (\ref{equ:couple initial}) leads to the discrete space fractional DNLS equations. On each time level $t_n$, $1 < n \le N$, a complex symmetric linear system of the form (\ref{equ3}) needs to be solved, which is equivalent to solve a block linear system of the form (\ref{positiveBlockForm}). In the experiments, the properties of the DNTB/DNCB preconditioned system matrices are verified, and the GMRES method, the DNCB preconditioned GMRES method (DNCB-GMRES), and the CPMHSS-preconditioned GMRES method (CPMHSS-GMRES) are implemented and compared.

\begin{figure}[htbp]
	\centering
	\begin{tabular}{c}
		\includegraphics[scale=0.30]{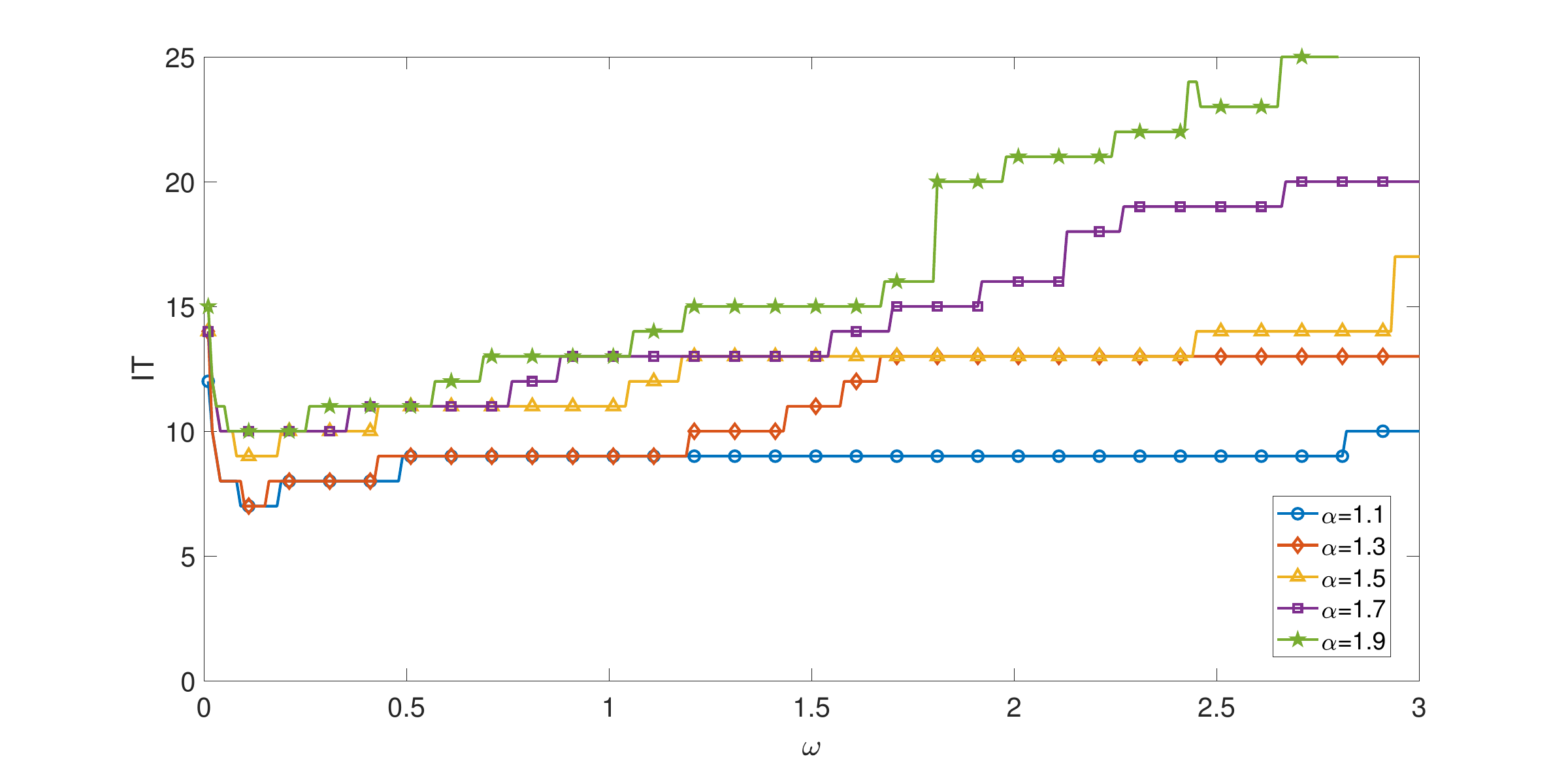}
		\\
	\end{tabular}
	\caption{The curves of IT of DNCB-GMRES versus the parameter $\omega\in (0,3]$ of $\mathcal{F}_{\mbox{\tiny DNCB}}$ in the DNLS case when $\alpha=1.1:0.2:1.9$ and $M=6400$: blue solid line with circle mark for $\alpha=1.1$, red solid line with diamond mark for $\alpha=1.3$, orange solid line with triangle mark for $\alpha=1.5$, purple solid line with square mark for $\alpha=1.7$, green solid line with pentagram mark for $\alpha=1.9$.}
	\label{fig:singal-WW-it}
\end{figure}

Figure \ref{fig:singal-WW-it} depicts the curves of IT of DNCB-GMRES versus the parameter $\omega \in (0,3]$ when $M=6400$. Here, we take the different fractional orders $\alpha=1.1:0.2:1.9$. It shows that IT increases rapidly as $\omega$ approaches zero in Figure \ref{fig:singal-WW-it}. As $\omega$ increases, IT reaches its minimum value quickly and then grows slowly. The empirical optimal value $\omega$ stays in $[0.05,0.50]$ for all cases of $\alpha$, which implies that DNCB-GMRES is not very sensitive to the preconditioning parameter $\omega$ if it is slightly away from 0. In addition, the curve related to a larger value of $\alpha$ is located at a higher position, in another word, the larger the value of $\alpha$ the more difficult it is to solve the related linear system.  

\begin{figure}[htbp]
	\centering
	\begin{tabular}{c}
		\includegraphics[scale=0.37]{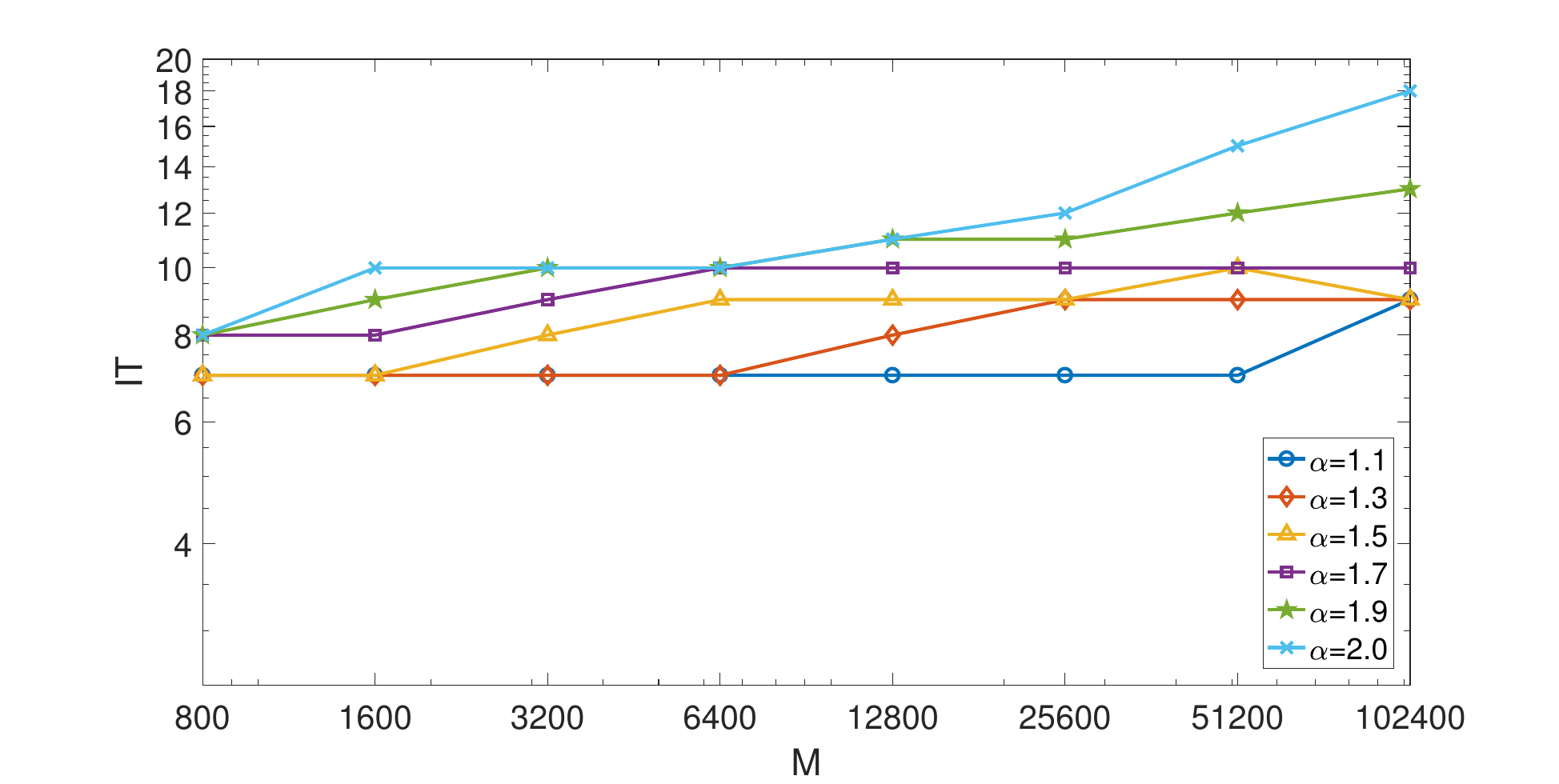}
		\\
	\end{tabular}
	\caption{ The curves of IT of DNCB-GMRES versus the number of the inner spatial discrete points $M$ of the LICD scheme in the DNLS case when $\alpha=1.1:0.2:1.9$ and $\alpha=2.0$: blue solid line with circle mark for $\alpha=1.1$, red solid line with diamond mark for $\alpha=1.3$, orange solid line with triangle mark for $\alpha=1.5$, purple solid line with square mark for $\alpha=1.7$, green solid line with pentagram mark for $\alpha=1.9$, cyan solid line with cross mark for $\alpha=2.0$.}
	\label{fig:singal-alp-IT}
\end{figure}

Figure \ref{fig:singal-alp-IT} depicts the curves of IT of DNCB-GMRES versus the number $M$ of the inner spatial discrete points. Since $M$ doubles consecutively, we use logarithmic scales for both coordinate axes in the plot. Here, we take the standard Laplacian ($\alpha=2.0$), the fractional orders $\alpha=1.1:0.2:1.9$, the values of $M$ ranging from 800 to 102400, and adopt the empirical optimal value $\omega$ for DNCB-GMRES. As can be seen from the curves in the plot, when $\alpha=2.0$, IT of DNCB-GMRES depends linearly on $M$, which aligns with the well-known condition number estimates of order $\mathcal{O}(M)$ when circulant preconditioning is applied to the standard Laplacian, see in \cite{CC1992}. When $\alpha=1.9$, the curve also indicates that IT of DNCB-GMRES depends linearly on $M$, but the slope of this curve is significantly smaller than that when $\alpha=2$. When $\alpha$ takes smaller values, the slope of the corresponding curve is also relatively smaller, or in other words, the linear dependence of IT of DNCB-GMRES on $M$ is relatively weaker. Therefore, the experiments validated the theoretical results presented in \cite{CC1992}, which support the linear dependence of IT of DNCB-GMRES on $M$ in the standard case ($\alpha=2$). Furthermore, they provided evidences that IT of DNCB-GMRES is also linearly dependent on $M$ in the fractional case (especially when $\alpha$ is large). As can be seen from Figure \ref{fig:singal-alp-IT}, the curve for $\alpha=1.1$ is at the bottom, and the larger value of $\alpha$ corresponds to the curve of higher position. Similarly to Figure \ref{fig:singal-WW-it}, it indicates that a larger value of fractional order $\alpha$ leads to a more difficult linear system.

\begin{figure}[htbp]
	\centering
	\begin{tabular}{cc}
		\includegraphics[scale=0.15]{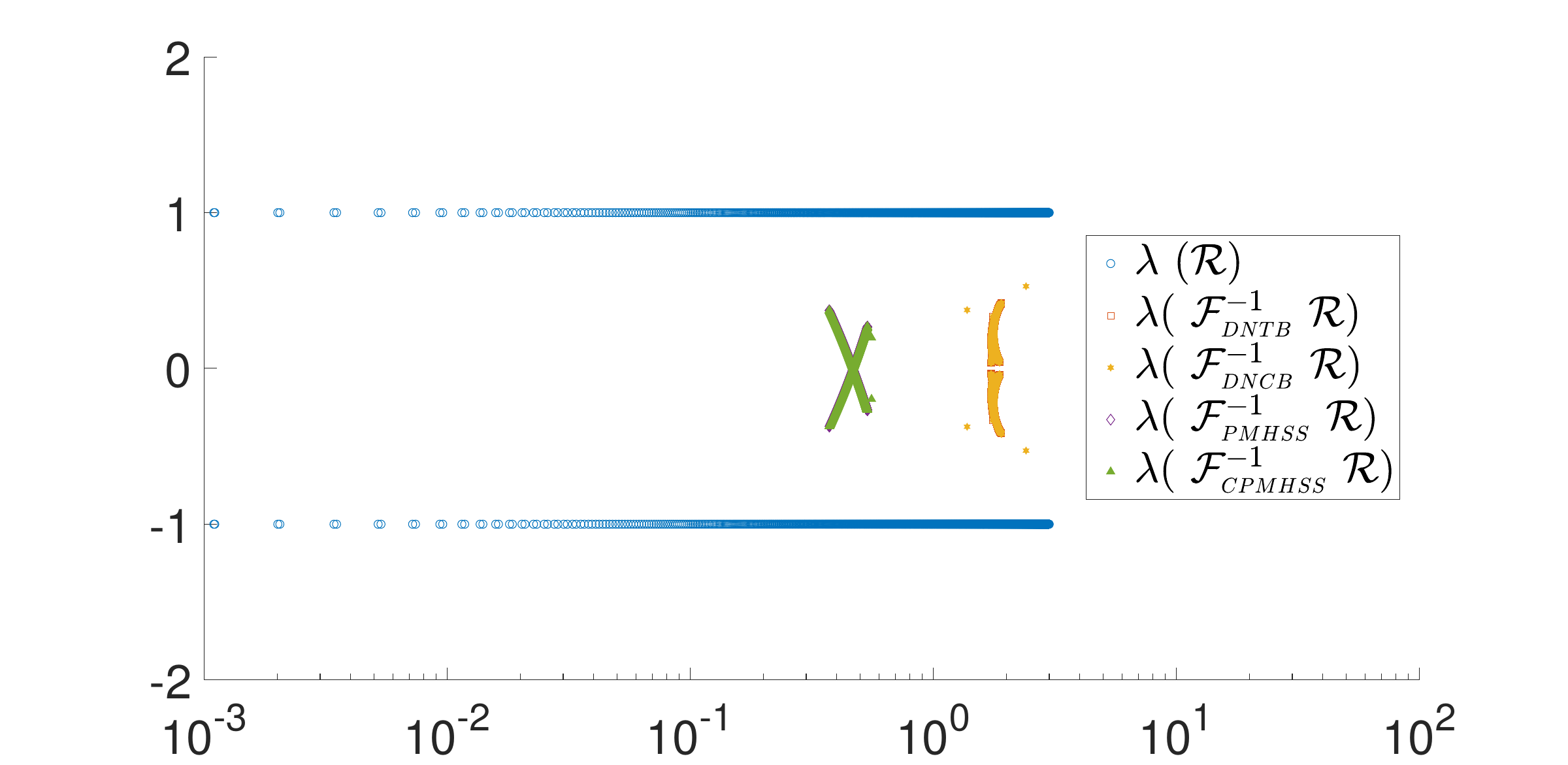}  &
		\includegraphics[scale=0.15]{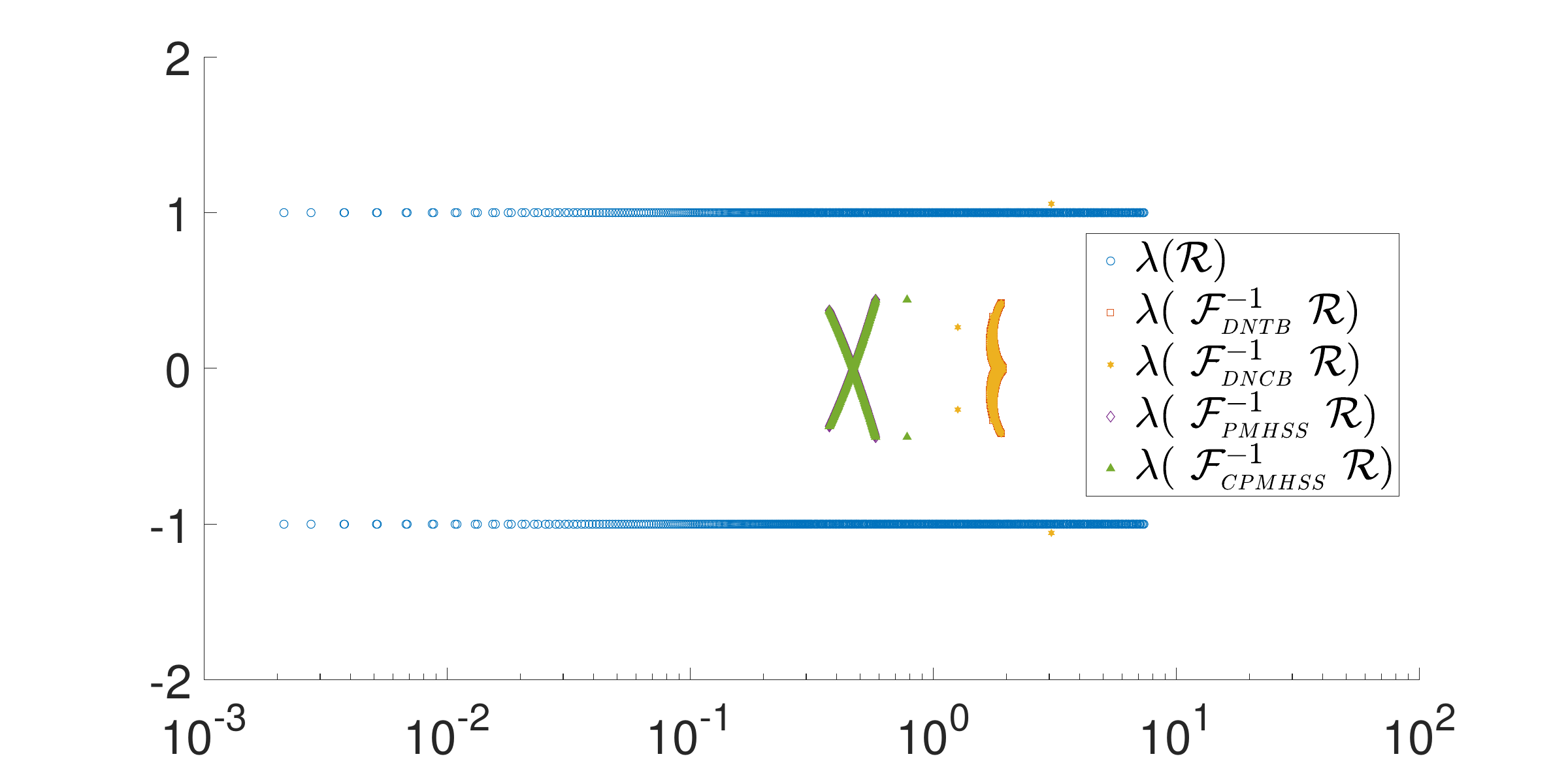}
		
		\\
	\end{tabular}
	\caption{The eigenvalue distribution of $\mathcal{R}$, $\mathcal{F}_{\mbox{\tiny DNTB}}^{-1}\mathcal{R}$, $\mathcal{F}_{\mbox{\tiny DNCB}}^{-1}\mathcal{R}$, $\mathcal{F}_{\mbox{\tiny PMHSS}}^{-1}\mathcal{R}$ and $\mathcal{F}_{\mbox{\tiny CPMHSS}}^{-1}\mathcal{R}$ in the case of $\alpha=1.3$ for $M=1600$ (left), $3200$ (right): blue circle mark for $\mathcal{R}$, red square mark for $\mathcal{F}_{\mbox{\tiny DNTB}}^{-1}\mathcal{R}$, orange solid hexagram mark for $\mathcal{F}_{\mbox{\tiny DNCB}}^{-1}\mathcal{R}$, purple diamond mark for $\mathcal{F}_{\mbox{\tiny PMHSS}}^{-1}\mathcal{R}$, green solid triangle mark for $\mathcal{F}_{\mbox{\tiny CPMHSS}}^{-1}\mathcal{R}$.}
	\label{fig:eig-singal-1.3}
\end{figure}

\begin{figure}[htbp]
	\centering
	\begin{tabular}{cc}
		\includegraphics[scale=0.15]{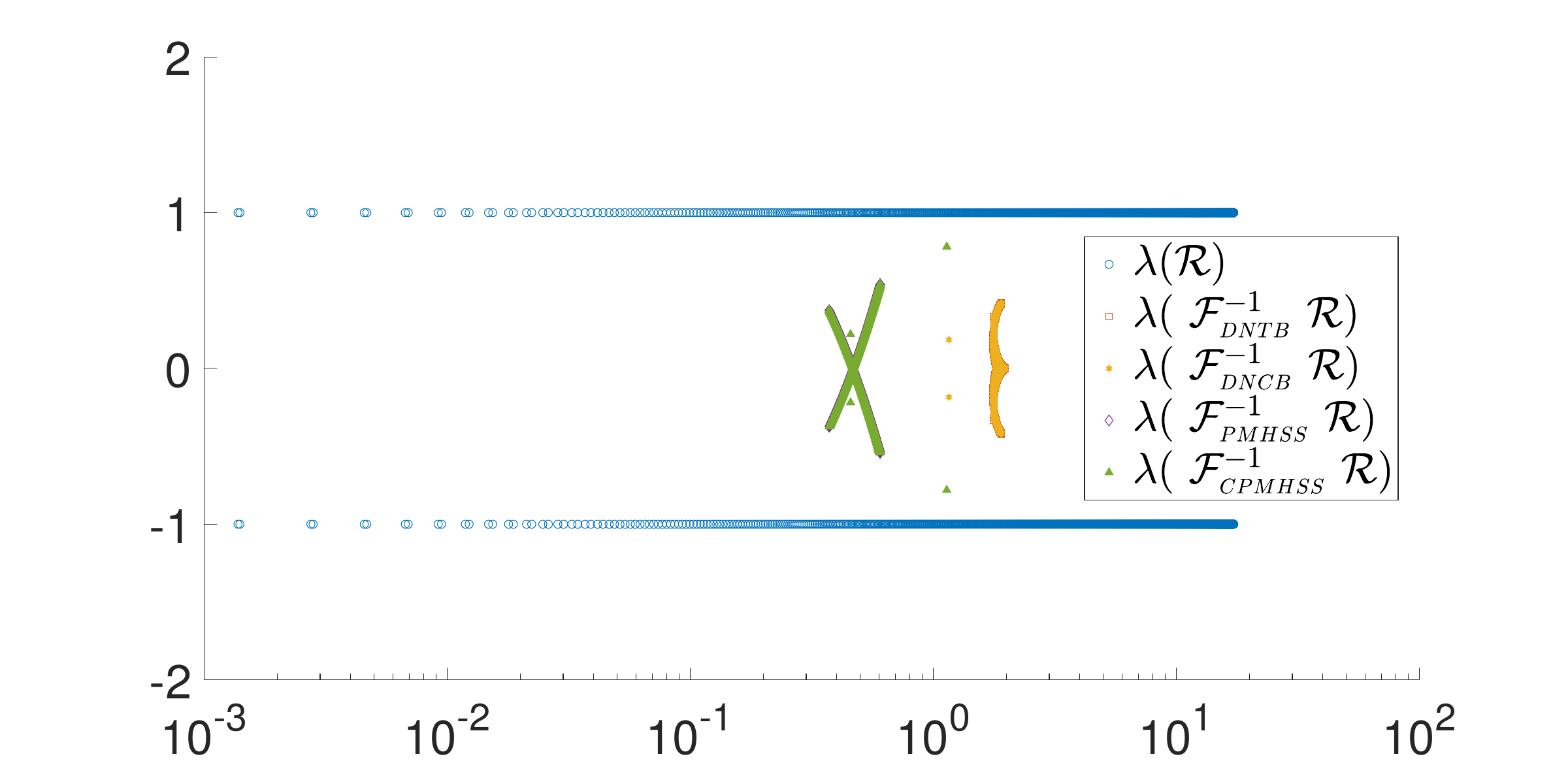}  &
		\includegraphics[scale=0.15]{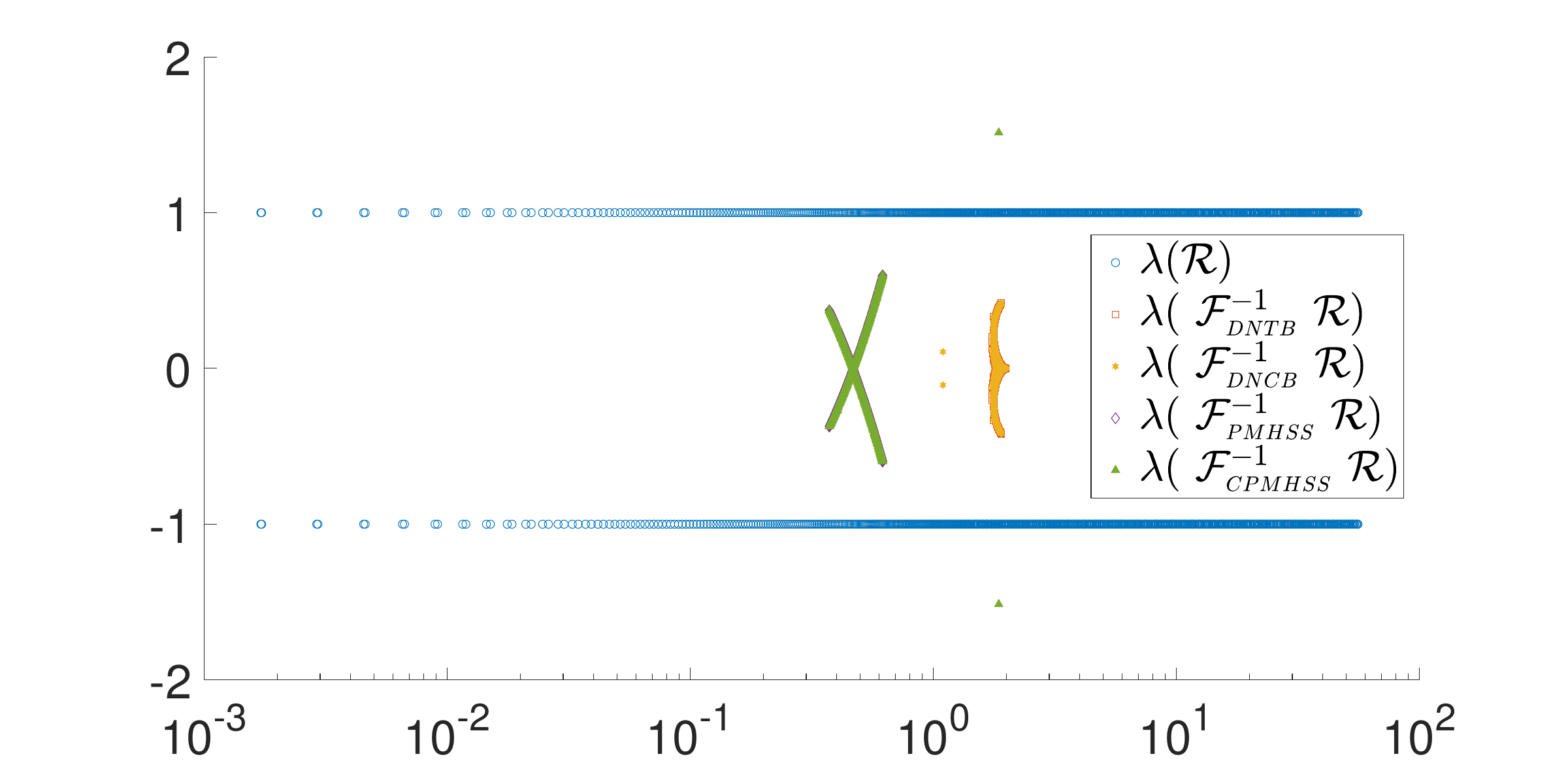}
		
		\\
	\end{tabular}
	\caption{ The eigenvalue distribution of $\mathcal{R}$, $\mathcal{F}_{\mbox{\tiny DNTB}}^{-1}\mathcal{R}$, $\mathcal{F}_{\mbox{\tiny DNCB}}^{-1}\mathcal{R}$, $\mathcal{F}_{\mbox{\tiny PMHSS}}^{-1}\mathcal{R}$ and $\mathcal{F}_{\mbox{\tiny CPMHSS}}^{-1}\mathcal{R}$ in the case of $\alpha=1.7$ for $M=1600$ (left), $3200$ (right): blue circle mark for $\mathcal{R}$, red square mark for $\mathcal{F}_{\mbox{\tiny DNTB}}^{-1}\mathcal{R}$, orange solid hexagram mark for $\mathcal{F}_{\mbox{\tiny DNCB}}^{-1}\mathcal{R}$, purple diamond mark for $\mathcal{F}_{\mbox{\tiny PMHSS}}^{-1}\mathcal{R}$, green solid triangle mark for $\mathcal{F}_{\mbox{\tiny CPMHSS}}^{-1}\mathcal{R}$.}
	\label{fig:eig-singal-1.7}
\end{figure}

Figures \ref{fig:eig-singal-1.3}-\ref{fig:eig-singal-1.7} depict the eigenvalue distribution of the system matrix $\mathcal{R}$, the DNTB/DNCB preconditioned system matrix $\mathcal{F}_{\mbox{\tiny DNTB}}^{-1}\mathcal{R}$/$\mathcal{F}_{\mbox{\tiny DNCB}}^{-1}\mathcal{R}$, and the PMHSS/CPMHSS preconditioned system matrix $\mathcal{F}_{\mbox{\tiny PMHSS}}^{-1}\mathcal{R}$/$\mathcal{F}_{\mbox{\tiny CPMHSS}}^{-1}\mathcal{R}$. The Strang's circulant approximation is adopted in $\mathcal{F}_{\mbox{\tiny DNCB}}$ and $\mathcal{F}_{\mbox{\tiny CPMHSS}}$, and the fractional order $\alpha=1.3, 1.7$ and the number of the spatial discrete points $M=1600, 3200$ are selected. In each figure, the plots of the cases for $M=1600$ and $M=3200$ are on the left and right, respectively. As can be seen from Figures \ref{fig:eig-singal-1.3}-\ref{fig:eig-singal-1.7}, the real parts of the eigenvalues of $\mathcal{R}$ are distributed over a wide range of scales from $\mathcal{O}(10^{-3})$ to $\mathcal{O}(10^1)$, and the larger the value of $M$, the wider the range of scales. The real parts of the eigenvalues of $\mathcal{F}_{\mbox{\tiny PMHSS}}^{-1}\mathcal{R}$/$\mathcal{F}_{\mbox{\tiny CPMHSS}}^{-1}\mathcal{R}$ are distributed over a range of scales from $\mathcal{O}(10^{-1})$ to $\mathcal{O}(10^{0})$, and the related imaginary parts stays in $[-0.5, 0.5]$. The real parts of the eigenvalues of $\mathcal{F}_{\mbox{\tiny DNTB}}^{-1}\mathcal{R}$/$\mathcal{F}_{\mbox{\tiny DNCB}}^{-1}\mathcal{R}$ are distributed around $\mathcal{O}(10^0)$, and the related imaginary parts also stays in $[-0.5, 0.5]$. Therefore, it shows that the eigenvalues of $\mathcal{F}_{\mbox{\tiny DNTB}}^{-1}\mathcal{R}$/$\mathcal{F}_{\mbox{\tiny DNCB}}^{-1}\mathcal{R}$ and $\mathcal{F}_{\mbox{\tiny PMHSS}}^{-1}\mathcal{R}$/$\mathcal{F}_{\mbox{\tiny CPMHSS}}^{-1}\mathcal{R}$ are more clustered than those of $\mathcal{R}$, and the eigenvalues of $\mathcal{F}_{\mbox{\tiny DNTB}}^{-1}\mathcal{R}$/$\mathcal{F}_{\mbox{\tiny DNCB}}^{-1}\mathcal{R}$ are clustered even more compact than those of $\mathcal{F}_{\mbox{\tiny PMHSS}}^{-1}\mathcal{R}$/$\mathcal{F}_{\mbox{\tiny CPMHSS}}^{-1}\mathcal{R}$. When $M$ increases from 1600 to 3200, the distribution of the eigenvalues of $\mathcal{F}_{\mbox{\tiny DNTB}}^{-1}\mathcal{R}$/$\mathcal{F}_{\mbox{\tiny DNCB}}^{-1}\mathcal{R}$ remains almost the same. Furthermore, most of the eigenvalues of $\mathcal{F}_{\mbox{\tiny DNCB}}^{-1}\mathcal{R}$ are clustered around those of $\mathcal{F}_{\mbox{\tiny DNTB}}^{-1}\mathcal{R}$, which indicates that $\mathcal{F}_{\mbox{\tiny DNCB}}$ is a good approximation to $\mathcal{F}_{\mbox{\tiny DNTB}}$.

Table \ref{tab:BDNS-different-circulant-alp=1.5} list IT of DNCB-GMRES in conjunction with different circulant approximations. In these experiments, we take the empirical optimal value of the preconditioning parameter $\omega$ (searching in $(0,3]$), the fractional order $\alpha=1.5$, and the number of the spatial discrete points $M=6400$. The results of several representative circulant approximations are listed, including Strang's circulant approximation, T. Chan's circulant approximation, R. Chan's circulant approximation, circulant approximations constructed from some famous kernels (e.g., Modified Dirichlet kernel, von Hann kernel, Hamming kernel), and superoptimal circulant approximation \cite{CRH2007SIAMbook}. As shown in this table, IT of DNCB-GMRES equals to $9$ in most cases of the circulant approximations (except for the superoptimal circulant approximation which takes $32$ iterations), and the corresponding optimal empirical values of $\omega$ are almost the same, i.e., $[0.07,0.20]$ (except for the T. Chan's circulant approximation $[0.08,0.18]$ and the superoptimal circulant approximation $[1.50,1.74]$). Therefore, DNCB-GMRES is easy to apply due to the fact that most of the circulant approximations are efficient. The relatively poor performance of the superoptimal circulant approximation may be due to the failure of finding a good enough empirical optimal value of $\omega$ in $(0,3]$.

\begin{figure}[htbp]
	\centering
	\begin{tabular}{cc}
		\includegraphics[scale=0.15]{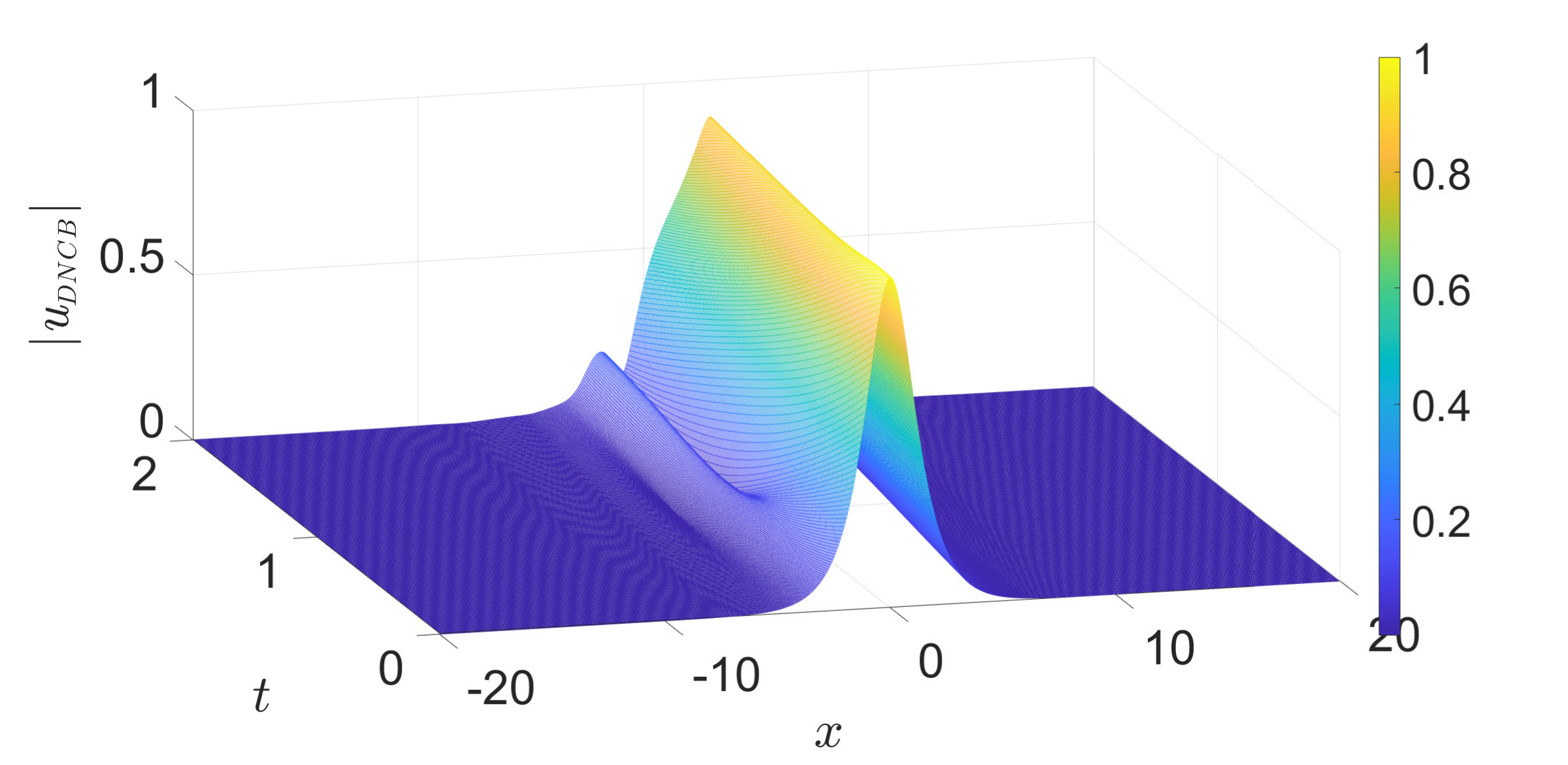}  &
		\includegraphics[scale=0.15]{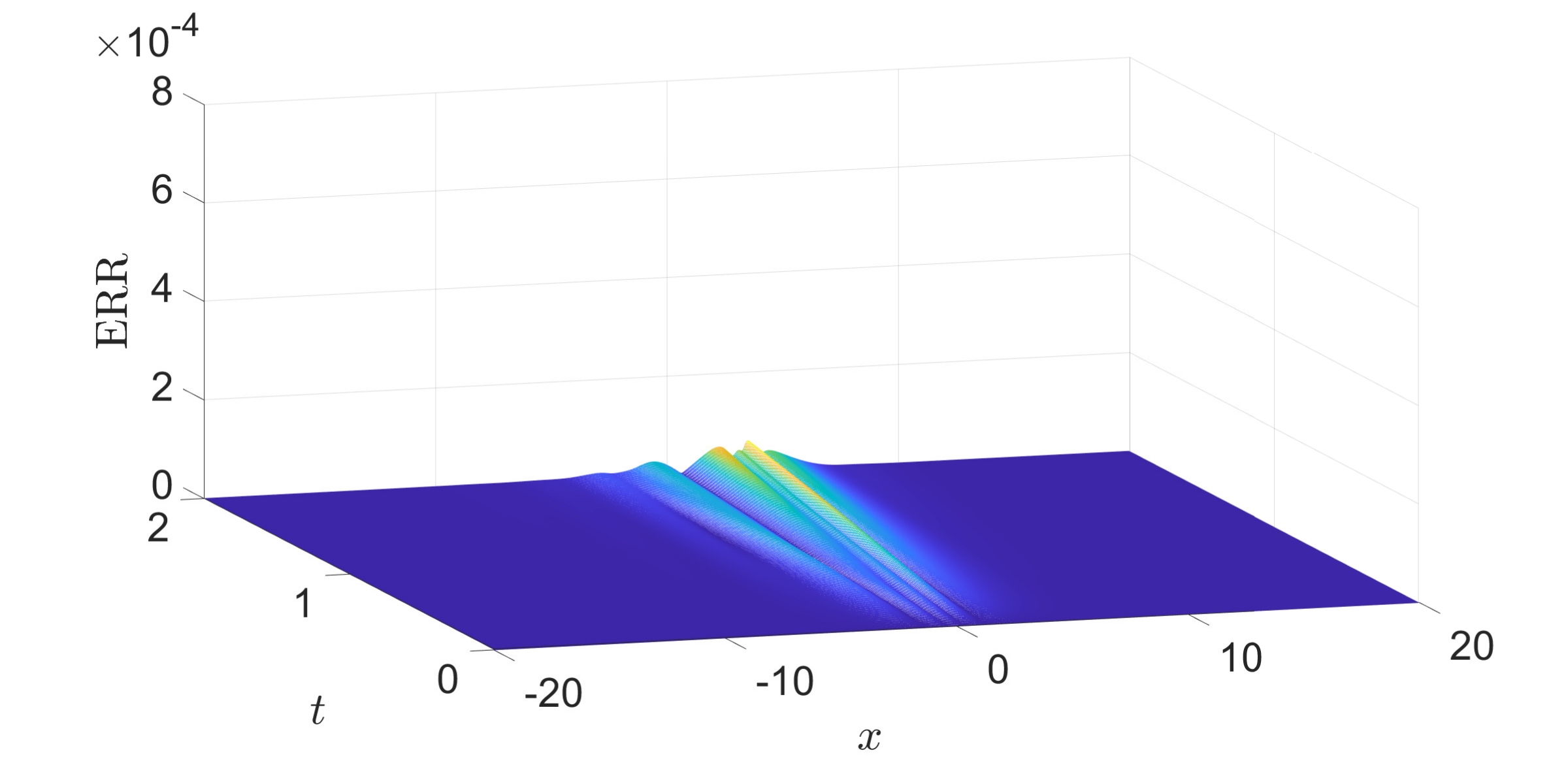}
		
		\\
	\end{tabular}
	\caption{The numerical solution (left) and its error (right) with the exact solution of the LICD scheme in the DNLS case when $\alpha=1.1$ and $M=800$.}
	\label{fig:singal-1.1}
\end{figure}

\begin{figure}[htbp]
	\centering
	\begin{tabular}{cc}
		\includegraphics[scale=0.15]{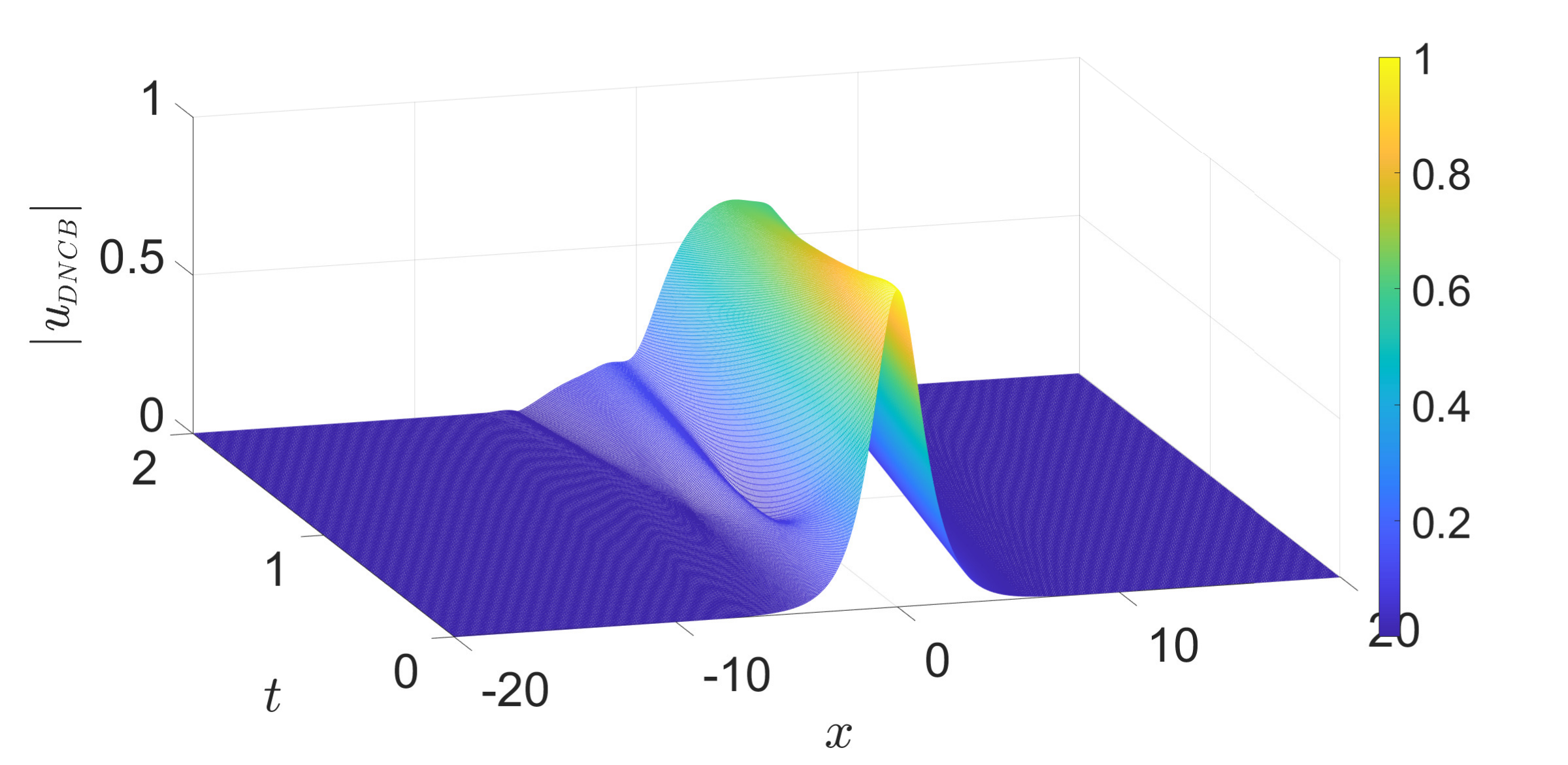}  &
		\includegraphics[scale=0.15]{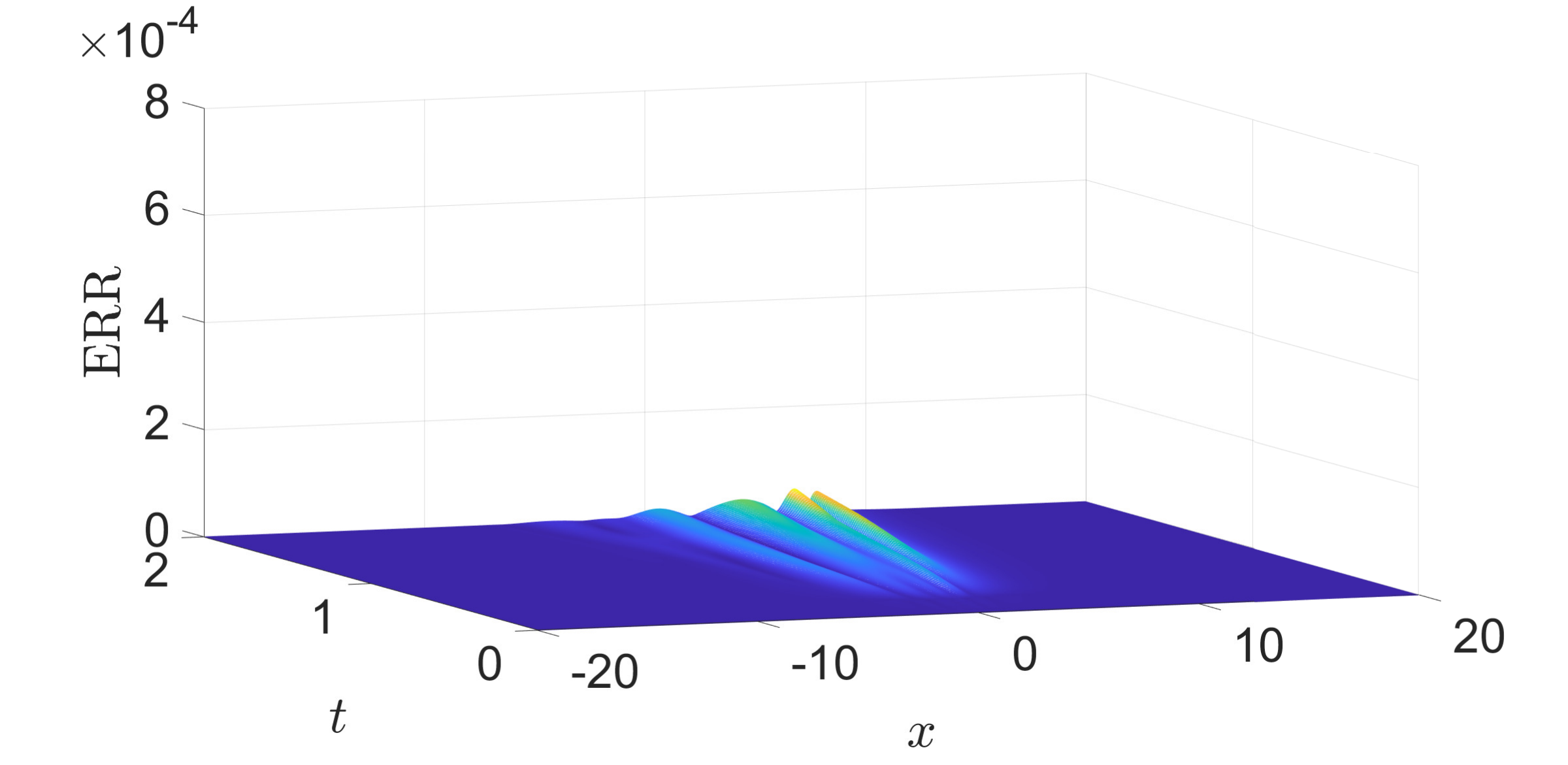}
		
		\\
	\end{tabular}
	\caption{The numerical solution (left) and its error (right) with the exact solution of the LICD scheme in the DNLS case when $\alpha=1.5$ and $M=800$.}
	\label{fig:singal-1.5}
\end{figure}

\begin{figure}[htbp]
	\centering
	\begin{tabular}{cc}
		\includegraphics[scale=0.15]{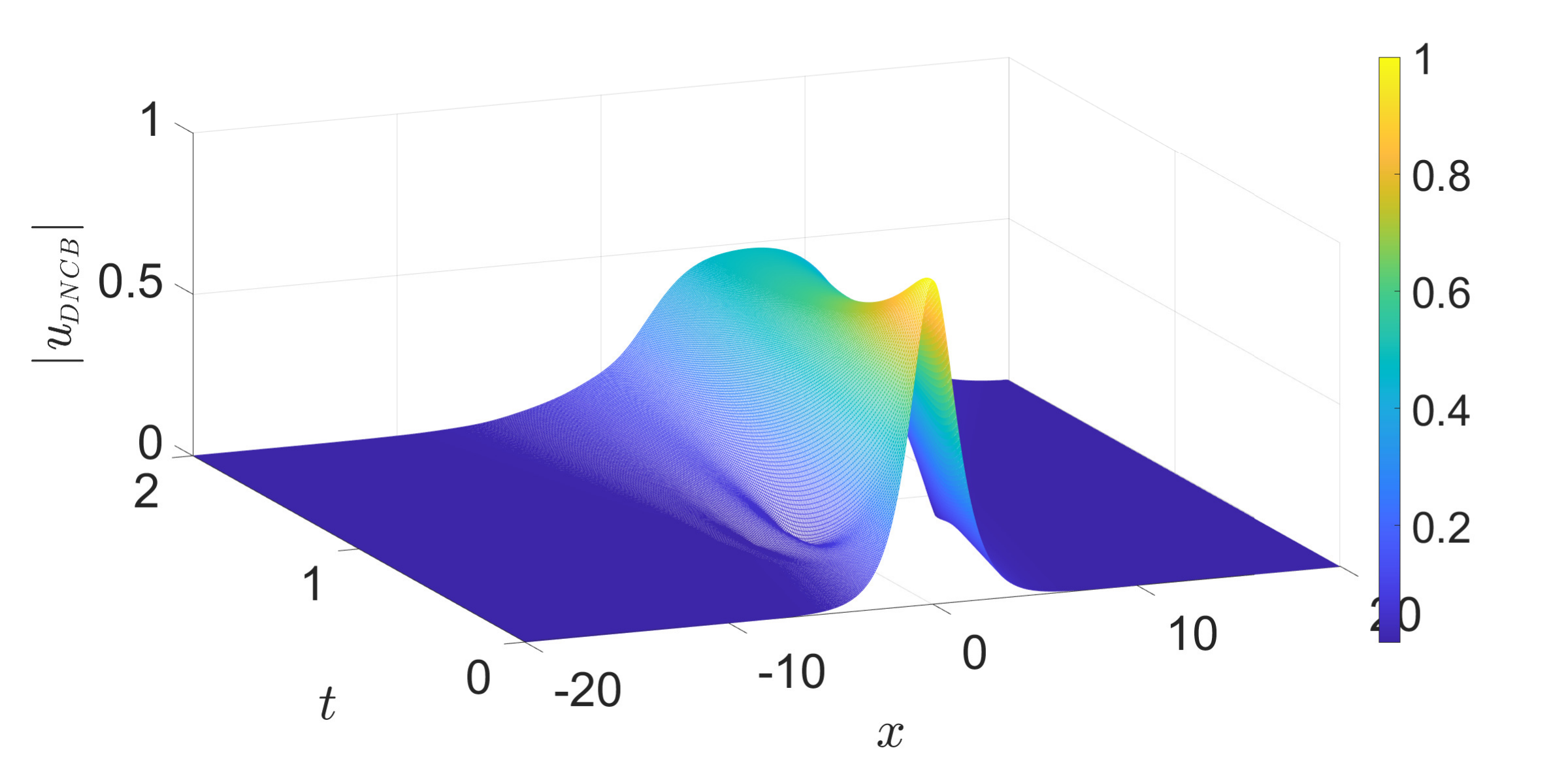}  &
		\includegraphics[scale=0.15]{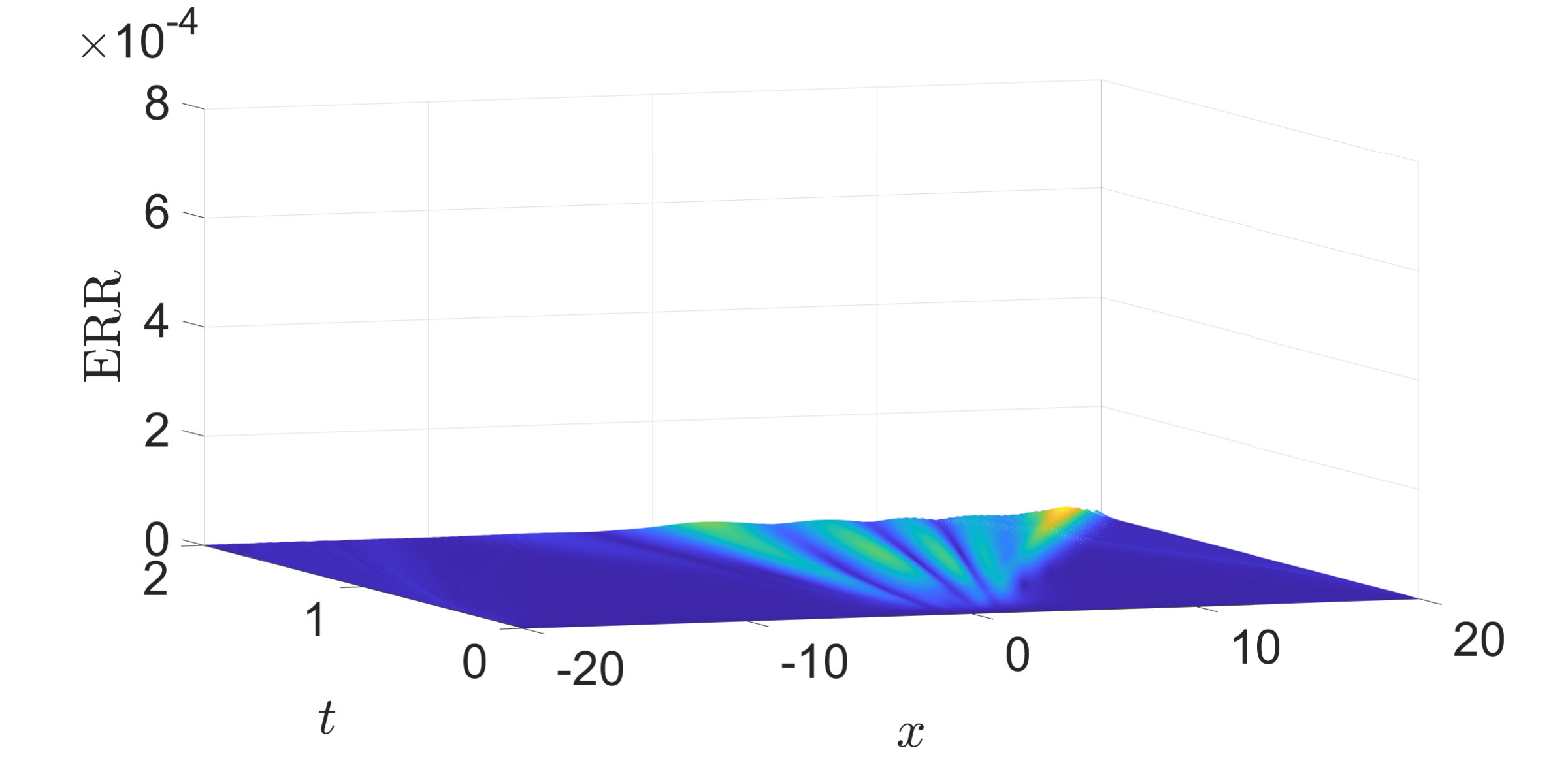}
		\\
	\end{tabular}
	\caption{The numerical solution (left) and its error (right) with the exact solution of the LICD scheme in the DNLS case when $\alpha=1.9$ and $M=800$.}
	\label{fig:singal-1.9}
\end{figure}

\begin{figure}[htbp]
	\centering
	\begin{tabular}{cc}
		\includegraphics[scale=0.15]{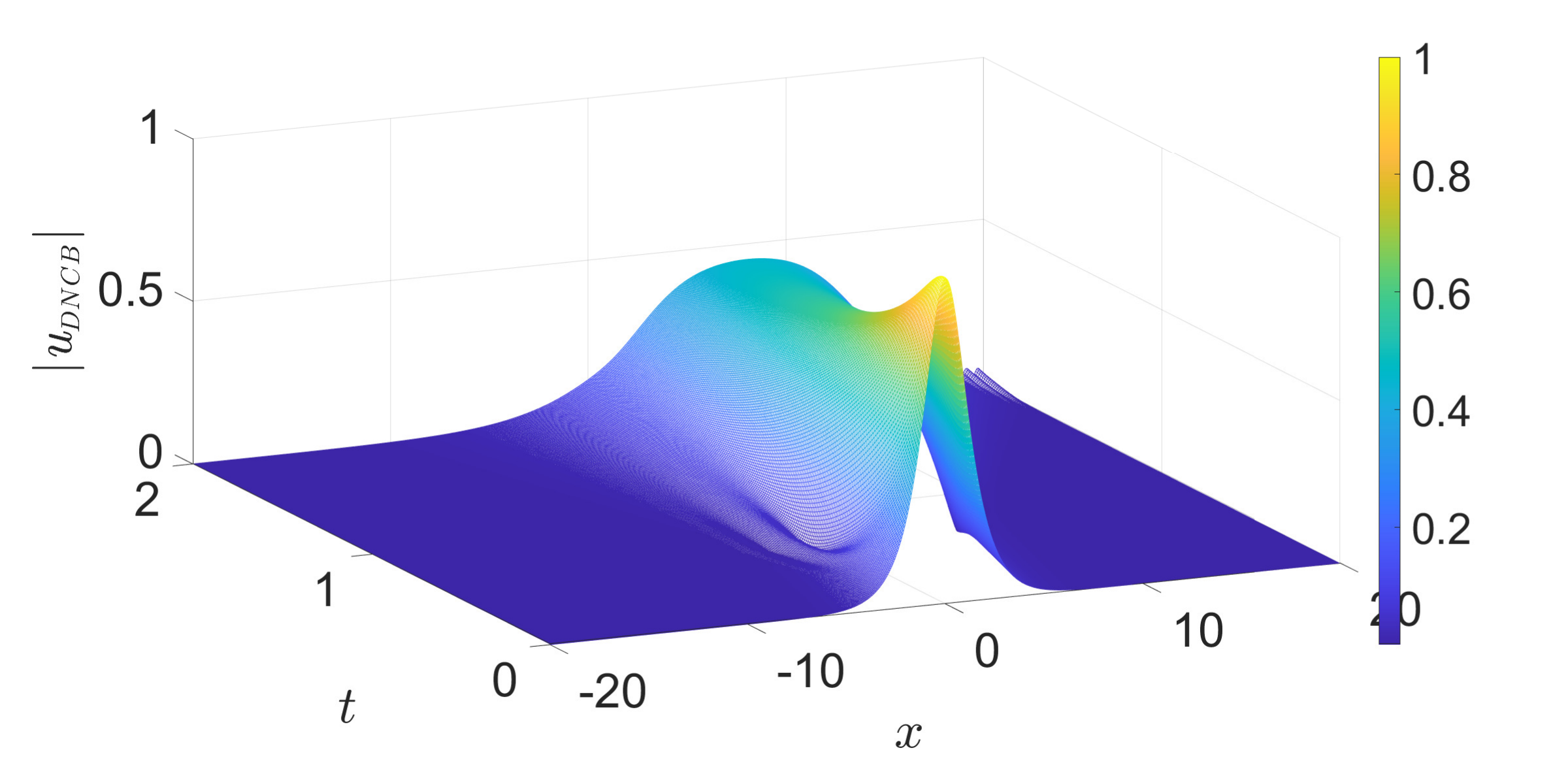}  &
		\includegraphics[scale=0.15]{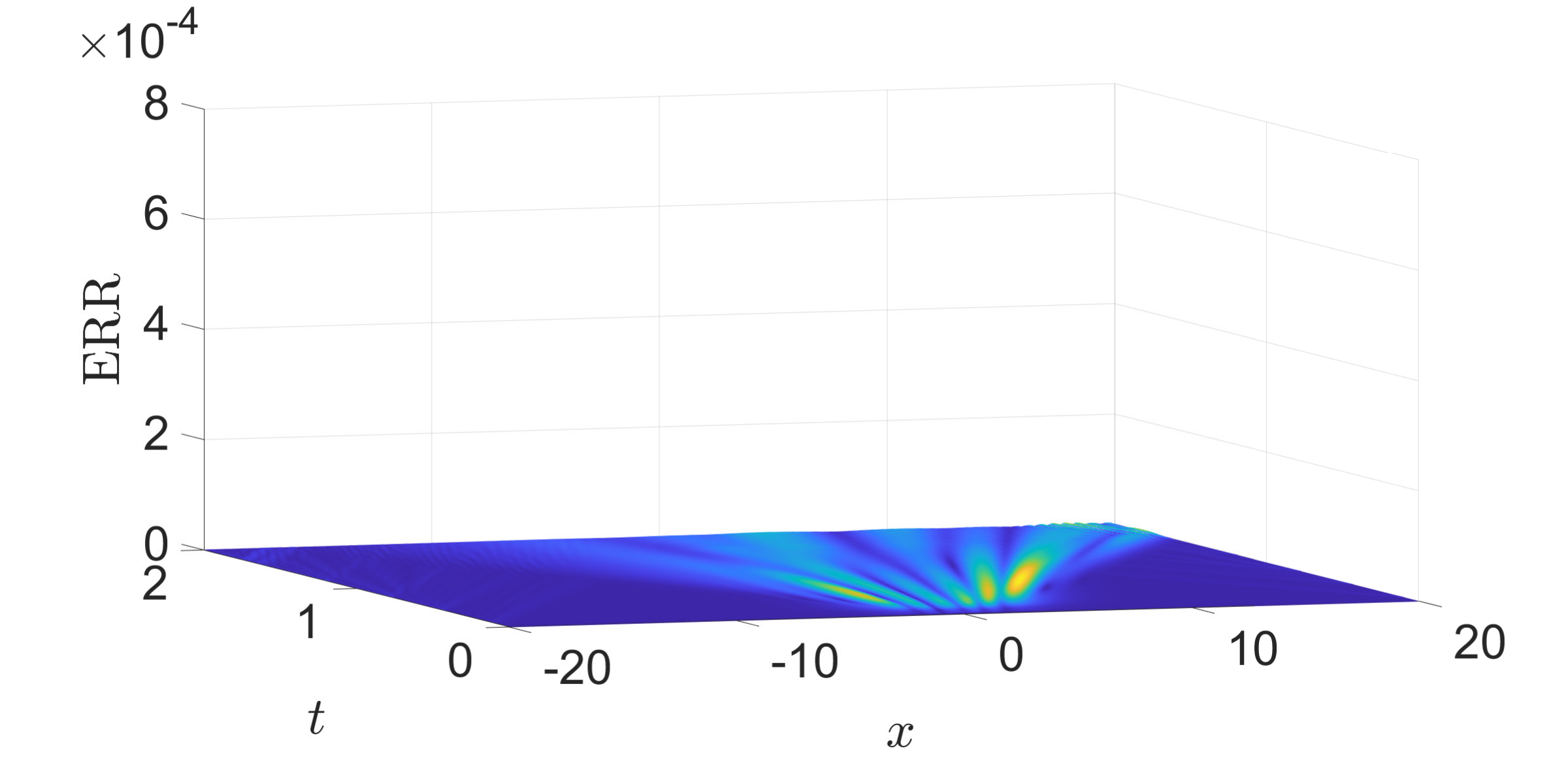}
		\\
	\end{tabular}
	\caption{The numerical solution (left) and its error (right) with the exact solution of the LICD scheme in the DNLS case when $\alpha=2$ and $M=800$.}
	\label{fig:singal-2}
\end{figure}

Figures \ref{fig:singal-1.1}-\ref{fig:singal-2} depict the numerical solutions $u_{\text{\tiny DNCB}}$ and the related absolute errors
$\text{ERR} = |u_{\text{\tiny DNCB}}-u_{\text{\tiny GE}}|$. Here, the exact solution $u_{\text{\tiny GE}}$ is obtained by solving the discrete space fractional DNLS equations with Guassian elimination (GE) based on the LICD scheme. We take the fractional order $\alpha=1.1:0.4:1.9$ and $\alpha=2$, and the number of the spatial discrete points $M=800$. In these figures, the numerical solutions $u_{\text{\tiny DNCB}}$ obtained by DNCB-GMRES are plotted on the left, and the errors between $u_{\text{\tiny DNCB}}$ and $u_{\text{\tiny GE}}$ are plotted on the right. As shown in these figures, the shape of $u_{\text{\tiny DNCB}}$ in hight and width is affected by the values of the fractional order $\alpha$. When $\alpha$ approaches to $2$, the shape of $u_{\text{\tiny DNCB}}$ tends to converge to the shape of the solution of the standard DNLS equation (the case of $\alpha=2$). Moreover, ERR stays as small as around $\mathcal{O}(10^{-4})$ in the whole computational space-time domain, which confirms that the numerical solution $u_{\text{\tiny DNCB}}$ is reliable. Obviously, even more exact numerical solution $u_{\text{\tiny DNCB}}$ can be computed by improving the stopping criterion of DNCB-GMRES.

\begin{table}[htbp]
	\setlength{\abovecaptionskip}{0pt}
	\setlength{\belowcaptionskip}{10pt} \centering{
		\caption{\label{tab:BDNS-different-circulant-alp=1.5}
			IT of DNCB-GMRES with different circulant approximations, and the empirical optimal parameter $\omega$ of DNCB-GMRES in the case of $\alpha=1.5$ and $M=6400$.}\scriptsize
		\begin{tabular}{lrc}\specialrule{0em}{2pt}{2pt}\hline\specialrule{0em}{2pt}{2pt}
			
			Circulant Approximation &   IT &  $\omega$  \\\specialrule{0em}{1pt}{1pt}\hline\specialrule{0em}{3pt}{3pt}
			
			Strang &9 & 	[0.07,0.20]
			\\\specialrule{0em}{3pt}{3pt}
			
				T. Chan & 9 & 	[0.08,0.18]
			\\\specialrule{0em}{3pt}{3pt}
			
			R. Chan &	9 & 	[0.07,0.20]
			\\\specialrule{0em}{3pt}{3pt}
			
			Modified Dirichlet kernel &  9	& 	[0.07,0.20]
			\\\specialrule{0em}{3pt}{3pt}		
			
			von Hann kernel & 	9 &	[0.07,0.20]
			
			\\\specialrule{0em}{3pt}{3pt}
			
			Hamming kernel & 	9 &	[0.07,0.20]
			\\\specialrule{0em}{3pt}{3pt}
			
			Superoptimal  & 	32 & 	[1.50,1.74]
			\\\specialrule{0em}{3pt}{3pt}\hline
	\end{tabular}}
\end{table}

\subsection{The coupled nonlinear case }
    Let $\gamma=1$, $\rho=-2$, $\beta=1$, $1<\alpha\le 2$, we consider the truncated space fractional CNLS equations as follows
	\begin{eqnarray}\label{couple equation}
	\begin{cases}
	\imath u_t-\gamma(-\Delta)^{\frac{\alpha}{2}}u+\rho(\vert u \vert^2+\beta\vert v \vert^2)u=0, \qquad  \\
	\imath v_t-\gamma(-\Delta)^{\frac{\alpha}{2}}v+\rho(\vert v \vert^2+\beta\vert u \vert^2)v=0, \qquad
	\end{cases}
	-20\le x \le 20,\quad 0 < t \le 2,	
	\end{eqnarray}
	with the initial and boundary conditions 	
	\begin{equation}\label{couple condition}
	\begin{cases}
	u(x,0)=\text{sech}(x+1)  e^{2\imath x}, \quad v(x,0)=\text{sech}(x-1)  e^{-2\imath x},\\
	u(-20,t)=u(20,t)=0, \quad v(-20,t)=v(20,t)=0.
	\end{cases}
	\end{equation}

The LICD scheme applied to (\ref{couple equation}) and (\ref{couple condition}) leads to the discrete space fractional CNLS equations. On each time level $t_n$, $1 < n \le N$, two coupled complex symmetric linear systems of the form (\ref{equ3}) need to be solved sequentially, which are equivalent to solve two coupled block linear systems of the form (\ref{positiveBlockForm}). In the experiments, Gaussian elimination (GE) is used to solve the aforementioned complex symmetric linear systems, and GMRES, DNCB-GMRES and CPMHSS-GMRES are adopted to solve the corresponding coupled block linear systems. In this subsection, `IT' is denoted by the total number of iterations of the preconditioned GMRES method for solving two coupled block linear systems on time level $t_2$, and `CPU' is denoted by the total computing time in seconds for solving the coupled (complex symmetric or block) linear systems.

\begin{figure}[htbp]
	\centering
	\begin{tabular}{c}
		\includegraphics[scale=0.37]{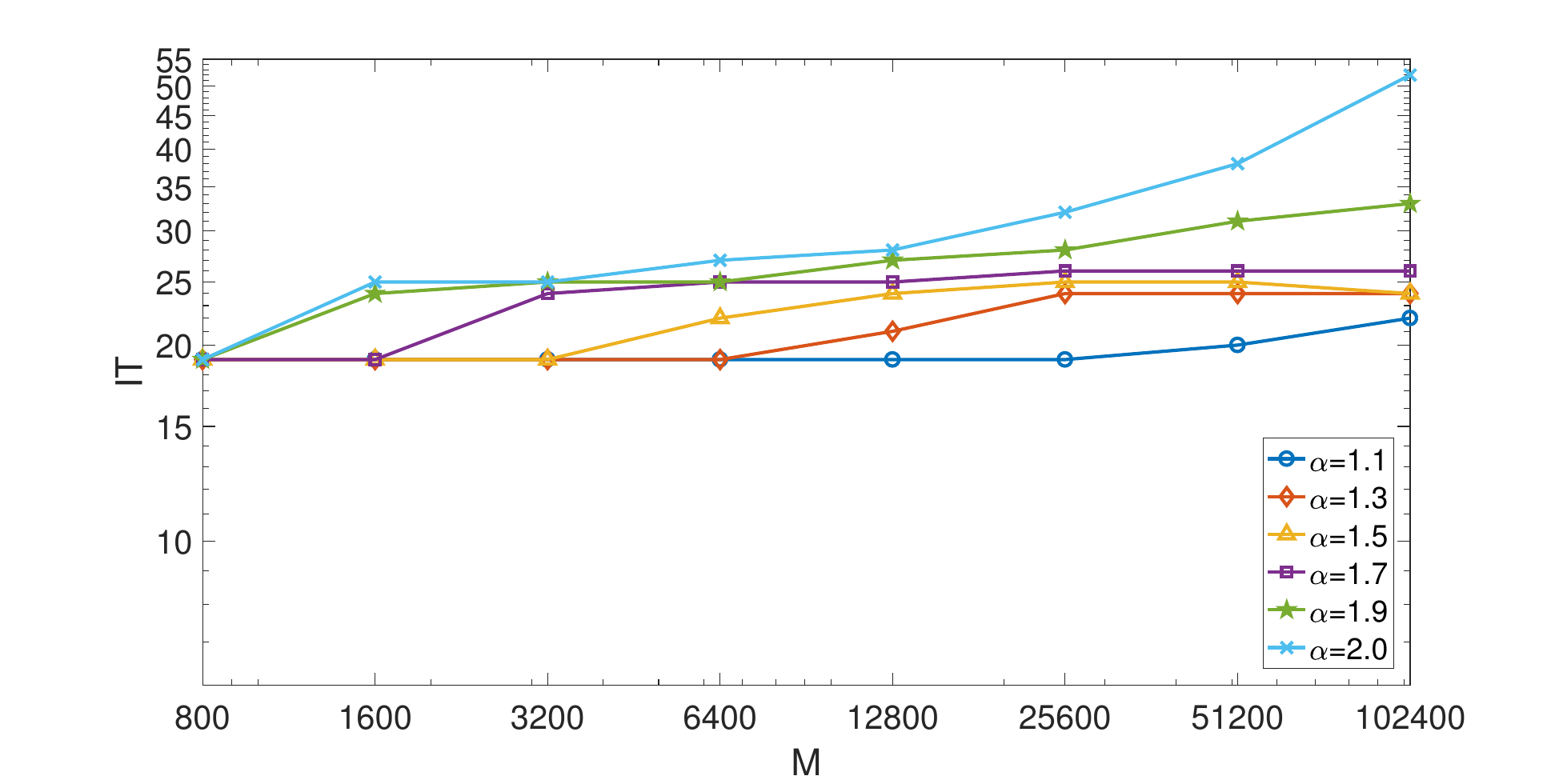}
		\\
	\end{tabular}
	\caption{ The curves of IT of DNCB-GMRES versus the number of the inner spatial discrete points $M$ of the LICD scheme in the CNLS case when $\alpha=1.1:0.2:1.9$ and $\alpha=2.0$: blue solid line with circle mark for $\alpha=1.1$, red solid line with diamond mark for $\alpha=1.3$, orange solid line with triangle mark for $\alpha=1.5$, purple solid line with square mark for $\alpha=1.7$, green solid line with pentagram mark for $\alpha=1.9$, cyan solid line with cross mark for $\alpha=2.0$.}
	\label{fig:couple space-it}
\end{figure}

Figure \ref{fig:couple space-it} depicts the curves of IT of DNCB-GMRES versus the number $M$ of the spatial discrete points. Logarithmic scales for both coordinate axes are used in the plot since $M$ doubles consecutively. Here, we take the standard Laplacian ($\alpha=2.0$), the fractional orders $\alpha=1.1:0.2:1.9$, the values of $M$ ranging from 800 to 102400, and adopt the empirical optimal value $\omega$ for DNCB-GMRES. Similarly to Figure \ref{fig:singal-alp-IT}, IT of DNCB-GMRES exhibits a linear dependence on $M$, which weakens as $\alpha$ decreases. As shown in Figure \ref{fig:couple space-it}, the curve for $\alpha=1.1$ stays at the bottom, and a larger value of $\alpha$ leads to a higher position of the corresponding curve, meaning that smaller $\alpha$ makes the linear systems easier to solve.

\begin{figure}[htbp]
	\centering
	\begin{tabular}{c}
		\includegraphics[scale=0.30]{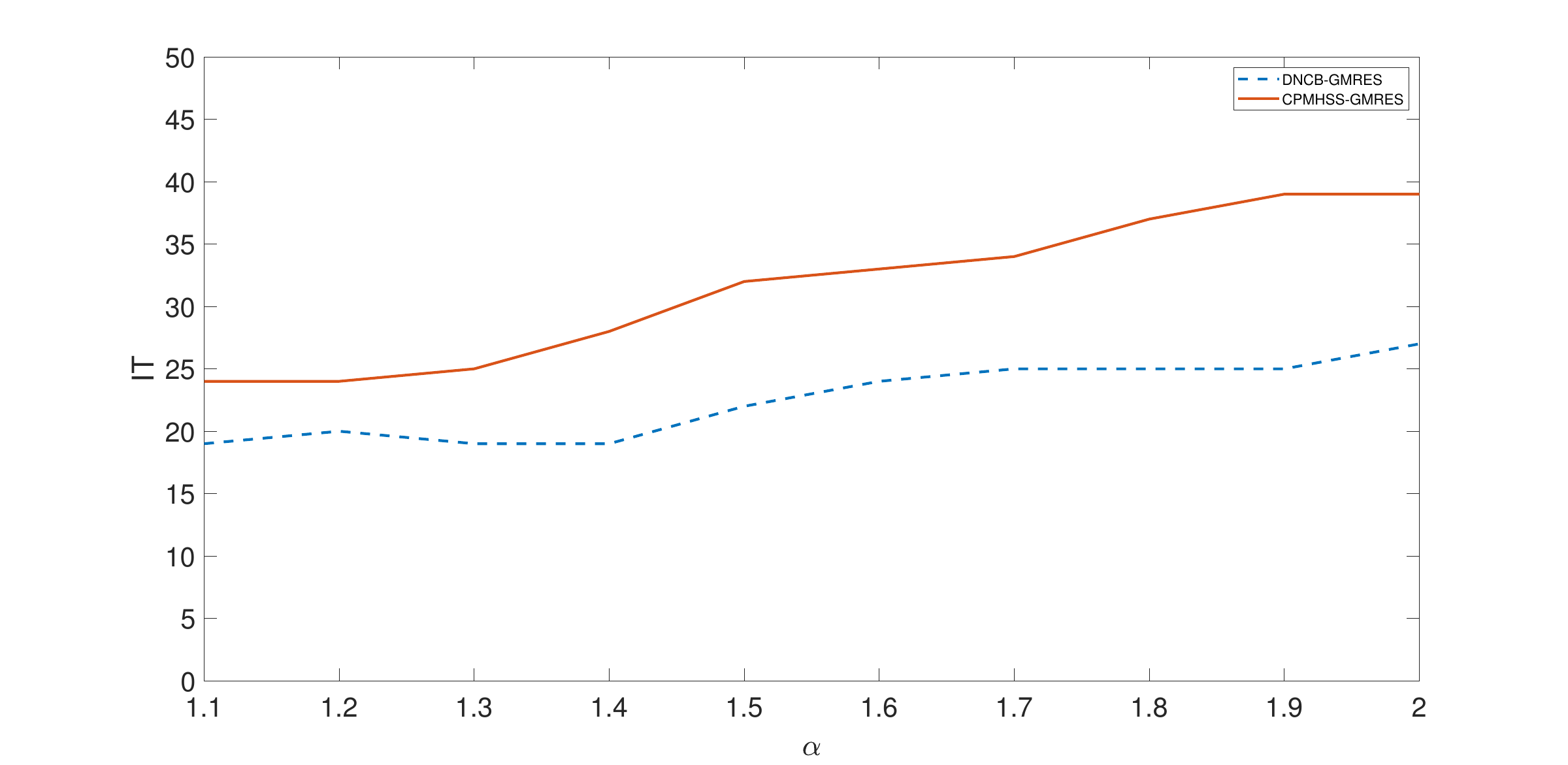}
		\\
	\end{tabular}
	\caption{ The curves of IT of DNCB-GMRES and CPMHSS-GMRES versus the fractional order $\alpha=1.1:0.1:2$ in the CNLS case when $M=6400$: blue dashed line for DNCB-GMRES, red solid line for CPMHSS-GMRES.}
	\label{fig:couple-alp-IT}
\end{figure}

For the case of $M=6400$, the curves of IT of DNCB-GMRES and CPMHSS-GMRES versus the fractional order ($\alpha=1.1:0.1:2$) are plotted in Figure \ref{fig:couple-alp-IT}. The blue dashed line represents DNCB-GMRES, and the red solid line represents CPMHSS-GMRES. As shown in Figure \ref{fig:couple-alp-IT}, IT grows approximately linearly as the fractional order $\alpha$ increases, and the slopes of both curves for DNCB-GMRES and CPMHSS-GMRES stay very low. It indicates that the convergence rates of DNCB-GMRES and CPMHSS-GMRES both have less dependence on or sensitivity to the fractional order $\alpha$. In addition, the DNCB-GMRES curve stays below the CPMHSS-GMRES curve, and the gap between them becomes larger as $\alpha$ tends to $2$, which means that DNCB-GMRES outperforms CPMHSS-GMRES in terms of IT, and CBDN-GMRES is less sensitive to $\alpha$ than CPMHSS-GMRES.

Tables \ref{tab:alp=1.1-beta=1}-\ref{tab:alp=1.9-beta=1} list CPU of GE, and IT and CPU of GMRES, DNCB-GMRES and CPMHSS-GMRES. Here, we take fractional orders $\alpha=1.1:0.2:1.9$, and the number of the spatial discrete points $M=3200$, $6400$, $12800$ and $25600$. In these tables, ``N/A'' means that IT is not available for GE, the notation ``--'' means that GE fails to find the solution or the related iteration method reaches the maximum number of iterations without convergence in that case. The related empirical optimal parameters of DNCB-GMRES and CPMHSS-GMRES are listed in Tables \ref{tab:omega-alp=1.1}-\ref{tab:omega-alp=1.9}. Specifically, the empirical optimal parameters are denoted by ``$\omega_u$'' and ``$\omega_v$'' when the preconditioned GMRES methods are used to solve the block linear systems (\ref{positiveBlockForm}) related to $u$ and $v$ respectively. According to the data in Tables \ref{tab:alp=1.1-beta=1}-\ref{tab:alp=1.9-beta=1}, we can further compute and list the speed-up (SU) of DNCB-GMRES against CPMHSS-GMRES in Table \ref{tab:sp-3200} in which SU is defined as follows
\begin{align} \nonumber
\mbox{SU} &=
\frac{\mbox{CPU of CPMHSS-GMRES}}{\mbox{CPU of DNCB-GMRES}},
\end{align}
and the ratio of IT of CPMHSS-GMRES to that of DNCB-GMRES is shown in the parentheses next to SU.
According to Tables \ref{tab:alp=1.1-beta=1}-\ref{tab:alp=1.9-beta=1}, when the number of the spatial discrete points increases to $M=12800$ and above, GE encounters the problem of ``Out of Memory'', and it fails to solve the linear system. For the case of $M$ smaller than 12800, GE is much more time consuming than all the tested iteration methods. This is because the storage requirements and computational costs of GE cannot be effectively reduced when the linear system is dense even though it has Toeplitz-plus-diagonal structure. Hence, GE can only handle the situation of small dense linear systems for the case of limited storage and computing resources. According to the results of the iteration methods, GMRES fails to converge in the case of $\alpha=1.9$ for $M=25600$, and it requires the largest IT and CPU in all the other tests. Moreover, when the fractional order $\alpha$ and the number of spatial discrete points $M$ grows, IT of GMRES increases very rapidly, which means that the discrete space fractional CNLS equations become more difficult to solve for larger value of $\alpha$ and $M$. Meanwhile, DNCB-GMRES and CPMHSS-GMRES converge successfully in all the tests, and the former outperforms the latter in terms of both IT and CPU. To further demonstrate the advantages of DNCB-GMRES, it is reasonable to check the speed-up of DNCB-GMRES against CPMHSS-GMRES listed in Table \ref{tab:sp-3200}, which ranges from $110.43\%$ to $150.74\%$, roughly in line with the ratio of IT of CPMHSS-GMRES to that of DNCB-GMRES ranging from 1.17 to 1.53. Therefore, the computational efficiency of DNCB-GMRES is further improved compared to that of CPMHSS-GMRES.

\begin{figure}[htbp]
	\centering
	\begin{tabular}{cc}
		\includegraphics[scale=0.15]{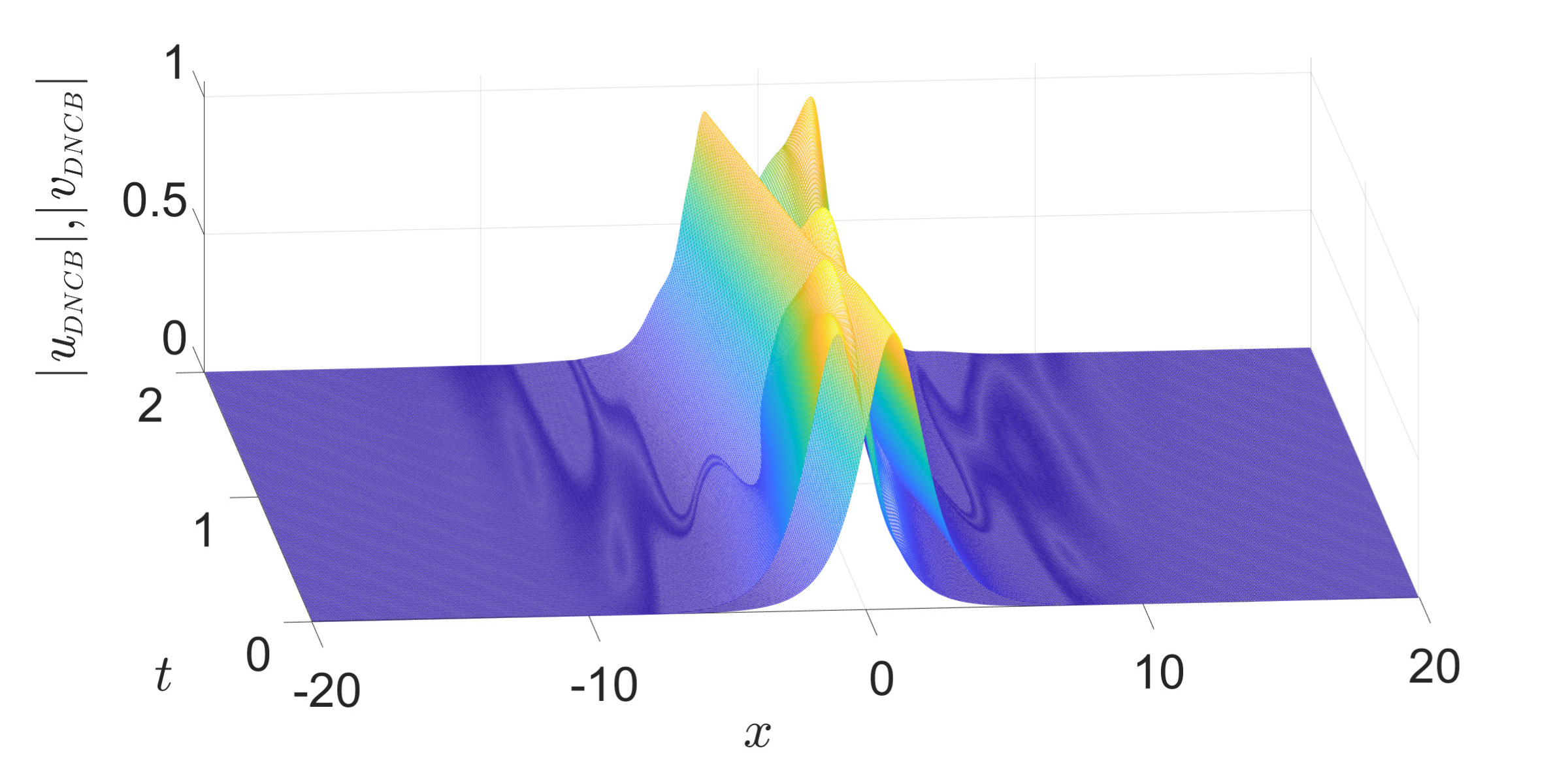}  &
		\includegraphics[scale=0.15]{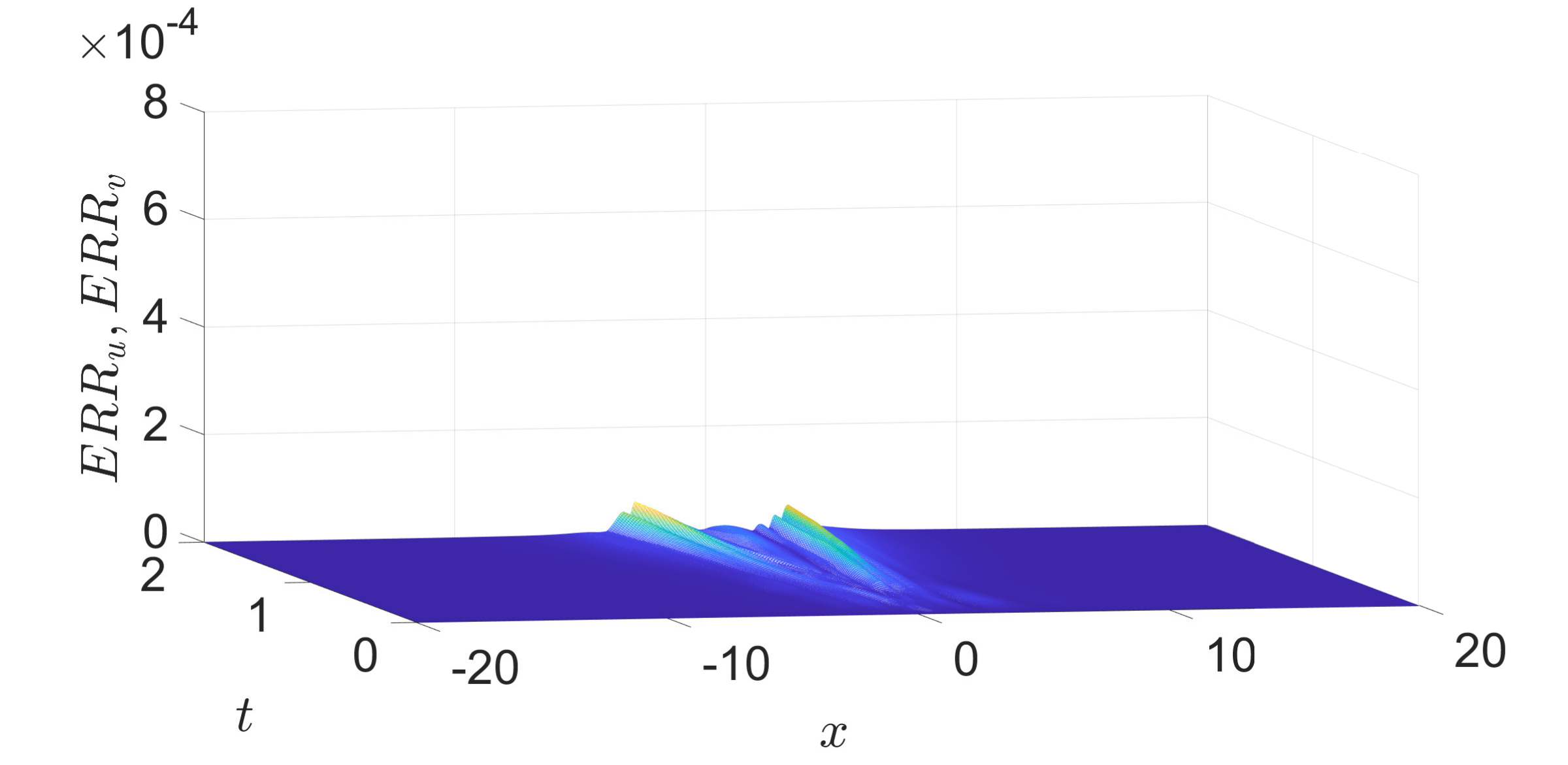}
		
		\\
	\end{tabular}
	\caption{The numerical solution (left) and its error (right) with the exact solution of the LICD scheme in the CNLS case when $\alpha=1.1$ and $M=800$.}
	\label{fig:couple-alp=1.1-beta=1}
\end{figure}

\begin{figure}[htbp]
	\centering
	\begin{tabular}{cc}
		\includegraphics[scale=0.15]{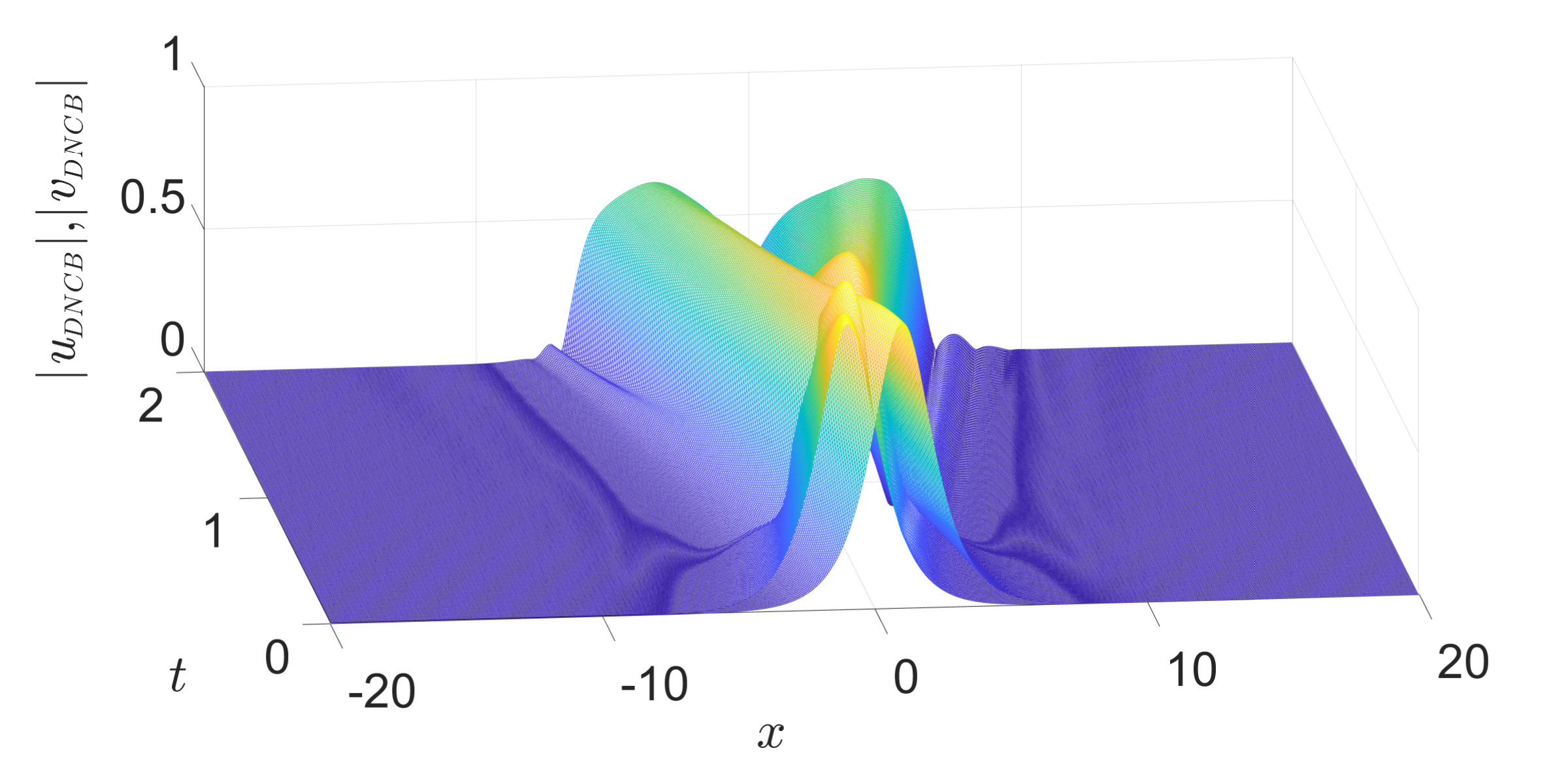}  &
		\includegraphics[scale=0.15]{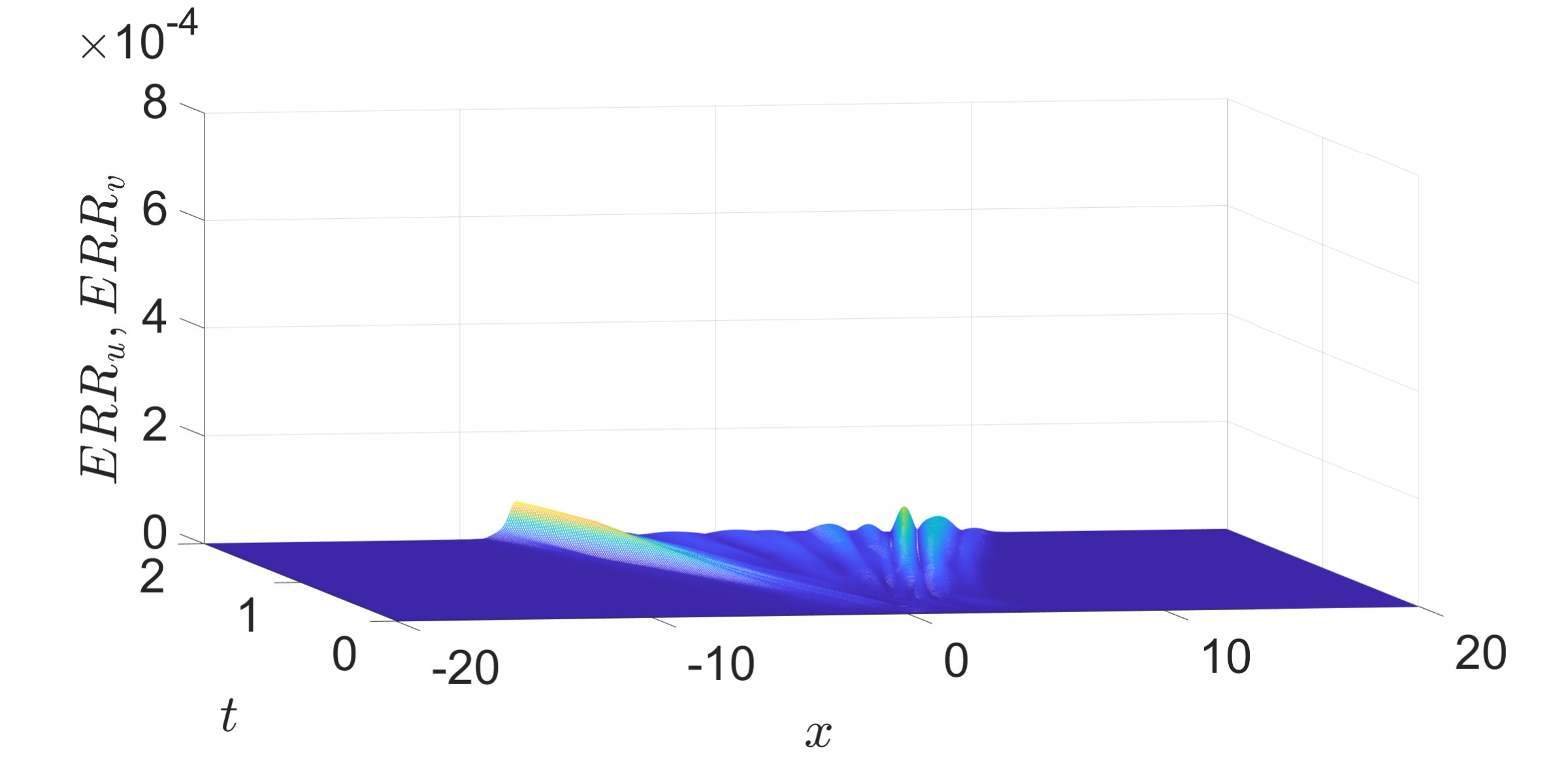}
		
		\\
	\end{tabular}
	\caption{The numerical solution (left) and its error (right) with the exact solution of the LICD scheme in the CNLS case when $\alpha=1.5$ and $M=800$.}
	\label{fig:couple-alp=1.5-beta=1}
\end{figure}

\begin{figure}[htbp]
	\centering
	\begin{tabular}{cc}
		\includegraphics[scale=0.15]{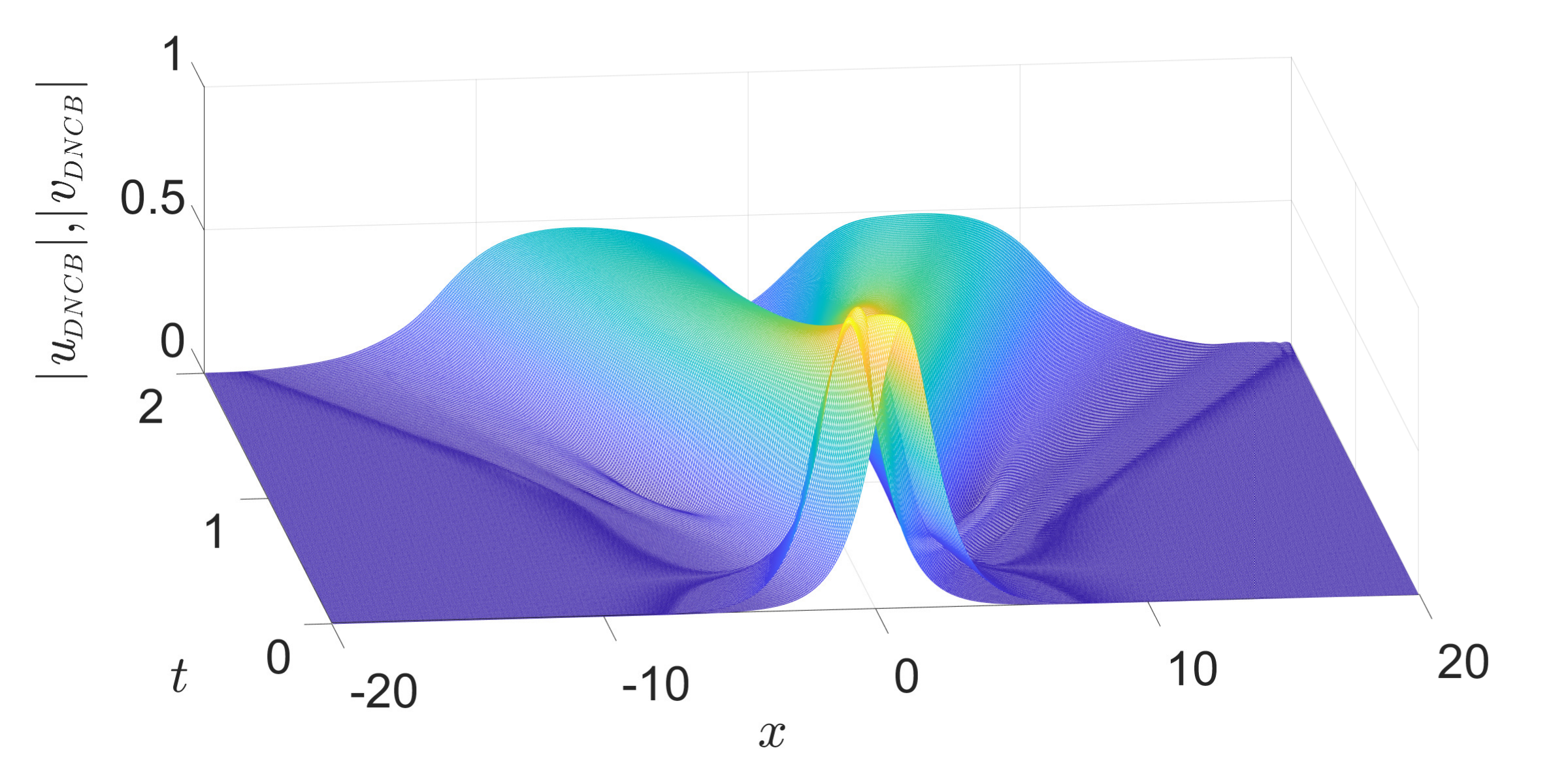}   &
		\includegraphics[scale=0.15]{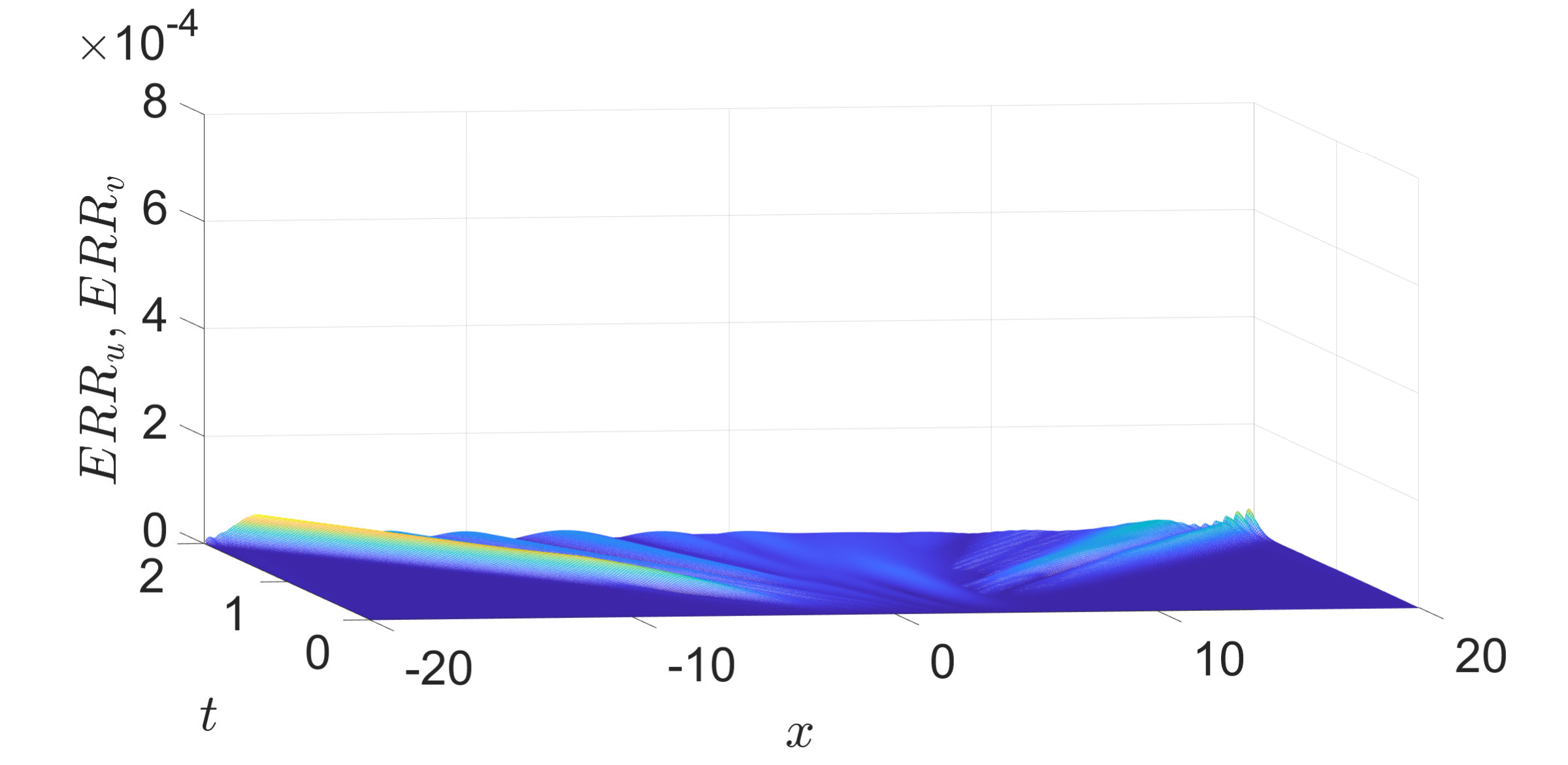}
		\\
	\end{tabular}
	\caption{The numerical solution (left) and its error (right) with the exact solution of the LICD scheme in the CNLS case when $\alpha=1.9$ and $M=800$.}
	\label{fig:couple-alp=1.9-beta=1}
\end{figure}

\begin{figure}[htbp]
	\centering
	\begin{tabular}{cc}
		\includegraphics[scale=0.15]{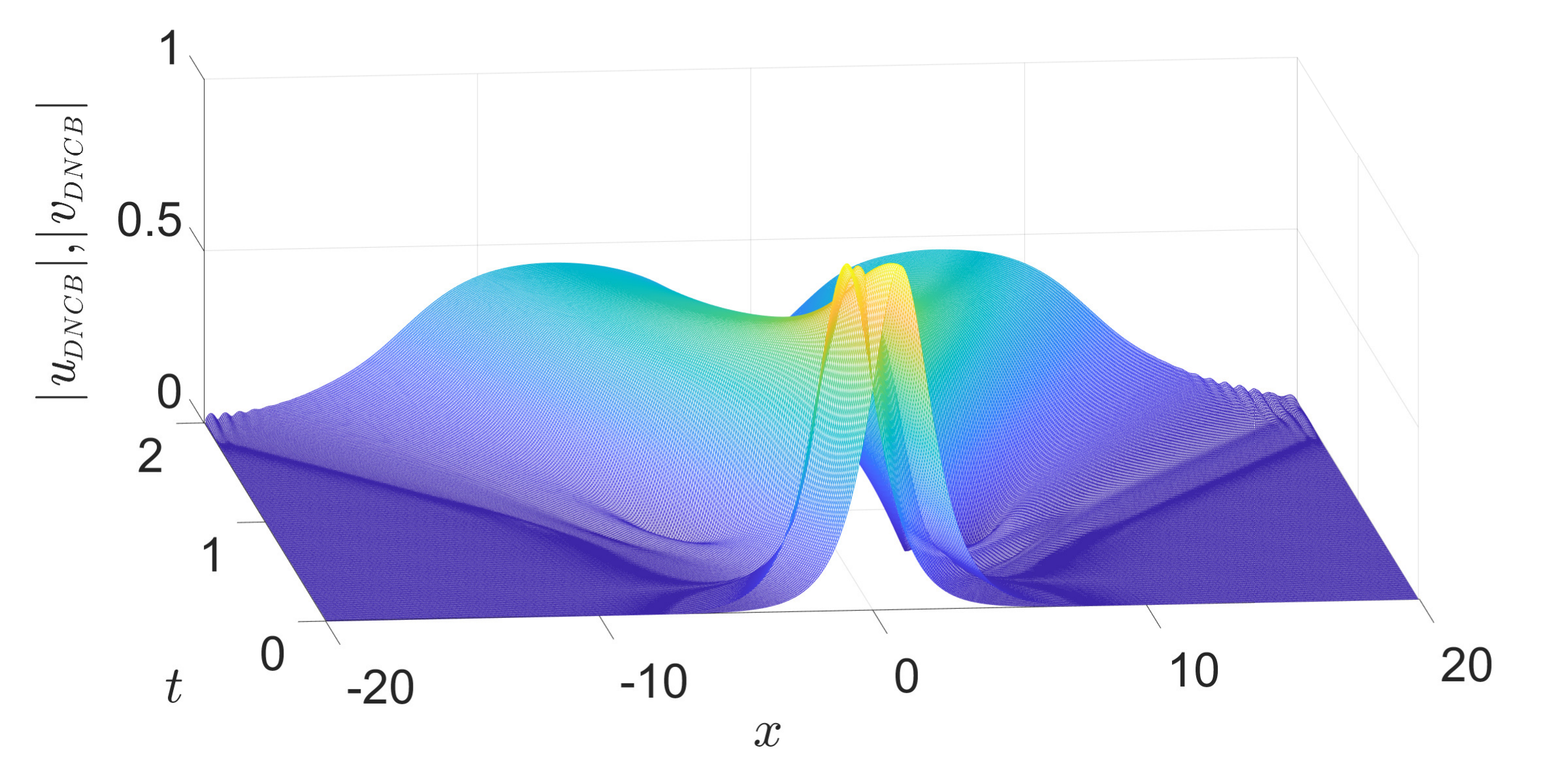}   &
		\includegraphics[scale=0.15]{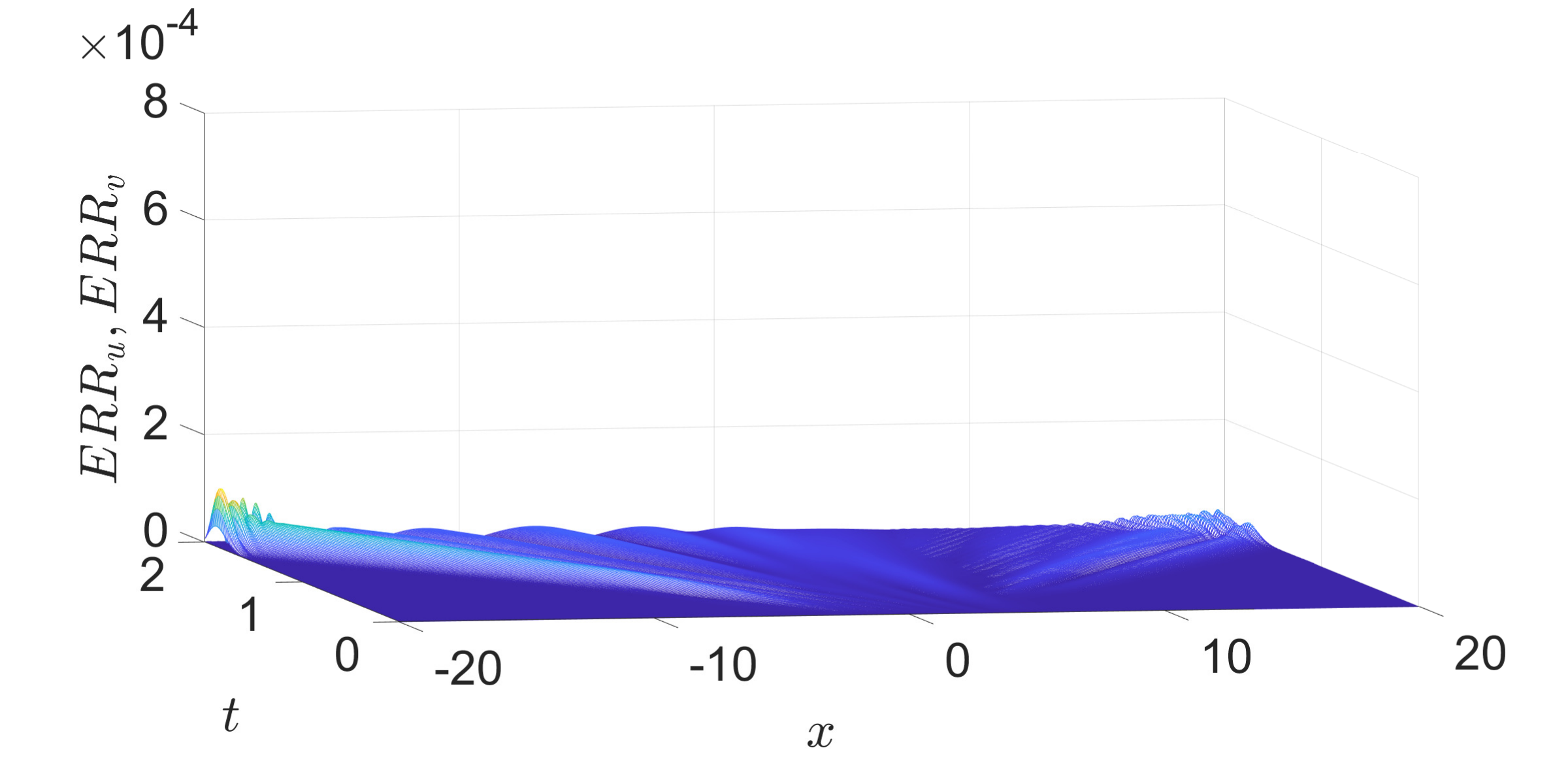}
		\\
	\end{tabular}
	\caption{The numerical solution (left) and its error (right) with the exact solution of the LICD scheme in the CNLS case when $\alpha=2$ and $M=800$.}
	\label{fig:couple-alp=2-beta=1}
\end{figure}

In the case of $M=800$, and $\alpha=1.1:0.4:1.9$ and $\alpha=2$, Figures \ref{fig:couple-alp=1.1-beta=1}-\ref{fig:couple-alp=2-beta=1} show the plots of the numerical solutions $u_{\text{\tiny DNCB}}$ and $v_{\text{\tiny DNCB}}$ of the space fractional CNLS equations (\ref{couple equation}) on the left, and their errors $\text{ERR}_{u} = |u_{\text{\tiny DNCB}}-u_{\text{\tiny GE}}|$ and $\text{ERR}_{v} = |v_{\text{\tiny DNCB}}-v_{\text{\tiny GE}}|$ on the right. Here, $u_{\text{\tiny DNCB}}$ and $v_{\text{\tiny DNCB}}$ are obtained by solving the discrete space fractional CNLS equations with DNCB-GMRES on each time level of the LICD scheme, and $u_{\text{\tiny GE}}$ and $v_{\text{\tiny GE}}$ are the exact solutions of the original LICD scheme computed by GE. It is observed that the shape of the numerical solution is affected by the fractional order $\alpha$. When $\alpha$ approaches to $2$, the shape of the solution of the space fractional CNLS equations (\ref{couple equation}) tends to converge to that of the solution of the standard CNLS equations ($\alpha=2$). In addition, the error between the numerical solution of (\ref{couple equation}) and the exact solution of the original LICD scheme remains as small as around $\mathcal{O}(10^{-4})$ in the whole space-time domain, which shows that the numerical solution obtained by DNCB-GMRES is reliable. In fact, numerical solutions $u_{\text{\tiny DNCB}}$ and $v_{\text{\tiny DNCB}}$ with higher precision can be obtained by improving the stopping criterion of DNCB-GMRES.

\begin{table}[htbp]
	\setlength{\abovecaptionskip}{0pt}
	\setlength{\belowcaptionskip}{10pt} \centering{
		\caption{\label{tab:alp=1.1-beta=1}
			CPU of GE, and IT and CPU of GMRES, DNCB-GMRES and CPMHSS-GMRES in the CNLS case when $\alpha=1.1$ and $M=3200$, $6400$, $12800$, $25600$.}\scriptsize
\begin{tabular}{llrcccccc}\specialrule{0em}{2pt}{2pt}\hline\specialrule{0em}{2pt}{2pt}
			$M$  &  & 3200 &   & 6400 &     & 12800 &     & 25600     \\\cmidrule(l){2-3}\cmidrule(l){4-5}\cmidrule(l){6-7}\cmidrule(l){8-9}
			&  IT &  CPU & IT &  CPU  & IT &  CPU & IT &  CPU   \\\specialrule{0em}{1pt}{1pt}\hline\specialrule{0em}{3pt}{3pt}			
			GE &  N/A &	3.66E+01  &  N/A	 &	2.87E+02 &  N/A	 &  --	  &  N/A		&  --	
			\\\specialrule{0em}{3pt}{3pt}				
			GMRES &  	37 &	2.82E+00&	51 &	1.03E+00 & 68 & 5.36E+00 & 82	& 1.62E+01
			\\\specialrule{0em}{3pt}{3pt}			
			DNCB-GMRES &  	19 &	7.01E-02&	19 &	1.30E-01 & 19 & 3.53E-01 & 19 &	7.67E-01			
			\\\specialrule{0em}{3pt}{3pt}				
			CPMHSS-GMRES &  	24 &	8.47E-02&	24 &	1.53E-01 &  24   &4.22E-01 &24 &	1.01E+00			
			\\\specialrule{0em}{3pt}{3pt}\hline			
	\end{tabular}
}
\end{table}

\begin{table}[htbp]
	\setlength{\abovecaptionskip}{0pt}
	\setlength{\belowcaptionskip}{10pt} \centering{
		\caption{\label{tab:alp=1.3-beta=1}
			CPU of GE, and IT and CPU of GMRES, DNCB-GMRES and CPMHSS-GMRES in the CNLS case when $\alpha=1.3$ and $M=3200$, $6400$, $12800$, $25600$.}\scriptsize
\begin{tabular}{llrcccccc}\specialrule{0em}{2pt}{2pt}\hline\specialrule{0em}{2pt}{2pt}
			$M$  &  & 3200 &   & 6400 &     & 12800 &     & 25600     \\\cmidrule(l){2-3}\cmidrule(l){4-5}\cmidrule(l){6-7}\cmidrule(l){8-9}
			&    IT &  CPU & IT &  CPU  & IT &  CPU & IT &  CPU   \\\specialrule{0em}{1pt}{1pt}\hline\specialrule{0em}{3pt}{3pt}			
			GE &  N/A &	3.77E+01  &  N/A	 &	2.90E+02 &  N/A	 &  --	  &  N/A		&  --	
			\\\specialrule{0em}{3pt}{3pt}			
			GMRES &  	62 &	4.84E-01&	91 &	1.18E+00 & 137 & 4.60E+00 & 196	& 2.16E+01
			\\\specialrule{0em}{3pt}{3pt}			
			DNCB-GMRES & 	19 &	6.86E-02&	19 &	1.24E-01 & 21 & 3.83E-01 & 24 &	9.35E-01			
			\\\specialrule{0em}{3pt}{3pt}				
			CPMHSS-GMRES &  	24 &	8.53E-02&	26 &	1.68E-01 &  29   &5.20E-01 &28 &	1.03E+00			
			\\\specialrule{0em}{3pt}{3pt}\hline			
	\end{tabular}
}
\end{table}

\begin{table}[htbp]
	\setlength{\abovecaptionskip}{0pt}
	\setlength{\belowcaptionskip}{10pt} \centering{
		\caption{\label{tab:alp=1.5-beta=1}
			CPU of GE, and IT and CPU of GMRES, DNCB-GMRES and CPMHSS-GMRES in the CNLS case when $\alpha=1.5$ and $M=3200$, $6400$, $12800$, $25600$.}\scriptsize
		\begin{tabular}{llrcccccc}\specialrule{0em}{2pt}{2pt}\hline\specialrule{0em}{2pt}{2pt}
			$M$  &  & 3200 &   & 6400 &     & 12800 &     & 25600     \\\cmidrule(l){2-3}\cmidrule(l){4-5}\cmidrule(l){6-7}\cmidrule(l){8-9}
			&    IT &  CPU & IT &  CPU  & IT &  CPU & IT &  CPU   \\\specialrule{0em}{1pt}{1pt}\hline\specialrule{0em}{3pt}{3pt}			
			GE &  N/A &	3.61E+01  &  N/A	 &	2.89E+02 &  N/A	 &  --	  &  N/A		&  --	
			\\\specialrule{0em}{3pt}{3pt}			
			GMRES &  	111 &	1.10E+00&	182 &	3.41E+00 & 298 & 1.69E+01 & 442	& 7.02E+01
			\\\specialrule{0em}{3pt}{3pt}			
			DNCB-GMRES & 	19 &	6.78E-02&	22 &	1.52E-01 & 24 & 4.65E-01 & 25 &	9.84E-01			
			\\\specialrule{0em}{3pt}{3pt}				
			CPMHSS-GMRES &  	29 &	1.02E-01&	32 &	2.16E-01 &  32   &5.71E-01 & 31 &	1.16E+00			
			\\\specialrule{0em}{3pt}{3pt}\hline			
	\end{tabular}
}
\end{table}

\begin{table}[htbp]
	\setlength{\abovecaptionskip}{0pt}
	\setlength{\belowcaptionskip}{10pt} \centering{
		\caption{\label{tab:alp=1.7-beta=1}
			CPU of GE, and IT and CPU of GMRES, DNCB-GMRES and CPMHSS-GMRES in the CNLS case when $\alpha=1.7$ and $M=3200$, $6400$, $12800$, $25600$.}\scriptsize
\begin{tabular}{llrcccccc}\specialrule{0em}{2pt}{2pt}\hline\specialrule{0em}{2pt}{2pt}
			$M$  &  & 3200 &   & 6400 &     & 12800 &     & 25600     \\\cmidrule(l){2-3}\cmidrule(l){4-5}\cmidrule(l){6-7}\cmidrule(l){8-9}
			&    IT &  CPU & IT &  CPU  & IT &  CPU & IT &  CPU   \\\specialrule{0em}{1pt}{1pt}\hline\specialrule{0em}{3pt}{3pt}			
			GE &  N/A &	3.76E+01  &  N/A	 &	2.97E+02 &  N/A	 &  --	  &  N/A		&  --	
			\\\specialrule{0em}{3pt}{3pt}			
			GMRES &  	204 &	2.24E+00&	364 &	1.11E+01 & 647 & 6.05E+01 & 1126	& 3.54E+02
			\\\specialrule{0em}{3pt}{3pt}			
			DNCB-GMRES & 	24 &	9.48E-02&	25 &	1.79E-01 & 25 & 4.78E-01 & 26 &	1.05E+00			
			\\\specialrule{0em}{3pt}{3pt}				
			CPMHSS-GMRES & 	33 &	1.25E-01&	34 &	2.41E-01 &  34   &6.20E-01 & 33 &	1.26E+00			
			\\\specialrule{0em}{3pt}{3pt}\hline			
	\end{tabular}
}
\end{table}

\begin{table}[htbp]
	\setlength{\abovecaptionskip}{0pt}
	\setlength{\belowcaptionskip}{10pt} \centering{
		\caption{\label{tab:alp=1.9-beta=1}
			CPU of GE, and IT and CPU of GMRES, DNCB-GMRES and CPMHSS-GMRES in the CNLS case when $\alpha=1.9$ and $M=3200$, $6400$, $12800$, $25600$.}\scriptsize
\begin{tabular}{llrcccccc}\specialrule{0em}{2pt}{2pt}\hline\specialrule{0em}{2pt}{2pt}
			$M$  &  & 3200 &   & 6400 &     & 12800 &     & 25600     \\\cmidrule(l){2-3}\cmidrule(l){4-5}\cmidrule(l){6-7}\cmidrule(l){8-9}
			&    IT &  CPU & IT &  CPU  & IT &  CPU & IT &  CPU    \\\specialrule{0em}{1pt}{1pt}\hline\specialrule{0em}{3pt}{3pt}			
			GE &  N/A &	3.62E+01  &  N/A	 &	2.95E+02 &  N/A	 &  --	  &  N/A		&  --	
			\\\specialrule{0em}{3pt}{3pt}			
			GMRES & 	414 &	7.45E+00&	786 &	4.41E+01 & 1497 & 2.92E+02 & --	& 8.74E+02
			\\\specialrule{0em}{3pt}{3pt}			
			DNCB-GMRES &  	25 &	9.75E-02&	25 &	1.79E-01  &  27  & 5.24E-01 &28&	1.15E+00			
			\\\specialrule{0em}{3pt}{3pt}				
			CPMHSS-GMRES &  	36&	1.39E-01&	37&	2.51E-01  &  37  & 6.70E-01 &37 &	1.39E+00			
			\\\specialrule{0em}{3pt}{3pt}\hline			
	\end{tabular}
}
\end{table}

\begin{table}[htbp]
	\setlength{\abovecaptionskip}{0pt}
	\setlength{\belowcaptionskip}{10pt} \centering{
		\caption{\label{tab:sp-3200}
			The speed-up of DNCB-GMRES against CPMHSS-GMRES in the CNLS case when $\alpha=1.1:0.2:1.9$ and $M=3200$, $6400$, $12800$, $25600$.}\scriptsize
\begin{tabular}{llrcccc}\specialrule{0em}{2pt}{2pt}\hline\specialrule{0em}{2pt}{2pt}
			$M$  &  & $\alpha$=1.1 &   $\alpha$=1.3 &     $\alpha$=1.5 &     $\alpha$=1.7 &  $\alpha$=1.9
			\\\specialrule{0em}{1pt}{1pt}\hline\specialrule{0em}{3pt}{3pt}			
			3200 &	SU &  120.83$\%$ (1.26)  &	124.34$\%$ (1.26)  &	150.74$\%$ (1.53)  &	132.07$\%$ (1.42)  &	142.97$\%$ (1.44)
			\\\specialrule{0em}{3pt}{3pt}			
			6400	& SU &  117.87$\%$ (1.26)  &	136.14$\%$ (1.37)  &	142.29$\%$ (1.45)  &	134.56$\%$ (1.36)  &	140.12$\%$ (1.48)
			\\\specialrule{0em}{3pt}{3pt}			
			12800 & SU &  119.54$\%$ (1.26)  &	135.59$\%$ (1.38) & 	122.91$\%$ (1.33) &	129.56$\%$ (1.36) &	127.79$\%$ (1.37)				
			\\\specialrule{0em}{3pt}{3pt}			
			25600 & SU & 131.22$\%$ (1.26)	& 110.43$\%$ (1.17) &	118.29$\%$ (1.24) &	120.27$\%$ (1.27) &	121.76$\%$ (1.32)
			\\\specialrule{0em}{3pt}{3pt}\hline			
	\end{tabular}
}

\end{table}

\begin{table}[htbp]
	\setlength{\abovecaptionskip}{0pt}
	\setlength{\belowcaptionskip}{10pt} \centering{
		\caption{\label{tab:omega-alp=1.1}
			The empirical optimal parameters of DNCB-GMRES and CPMHSS-GMRES in the CNLS case when $\alpha=1.1$ and $M=3200$, $6400$, $12800$, $25600$.}\scriptsize
		\begin{tabular}{llrccc}\specialrule{0em}{2pt}{2pt}\hline\specialrule{0em}{2pt}{2pt}
			
			$M$ & &  3200 &  6400 &12800  & 25600   \\\specialrule{0em}{1pt}{1pt}\hline\specialrule{0em}{3pt}{3pt}
			
			\multirow{2}{*}{DNCB-GMRES}  & $\omega_u$  &[0.12,0.52] &[0.12,0.52] &[0.12,0.52] & [0.12,0.56]	
			\\\cmidrule(l){2-6}
			& $\omega_v$ &[0.12,1.52]  &[0.12,1.60] &[0.12,0.20]
			&[0.12,0.20]	
			\\\specialrule{0em}{1pt}{1pt}\hline\specialrule{0em}{3pt}{3pt}
			
			\multirow{2}{*}{CPMHSS-GMRES}  & $\omega_u$   &[0.80,2.00]  &[0.80,2.00] &[0.70,2.00] &[1.15,1.45] 	
			\\\cmidrule(l){2-6}
			& $\omega_v$  &[0.80,2.00]  &[0.80,2.00]  &[0.75,1.95] & [0.70,1.25]	
			\\\specialrule{0em}{1pt}{1pt}\hline
	\end{tabular}}
\end{table}

\begin{table}[htbp]
	\setlength{\abovecaptionskip}{0pt}
	\setlength{\belowcaptionskip}{10pt} \centering{
		\caption{\label{tab:omega-alp=1.3}
			The empirical optimal parameters of DNCB-GMRES and CPMHSS-GMRES in the CNLS case when $\alpha=1.3$ and $M=3200$, $6400$, $12800$, $25600$.}\scriptsize
		\begin{tabular}{llrccc}\specialrule{0em}{2pt}{2pt}\hline\specialrule{0em}{2pt}{2pt}
			
			$M$  & &  3200 &  6400 &12800  & 25600   \\\specialrule{0em}{1pt}{1pt}\hline\specialrule{0em}{3pt}{3pt}
			
			\multirow{2}{*}{DNCB-GMRES}  & $\omega_u$  &[0.15,0.55] &[0.15,0.50] &[0.15,0.45] & [0.15,0.55]	
			\\\cmidrule(l){2-6}
			& $\omega_v$ &[1.10,1.35]  &[1.00,1.05] &[0.10,0.40]
			&[0.10,0.25]	
			\\\specialrule{0em}{1pt}{1pt}\hline\specialrule{0em}{3pt}{3pt}
			
			\multirow{2}{*}{CPMHSS-GMRES}  & $\omega_u$   &[1.55,1.95]  &[1.05,1.60] &[0.90,0.90] &[1.55,2.00] 	
			\\\cmidrule(l){2-6}
			& $\omega_v$  &[1.20,1.65]  &[1.10,1.60]  &[1.40,1.40] & [0.80,1.95]	
			\\\specialrule{0em}{1pt}{1pt}\hline
	\end{tabular}}
\end{table}

\begin{table}[htbp]
	\setlength{\abovecaptionskip}{0pt}
	\setlength{\belowcaptionskip}{10pt} \centering{
		\caption{\label{tab:omega-alp=1.5}
			The empirical optimal parameters of DNCB-GMRES and CPMHSS-GMRES in the CNLS case when $\alpha=1.5$ and $M=3200$, $6400$, $12800$, $25600$.}\scriptsize
		\begin{tabular}{llrccc}\specialrule{0em}{2pt}{2pt}\hline\specialrule{0em}{2pt}{2pt}
			
			$M$  & &  3200 &  6400 &12800  & 25600   \\\specialrule{0em}{1pt}{1pt}\hline\specialrule{0em}{3pt}{3pt}
			
			\multirow{2}{*}{DNCB-GMRES}  & $\omega_u$  &[0.15,0.45] &[0.15,0.40] &[0.15,0.45] & [0.15,0.40]	
			\\\cmidrule(l){2-6}
			& $\omega_v$ &[0.10,0.35]  &[0.10,0.15] &[0.10,0.15]
			&[1.05,1.30]	
			\\\specialrule{0em}{1pt}{1pt}\hline\specialrule{0em}{3pt}{3pt}
			
			\multirow{2}{*}{CPMHSS-GMRES}  & $\omega_u$   &[0.80,1.05]  &[1.00,2.00] &[1.05,2.00] &[1.65,2.00] 	
			\\\cmidrule(l){2-6}
			& $\omega_v$  &[0.80,1.80]  &[0.75,1.60]  &[0.65,2.00] & [1.10,2.00]	
			\\\specialrule{0em}{1pt}{1pt}\hline
			
	\end{tabular}}
\end{table}

\begin{table}[htbp]
	\setlength{\abovecaptionskip}{0pt}
	\setlength{\belowcaptionskip}{10pt} \centering{
		\caption{\label{tab:omega-alp=1.7}
			The empirical optimal parameters of DNCB-GMRES and CPMHSS-GMRES in the CNLS case when $\alpha=1.7$ and $M=3200$, $6400$, $12800$, $25600$.}\scriptsize
		\begin{tabular}{llrccc}\specialrule{0em}{2pt}{2pt}\hline\specialrule{0em}{2pt}{2pt}
			
			$M$  & &  3200 &  6400 &12800  & 25600  \\\specialrule{0em}{1pt}{1pt}\hline\specialrule{0em}{3pt}{3pt}
			
			\multirow{2}{*}{DNCB-GMRES}  & $\omega_u$  &[0.10,0.50] &[0.15,0.35] &[0.10,0.40] & [0.15,0.35]	
			\\\cmidrule(l){2-6}
			& $\omega_v$ &[0.10,0.15]  &[0.10,0.25] &[0.05,0.30]
			&[0.10,0.25]	
			\\\specialrule{0em}{1pt}{1pt}\hline\specialrule{0em}{3pt}{3pt}
			
			\multirow{2}{*}{CPMHSS-GMRES}  & $\omega_u$   &[1.65,2.00]  &[1.10,2.00] &[1.05,2.00] &[1.45,2.00] 	
			\\\cmidrule(l){2-6}
			& $\omega_v$  &[1.00,2.00]  &[0.75,1.95]  &[1.00,1.95] & [0.80,1.90]	
			\\\specialrule{0em}{1pt}{1pt}\hline
			
	\end{tabular}}
\end{table}

\begin{table}[htbp]
	\setlength{\abovecaptionskip}{0pt}
	\setlength{\belowcaptionskip}{10pt} \centering{
		\caption{\label{tab:omega-alp=1.9}
			The empirical optimal parameters of DNCB-GMRES and CPMHSS-GMRES in the CNLS case when $\alpha=1.9$ and $M=3200$, $6400$, $12800$, $25600$.}\scriptsize
		\begin{tabular}{llrccc}\specialrule{0em}{2pt}{2pt}\hline\specialrule{0em}{2pt}{2pt}
			
			$M$  & &  3200 &  6400 &12800  & 25600  \\\specialrule{0em}{1pt}{1pt}\hline\specialrule{0em}{3pt}{3pt}
			
			\multirow{2}{*}{DNCB-GMRES}  & $\omega_u$  &[0.10,0.40] &[0.10,0.30] &[0.15,0.35] & [0.15,0.35]	
			\\\cmidrule(l){2-6}
			& $\omega_v$ &[0.05,0.25]  &[0.10,0.20] &[0.10,0.50]
			&[0.15,0.30]	
			\\\specialrule{0em}{1pt}{1pt}\hline\specialrule{0em}{3pt}{3pt}
			
			\multirow{2}{*}{CPMHSS-GMRES}  & $\omega_u$   &[1.80,2.00]  &[1.10,2.00] &[1.30,1.90] &[1.35,1.75] 	
			\\\cmidrule(l){2-6}
			& $\omega_v$  &[0.85,1.85]  &[0.75,1.90]  &[0.85,1.75] & [0.60,1.65]	
			\\\specialrule{0em}{1pt}{1pt}\hline
			
	\end{tabular}}
\end{table}

\section{Concluding remarks}
\label{sec-conclusions}

The DNCB preconditioned Krylov subspace iteration methods (e.g., DNCB-GMRES) are effective and efficient linear solvers for the complex linear system arising from the system of repulsive space fractional CNLS equations. The DNCB preconditioner is a circulant approximation of the DNTB preconditioner naturally derived from the DNTB iteration method. In theory, the DNTB iteration method converges unconditionally, and the optimal iteration parameter is deducted. Moreover, good clustering properties of the spectra of the DNCB preconditioned system matrix are proved based on the sharp bounds for the eigenvalues of the discrete fractional Laplacian and its circulant approximation. The above theoretical results are sufficiently verified by the numerical experiments based on 1D repulsive space fractional DNLS and CNLS equations. Now, we consider possible future works. Firstly, the initial experiments show that the linear dependance on space mesh size of DNCB-GMRES weakens as the fractional order parameter decreases. This phenomenon can be further analyzed both theoretically and experimentally. Secondly, since the DNCB preconditioner is obtained by circulant approximation of the discrete fractional Laplacian $T$, we can also approximate $T$ by other techniques, for instance, the optimal circulant preconditioning technique based on axisymmetric extension \cite{LMV1997} or the sine transform based preconditioning technique \cite{HuangLinNgSun2022}, which may lead to even more efficient preconditioner. Thirdly, higher dimensional problems are always interesting, thus extending the DNCB preconditioner to the higher dimensional cases deserves further study. Fourthly, variable coefficient problems and non-uniform spatial discretization schemes will make the derived complex linear system lacking the explicit Toeplitz structure. Therefore, exploiting the possible implicit data-sparse structure and combining other approaches (e.g., hierarchical-matrix approach \cite{BM2008SpringerBook,HKL2016CPAM}) with our methodology to construct new preconditioners will be a useful and challenging work.

\vspace{2em}
\section*{Acknowledgments}
This work was funded by the National Natural Science Foundation (No. 11101213 and No. 12071215), China.


\end{document}